\documentclass[letterpaper,oneside]{amsart}
\pdfoutput=1 %arXiv, ensure pdflatex runs
\usepackage[T1]{fontenc}
\usepackage{color,soul}
\usepackage{amsrefs}
\usepackage[pagewise]{lineno}
\usepackage[dvipsnames]{xcolor}
\usepackage{todonotes}
\setuptodonotes{inline}
\usepackage[utf8]{inputenc} % ------ Pour utiliser des fichiers UTF-8
\usepackage{comment}
\usepackage{mathtools} 
\usepackage{colonequals}
\usepackage{amsmath}
\usepackage{amssymb}
\usepackage{tikz}
\usepackage{mathrsfs}
\usepackage{enumitem}
\usepackage[foot]{amsaddr}
\usepackage{mathtools}

\usepackage{mleftright}
\usepackage{graphicx}
\usepackage{hyperref}

\usepackage[english]{babel}
\usepackage[refpage, notocbasic]{nomencl}
\makenomenclature

\usepackage{tikz}
\usepackage{thmtools}
\usepackage{amssymb}

\usepackage[all]{xy} %para diagramas conmutativos

\numberwithin{equation}{section}

\newcommand{\Id}[1]{\operatorname{Id}_{#1}}

\newcommand{\Z}{\mathbb{Z}}

\newcommand{\N}{\mathbb{N}}

\newcommand{\Q}{\mathbb{Q}}

\newcommand{\pushoutdr}[1][dr]{\save*!/#1+1.2pc/#1:(1,-1)@^{|-}\restore}
\newcommand{\pushoutdl}[1][ur]{\save*!/#1-1.2pc/#1:(-1,1)@^{|-}\restore}
\newcommand{\pushoutur}[1][ur]{\save*!/#1+1.2pc/#1:(1,-1)@^{|-}\restore}
\newcommand{\pushoutul}[1][dr]{\save*!/#1-1.2pc/#1:(-1,1)@^{|-}\restore}

\newcommand{\cat}[1]{\mathbf{#1}}
\newcommand{\ocat}[1]{|\mathbf{#1}|}
\newcommand{\mcat}[3]{[#1,#2]_\mathbf{#3}}
\newcommand{\funt}[1]{\mathsf{#1}}
\newcommand{\mfunt}[3]{\mathsf{#1}:{\cat{#2}}\rightarrow{\cat{#3}}}
\newcommand{\dfunt}[4]{\mathsf{#1}:{\cat{#2}}\leftrightarrows{\cat{#3}}:\mathsf{#4}}

\title{Rational Homotopy in Pseudotopological Spaces}
\author{Jonathan Treviño-Marroquín}
\thanks{This work was supported by the SECIHTI postgraduate studies scholarship number 839062, and the results in this article form a part of the author's PhD thesis project supervised by Antonio Rieser.}

\theoremstyle{plain}
\newtheorem{teorema}{Theorem}[section]
\newtheorem{proposicion}[teorema]{Proposition}
\newtheorem{corolario}[teorema]{Corollary}
\newtheorem{lema}[teorema]{Lemma}

\theoremstyle{definition}
\newtheorem{definicion}[teorema]{Definition}
\newtheorem{definition}[teorema]{Definition}
\newtheorem{observacion}[teorema]{Remark}
\newtheorem{ejemplo}{Example}[teorema]

\newtheorem{innercustomthm}{Theorem}
\newenvironment{customthm}[1]
  {\renewcommand\theinnercustomthm{#1}\innercustomthm}
  {\endinnercustomthm}

\newtheorem{innercustomprop}{Proposition}
\newenvironment{customprop}[1]
  {\renewcommand\theinnercustomprop{#1}\innercustomprop}
  {\endinnercustomprop}

\newtheorem{innercustomcoro}{Corollary}
\newenvironment{customcoro}[1]
  {\renewcommand\theinnercustomcoro{#1}\innercustomcoro}
  {\endinnercustomcoro}

\declaretheoremstyle[
qed=\qedsymbol
]{mystyle}

\begin{document}
%\linenumbers

\begin{abstract}
	Pseudotopological spaces are the Cartesian closed hull of the category of \v{C}ech closure spaces.
	In this paper, we give a direct proof that the model category of the pseudotopological spaces constructed by Rieser is Quillen equivalent to the category of simplicial sets.
	In addition to noting that every pseudotopological space is weak homotopy equivalent to a topological CW complex, we prove that any weak equivalence of pseudotopological spaces can be converted to a weak equivalence of topological spaces.
	Finally, combining these ingredients, we construct rational homotopy for simply connected pseudotopological spaces. In this paper we also prove that the cartesian product of two pseudotopological CW complexes is a CW complex.
\end{abstract}

\maketitle

\tableofcontents

\pagenumbering{arabic}
\section*{Introduction}
\label{cap.introduccion}

	Setting a topological structure on the vertices of a finite set does not provide enough topological or algebraic topological information.
	This lack of information has generally been addressed in discrete homotopy in two approaches:
	developing homotopy in a certain class of graphs with certain product or embedding graphs in categories similar to topological spaces.
	Some examples are given below:
\begin{itemize}
	\item Work in a category which only includes graphs.
	For example, $A$-homotopy \cite{Babson_etal_2006}, a cubical-type homotopy is defined in graphs and uses the inductive product and the discrete interval to generate their homotopies, $\times$-homotopy \cite{Dochtermann_2009, Chih_Scull_2021}, introduced for studying not-necessarily reflexive graphs, and digital homotopy \cite{Lupton_etal_2022a,Lupton_Scoville_2022}, which allows both inductive and cartesian product, but where the graphs are determined by a lattice.
	
	\item Work in a category which contains more possible objects, mostly those categories contain the topological spaces.
	For example \v{C}ech closure spaces \cite{Rieser_2021} and pseudotopological spaces \cite{Rieser_arXiv_2022}, both of them could be used as categories in which to simultaneously develop the homotopy theory of topological spaces, graphs, digraphs, and point clouds (the latter for topological data analysis \cite{Carlsson_2009, Carlsson_et_al_2012}), extending earlier work of Demaria \cite{Demaria_1987} on graphs.
	In addition, semi-coarse spaces \cite{rieser2023semicoarse, Trevino_2025_FundGroupoidSC}, which is a generalization of coarse geometry which admits non-trivial structures on finite graphs.
	All of these examples can be studied with both the inductive product or the categorical product \cite{Bubenik_Milicevic_2021}.
\end{itemize}
	
	In particular, in this paper we work with the category of pseudotopological spaces.
	
	The category of pseudotopological spaces contains another category named the category of \v{C}ech closure spaces, which seems to have good homotopy properties to relate graphs and the geometric realization of its Vietoris-Rips complex.
	These two categories have taken importance in recent works because they have interesting objects to study, capture more information from discrete sets and the homotopy can be studied through powerful tools of category theory as it is the model categories. 

	In this paper we take the first steps to study the pseudotopological homotopy groups removing the torsion, i.e., we study the rational homotopy groups of these spaces.
	This homotopy theory branch requires to define objects in pseudotopological spaces and mimic results from topology. The tools to mimic the necessary theorems for rational homotopy are extensively linked to Rieser's paper \cite{Rieser_arXiv_2022}, where he study the Quillen Model category for the category of pseudotopological spaces, and in Ebel and Kapulkin's paper \cite{kapulkin_edel2023_SyntheticApproach}, where they observe as an example a Quillen equivalence between the Quillen model categories of topological and pseudotopological spaces.

	Supported in that papers, joint with the Quillen equivalence between topological spaces and simplicial sets, we can observe a third Quillen equivalence between pseudotopological spaces and simplicial sets.
	This simple observation leads to construct, for every continuous map $f: X\to Y$ of pseudotopological spaces, a continuous map of topological spaces such that one is weak homotopy equivalent if the other does as well, being to close to convert every pseudotopological weak homotopy question in a topological weak homotopy question. By pushing a little further, we obtain results ranging from the Hurewicz theorem to the bijective map \[ \left\lbrace \begin{array}{c}
\text{rational homotopy}\\ \text{types}	
\end{array}\right\rbrace \xrightarrow{\cong} \left\lbrace \begin{array}{c}
\text{isomorphism classes of}\\ \text{minimal Sullivan algebras}\\ \text{over }\Q	
\end{array} \right\rbrace . \]

	This paper also has a side objective which is developed in the second appendix: Observing that the cartesian product of two pseudotopological CW complexes is also a pseudotopological CW complex.
	 Being in general false in topology, but true in the category of compactly generated Hausdorff spaces (which is cartesian closed), it seemed feasible to achieve for pseudotopological spaces.
	 We show that is in general true, using strongly that the category is cartesian closed.

	Throughout this paper, we denote the category of sets with the functions as morphisms as $\cat{Set}$ and the category of topological spaces with continuous maps as $\cat{Top}$.
	Other important categories used in this document, for example the categories of simplicial sets and pseudotopological spaces, are introduced to the reader either in the content of the chapters or the \autoref{Appx:LimPushoutEmb}.

	The paper is organized in the following way:
	In \autoref{cap.PsTopCW} and \autoref{cap.SSetCat} we provide the translation of the CW complex and geometric realization from topological spaces to pseudotopological spaces.
	We observe that everything works in an analogous fashion, although we have to make subtle modifications.
	\autoref{sec.PsTopQuillenSSet} introduce the definitions of model category and Quillen equivalences in this paper, to provide the framework to verify that we actually have a Quillen equivalence between $\cat{PsTop}$ and $\cat{sSet}$. We finish that section observing that every pseudotopological spaces is weak homotopy equivalent to a pseudotopological CW complex.
	
	We define the rational pseudotopological spaces and the rational homotopy equivalences in \autoref{sec.SingHomHureTheo}, also including the map from $\cat{PsTop}$ morphisms to $\cat{Top}$ morphisms which satisfies that one is weak homotopy equivalent if the other does as well.
	The last section, \autoref{sec.RatHom}, recall the results from differential graded commutative algebra and make the necessary modifications of the proves of the affirmations to obtain the relations between the same rational homotopy type and isomorphic classes of Sullivan algebras.
	
	Lastly, in \autoref{Appx:LimPushoutEmb} and \autoref{Appx:ProductCWComplex}, we give some preliminaries and show that the cartesian product of two CW complexes is a CW complex in the category $\cat{PsTop}$.

\section{Pseudotopological CW Complexes}
\label{cap.PsTopCW}

	CW complexes joint with cell morphisms are really important objects in algebraic topology.
	This subcategory of the topological spaces allows us to study the homotopy and singular homology groups considering friendlier spaces and morphisms.
	In this chapter, we introduce the CW complexes in the category of pseudotopological spaces ($\cat{PsTop}$, \autoref{def:PsTopSpaces}) and we observe that, like in $\cat{Top}$, the collection of pseudotopological CW complexes are closed under retracts.
	
\begin{definicion}[\cite{Rieser_arXiv_2022}, 5.12]
\label{def:CellAttachment}
	If $X$ is a subspace of $Y$ and there is a pushout square
	\[ \xymatrix{
	S^{n-1} \ar[r]  \ar[d] & X \ar[d]\\
	D^n \ar[r] & \pushoutdr Y
	} \]
	for some $n\geq 0$, then we say that $Y$ is obtained from $X$ by \emph{attaching a cell} or a \emph{cell attachment}\index{cell attachment}. As a convention, when $n=0$, we let $S^{-1}\coloneqq \varnothing$ and $D^0\coloneqq\{*\}$, the one point set.
\end{definicion}

\begin{definicion}[\cite{Rieser_arXiv_2022}, 5.13]
\label{def:RelativeCellComplex}
	We say that a continuous map $f:X\rightarrow Y$ between pseudotopological spaces is a \emph{relative cell complex}\index{relative cell complex} if $f$ is an inclusion and $Y$ can be constructed from $X$ by a (possibly infinite, and even transfinite) sequence of cell attachments.
	If $Y$ may be constructed from $X$ by attaching a finite number of cells, then we say that it is a \emph{finite relative cell complex}\index{relative cell complex!finite} and if $\varnothing \rightarrow Y$ is a (finite) relative cell complex, then we say that $Y$ is a \emph{(finite) cell complex}\index{cell complex}.
\end{definicion}

\begin{definicion}
\label{def:PsTopCWComplex}
A \emph{(pseudotopological) CW complex}\index{CW complex} $X$ is a (pseudotopological) cell complex (from \autoref{def:RelativeCellComplex}) which can be written as the following sequence of cell attachments
\begin{align*}
\xymatrix{
	\varnothing \ar[r] \ar[d] \ar@{}[dr]^>{\ulcorner} & \bigsqcup\limits_{i\in I_0} D^0 \ar[d] & \bigsqcup\limits_{i\in I_2} S^2 \ar[r] \ar[d] \ar@{}[dr]^>{\ulcorner} & \bigsqcup\limits_{i\in I_2} D^2 \ar[d] \\
	\varnothing = X_{-1} \ar[r]  & X_0 \ar[r] & X_1 \ar[r] & X_2 \ar@{..>}[r] & X\\
	& \ar[u] \bigsqcup\limits_{i\in I_1} S^1 \ar[r] \ar@<-1ex>@{}[ur]^>{\llcorner} & \bigsqcup\limits_{i\in I_1} D^1 \ar[u]
}
\end{align*}
\end{definicion}

	Note that, for construction, $X$ has the weak structure with respect to the inclusions $X_i \to X$ because $X$ is the direct limit.
	Later in this document, we are going to obtain that every pseudotopological space is weak homotopy equivalent to a CW complex (see \autoref{theo:PsTopWECW}).

	In the following result, we claim that the product with any pseudotopological CW complex and the closed unit interval $[0,1]$, that is a CW complex, is also a pseudotopological CW complex.
	This going to be used in \autoref{prop:Hatcher_0.18}. 

\begin{proposicion}
\label{prop:CWtimesIisACWComplex}
	Let $X$ be a CW complex with decomposition $X_0\subset X_1\subset \cdots$.
	Then $X\times I$ is a CW complex with decomposition
\begin{align*}
Z_0=X_0\times\{0,1\},\ Z_n=((X_{n-1}\times I)\times (X_n\times \{0,1\}))/\sim
\end{align*}
such that the relation is between $(x,i)\in X_{n-1}\times I$ and $(x,i)\in X_n\times \{0,1\}$ for $i\in\{0,1\}$ and $x\in X_{n-1}$, $n\geq 1$.
\end{proposicion}

\begin{proof}
	Let $\Gamma_n$ be the set of $n$-cells of $X$ and $\Phi_{\alpha}^n:S_\alpha^n\rightarrow X_n$ for every $n\in\mathbb{N}_0$ and $\alpha\in\Gamma_n$.
	We denote by the $(\Phi_{\alpha}^n)^*$ the pushout of $\Phi_{\alpha}^n$.
	We can observe that
\begin{align*}
\xymatrix{
  \bigsqcup\limits_{\alpha\in \Gamma_{n-1}}\partial(D_\alpha^n\times I)\ar[r] \ar[d]_{(\Phi_\alpha^n)''} & \bigsqcup\limits_{\alpha\in \Gamma_{n-1}}D_\alpha^n\times I  \ar[d] & \\
 Z_n \ar[r] & X_n\times I \ar[r] & Z_{n+1} \\
 & \left(\bigsqcup\limits_{\alpha\in\Gamma_n}S_\alpha^n \right)\times\{0,1\}\ar[u] \ar[r] &  \left(\bigsqcup\limits_{\alpha\in\Gamma_n}D_\alpha^{n+1} \right)\times\{0,1\}\ar[u]
}
\end{align*}
Noting that
\begin{align*}
\partial(D_\alpha^n\times I) = (\partial D^n \times I)\cup (D^n\times\{0,1\})
\end{align*}
we define
\begin{align*}
(\Phi_\alpha^n)''\mid_{\partial D^n \times I} = & \bigsqcup\limits_{\alpha\in\Gamma_{n-1}}\Phi_{\alpha}^n \times \Id{I}\\
(\Phi_\alpha^n)''\mid_{D^n\times\{0,1\}} = & \bigsqcup\limits_{\alpha\in\Gamma_{n-1}}(\Phi_{\alpha}^n)^* \times \Id{\{0,1\}}
\end{align*}
	Concluding that $Z=\lim Z^n$ is a CW complex.
	
	By construction, the inclusions $Z_0 \subset X_0\times I \subset Z_1 \subset X_1 \subset Z_2\subset\ldots$ are subspaces. Let's name the morphisms
\begin{align*}
\alpha_i: & Z_i\rightarrow X_i\times I  & \beta_i: & X_i\times I\rightarrow Z_{i+1}\\
x_{i}': & X_i\times I\rightarrow X_{i+1}\times I  & z_i': & Z_i\rightarrow Z_{i+1} \\
x_i: & X_i\times I \rightarrow X\times I  & z_i: & Z_i\rightarrow Z
\end{align*}
	By the definition of colimits, we have unique $\Theta:X\times I\rightarrow Z$ and $\Psi:Z\rightarrow X\times I$ such that $z_{i+1} \beta_i = \Theta x_i$ and $ x_i\alpha_i= \Psi z_i$.
	We compute
\begin{align*}
\Psi\Theta x_i = & \Psi z_{i+1}\beta_i = x_{i+1}\alpha_{i+1}\beta_i = x_{i+1}x_i' = x_i \\
\Theta\Psi z_1 = & \Theta x_i\alpha_i = z_{i+1}\beta_i\alpha_i = z_{i+1}z'_i = z_i
\end{align*}
Concluding by uniqueness that $Z\cong X\times I$.
\end{proof}

	These result is extended in \autoref{Appx:ProductCWComplex}.
	In that section, we find out that the cartesian categorical product of every pair of pseudotopological CW complexes is a pseudotopological CW complex (\autoref{theo:CartesianProductCW}). 

\begin{definicion}
\label{def:CellularMap}
	Let $f:X\rightarrow Y$ be a map between CW complexes.
	We say that $f$ is a cellular map if $f(X^n)\subset Y^n$.
\end{definicion}

	Since $\cat{PsTop}$ has all its small limits and small colimits, for every map $f:X\rightarrow Y$ between pseudotopological spaces we can define the mapping cylinder $M_f$ as the quotient space of the disjoint union $(X\times I)\sqcup Y$ identifying $(x,1)\in X\times I$ with $f(x)\in Y$.
	Note that every mapping cylinder $M_f$ has a retraction to its base $Y\subset M_f$.
	
	A \emph{CW pair}\index{CW pair} is an order pair $(X,A)$ such that $X$ is a CW complex and $A$ is a \emph{subcomplex}\index{subcomplex}, i.e., $A$ is a closed set of $X$ with the closure operation induced by the pseudotopological structure on $X$.

\begin{lema}
\label{lem:MappingCylinder}
	Let $(X,A)$ be a CW pair, then $(X,M_i)$ is a CW pair with $i:A\rightarrow X$ the inclusion.
\end{lema}

\begin{proof}
	We can follow the proof of \autoref{prop:CWtimesIisACWComplex}, just adding the $n$-cells of $X$ that are not in $A$ just in the bottom face of the cylinder.
\end{proof}

	A \emph{deformation retraction} of a space $X$ onto a subspace $A$ is a continuous map $H:X\times I \rightarrow X$ such that $H(-,0) = \Id{X}$, $H(A,1)=A$ and $H\mid_{A,t}=\Id{A}$ for all $t$.

	We say that $(X,A)$ has the \emph{homotopy extension property}\index{homotopy extension property} if for each commutative solid arrow diagram
\begin{align*}
\xymatrix{
	& X \ar[rd]^{(\Id{X},0)} \ar@/^1pc/[rrd]^{g} \\
	A \ar[rd]_{(\Id{A},0)} \ar[ru]^{i} & & A\times I \ar@{..>}[r]^{H} & Y\\
	& A\times I \ar[ru]^{(i,\Id{I})} \ar@/_1pc/[rru]^{h}
}
\end{align*}
there exists a morphism $H:A\times I\rightarrow Y$ which makes the entire diagram commute.
	We can define that for every map replacing $i:A\to X$ for an arbitrary map $f:A\to X$.

\begin{lema}[\cite{Rieser_arXiv_2022}, 4.19]
	A map $f:A\rightarrow X$ has the homotopy extension property if and only if there exists a retract $r_k:A\times I \rightarrow X\times \{0\}\cup A\times I$ of the map $j_k: X\times \{0\}\cup A\times I\rightarrow A\times I$ in the pushout diagram
\begin{align*}
\xymatrix{
	& X \ar[rd]^{(\Id{X},0)} \ar@/^1pc/[rrd]^{i_0} \\
	A \ar[rd]_{(\Id{A},0)} \ar[ru]^{f} & & X\times \{0\}\cup A\times I \ar[r]^{j_k} & X\times I \\
	& A\times I \ar[ru]^{(i,\Id{I})} \ar@/_1pc/[rru]_{(f,\Id{I})}
}
\end{align*}
\end{lema}

\begin{proposicion}[\cite{Hatcher_2002}, Proposition 0.16]
\label{prop:Hatcher0.16}
	If $(X,A)$ is a CW pair, then $X\times \{0\}\cup A\times I$ is a deformation retract of $X\times I$, hence $(X,A)$ has the homotopy extension property
\end{proposicion}

\begin{proof}
	There is a retraction $r:D^n\times I\rightarrow D^n\times\{0\}$, for example \cite{Hatcher_2002} give us the radial projection from the point $(0,2)\in D^n\times\mathbb{R}$.
	Then setting $r_t=tr+(1-t)\Id{D^n\times I}$ gives a deformation retraction of $D^n\times I$ onto $D^n\times\{0\}\cup \partial D^n\times I$.
	This deformation gives rise to a deformation retraction of $X^n\times I$ onto $X^n\times \{0\}\cup(X^{n-1}\cup A^n)\times I$ since $X^n\times I$ is obtained from $X^n\times \{0\}\cup(X^{n-1}\cup A^n)\times I$ by attaching copies of $D^n\times I$ along $D^n\times\{0\}\cup \partial D^n \times I$.
	If we perform this deformation retraction during the interval $[1/2^{n+1},1/2^n]$ we obtain the desire result, remembering that $X$ has the weak structure with respect to the skeleton.
\end{proof}

\begin{definicion}
\label{def:Polyhedron}
	A polyhedron in $\mathbb{R}^n$ is a subspace which is the union of finitely many convex polyhedral.
	By a PL (piece-wise linear) map from a polyhedron to $\mathbb{R}^k$ we shall mean a map that is linear when restricted to each convex polyhedron in some such decomposition of the polyhedron into convex polyhedral.
\end{definicion}

\begin{lema}[\cite{Hatcher_2002}, Lemma 4.10]
\label{lem:Hatcher4.10}
	Let $f:I^n\rightarrow Z$ be a map, where $Z$ is obtained from a subspace $W$ by attaching a cell $e^k$.
	Then there is a homotopy $f_t: (I^n,f^{-1}(e^k))\rightarrow (Z,e^k) rel f^{-1}(W)$ from $f=f_0$ to a map $f_1$ for which there is a polyhedron $K\subset I^n$ such that
\begin{enumerate}
\item $f_1(K)\subset e^k$ and $f_1|K$ is PL with respect to some identification of $e^k$ with $\mathbb{R}^k$.
\item $K\supset f_1^{-1}(U)$ for some nonempty set $U$ in $e^k$.
\end{enumerate}
\end{lema}

\begin{proof}
	(The following proof is completely obtained by Hatcher without changes.)
	Identifying $e^k$ with $\mathbb{R}^k$, let $B_1,B_2\subset e^k$ be the closed balls of radius 1 and 2 centered at the origin.
	Since $f^{-1}(B_2)$ is closed and, therefore, compact in $I^n$, it follows that $f$ is uniformly continuous on $f^{-1}(B_2)$.
	Thus, there exists $\varepsilon > 0 $ such that $|x-y|<\varepsilon$ implies that $|f(x)-f(y)|<\frac{1}{2}$ for all $x,y\in f^{-1}(B_2)$.
	Subdivide the interval $I$ so that the induced subdivision of $I^n$ into cubes has each cube lying on a ball of diameter less than $\varepsilon$.
	Let $K_1$ be the union of every cube meeting $f^{-1}(B_1)$, and let $K_2$ be the union of all the cubes meeting $K_1$.
	We may assume that $\varepsilon$ is chosen smaller than half the distance between the compact sets $f^{-1}(B_1)$ and $I^n-f^{-1}(int(B_2))$, and then we have $K_2\subset f^{-1}(B_2)$.

	Now we subdivide all the cubes of $K_2$ into simplices. This can be done inductively.
	The boundary of each cube is a union of cubes of one lower dimension, so assuming these lower-dimensional cubes have already been subdivided into simplices, we obtain a subdivision of the cube itself by taking its center point as a new vertex and joining this by a cone to each simplex in the boundary of the cube.

	Let $g:K_2\rightarrow e^k=\mathbb{R}^k$ be the map that equals $f$ on all vertices of simplices of the subdivision and is linear on each simplex.
	Let $\varphi: K_2\rightarrow [0,1]$ be the map that is linear on simplices and has the value 1 on vertices in $K_1$ and 0 on vertices in $K_2-K_1$.
	Thus $\varphi(K_1)=1$.
	Define a homotopy $f_t:K_2\rightarrow e^k$ by the formula $(1-t\varphi)f+(t\varphi)g$, so $f_0=f$ and $f_1|K_1=g|K_1$. Since $f_t$ is the constant homotopy on simplices in $K_2$ disjoint from $K_1$, and in particular on simplices in the closure of $I^n-K_2$, we may extend $f_t$ to be constant homotopy of $f$ on $I^n-K_2$.

	The map $f_1$ takes the closure of $I^n-K_1$ to a compact set $C$ which, we claim, is disjoint from the center point 0 of $B_1$ and hence from neighborhood $U$ of $0$.
	This will prove the lemma, with $K=K_1$, since we will then have $f^{-1}(U)\subset K_1$ with $f_1$ PL on $K_1$ where it is equal to $g$.

	The verification of the claim has two steps:
\begin{enumerate}
	\item On $I^n-K_2$ we have $f_1=f$, and $f$ takes $I^n-K_2$ outside $B_1$ since $f^{-1}(B_1)\subset K_2$ by construction.
	
	\item For a simplex $\sigma$ of $K_2$ not in $K_1$ we have $f(\sigma)$ contained in some ball $B_\sigma$ of radius $\frac{1}{2}$ by the choice of $\varepsilon$ and the fact that $K_2\subset f^{-1}(B_2)$.
	Since $f(\sigma)\subset B_\sigma$ for all $t$, and in particular $f_1(\sigma)\subset B_\sigma$. We know that $B_\sigma$ is not contained in $B_1$ since $\sigma$ contains points outside $K_1$ hence outside $f^{-1}(B_1)$.
	The radius of $B_\sigma$ is half that of $B_1$, so it follows that 0 is not in $B_\sigma$, and hence 0 is not in $f_1(\sigma)$.\qedhere
\end{enumerate}
\end{proof}

\begin{teorema}[\cite{Hatcher_2002}, Theorem 4.8]
\label{theo:Hatcher_4.8}
	Every map $f:X\rightarrow Y$ of CW complexes is homotopic to a cellular map.
	If $f$ is already cellular on a subcomplex $A\subset X$, the homotopy may be taken to be stationary on $A$.
\end{teorema}

\begin{proof}
	Suppose inductively that $f:X\rightarrow Y$ is already cellular on the skeleton $X^{n-1}$ and let $e^n$ be an $n$-cell of $X$.
	Take the image of the characteristic map of $e^n$, $\Phi: D_\alpha^n\rightarrow X$, which it is compact by being the image of a compact, then $f\Phi_\alpha$ also is.

	By theorem 5.16 from \cite{Rieser_arXiv_2022}, $f\circ\Phi_\alpha (D_\alpha)$ intersect the interior of only finitely many cells.
	Thus $f(e^n)$ also meets only finitely many cells of $Y$.

	Let $e^k\subset Y$ be a cell of highest dimension meeting $f(e^n)$.
	We may assume that $k>n$, otherwise $f$ is already cellular on $e^n$.

	By \autoref{lem:Hatcher4.10}, composing the given map $f:X^{n-1}\cup e^n\rightarrow Y^k$ with a characteristic map $I^n\rightarrow X$ for $e^n$, we obtain a map $f$ as in the lemma, with $Z=Y^k$ and $W=Y^k-e^k$.
	The homotopy given by the lemma is fixed on $\partial I^n$, hence induces a homotopy $f_t$ of $f|_{X^{n-1}\cup e^n}$ staying fixed on $X^{n-1}$.
	The image of the resulting map $f_1$ intersects the open set $U$ in $e^k$ in a set contained in the union of finitely many hyperplanes of dimension at most $n$, so if $n<k$ there will be points $p\in U$ not in the image of $f_1$.
	Then we can deform $f|_{X^{n-1}\cup e^n}\ \text{rel}\ X^{n-1}$ so that $f(e^n)$ misses the whole cell $e^k$ by composing with a deformation retraction of $Y^k-{p}$ onto $Y^k-e^k$.

	By finitely many iterations of this process we eventually make $f(e^n)$ miss all cells of dimension greater than $n$.
	Doing this for all $n$-cells simultaneously , staying fixed on $n$-cells in $A$ where $f$ is already cellular, we obtain a homotopy of $f|_{X^n}~\text{rel} (X^{n-1}\cup A^n)$ to a cellular map. 
	The induction step is then completed by appealing to the homotopy extension 
\end{proof}

\begin{definicion}
\label{def:retraction}
	Let $r:X\rightarrow Y$ be a continuous map between pseudotopological spaces.
	We say that $r$ is a \emph{retraction}\index{retraction}, and $Y$ is a \emph{retract}\index{retract} of $X$, if there exists $i:Y\rightarrow X$ such that $r\circ i=\Id{Y}$.
	In general, we say that $X$ dominate $Y$, or $Y$ is \emph{dominated by}\index{dominated by} $X$, if $r\circ i\simeq \Id{X}$.
\end{definicion}

\begin{proposicion}[\cite{Hatcher_2002}, 0.18]
\label{prop:Hatcher_0.18}
	If $(X_1,A)$ is a pair such that there exists a deformation retraction of $X_1\times I$ onto $X_1\times\{0\}\cup A\times I$.
	If we have attaching maps $f,g:A\rightarrow X_0$ that are homotopic, then $X_0 \sqcup_f X_1 \simeq X_0\sqcup_g X_1 \text{rel} (X_0)$.
\end{proposicion}

\begin{proposicion}[\cite{Hatcher_2002}, A.11]
\label{prop:Hatcher_A.11}
	A space dominated by a CW complex is homotopy equivalent to a CW complex.
\end{proposicion}

\begin{proof}
	Let's define the mapping telescope $T(f_1,f_2,\ldots)$ of a sequence of maps $X_1 \xrightarrow{f_1}  X_2 \xrightarrow{f_2}  X_3 \xrightarrow{f_3} \cdots$ is the quotient space of $\sqcup_i (X_i\times [i,i+1])$ obtained by identifying $(x,i+1)\in X_i\times [i,i+1]$ with $(f(x),i+1)\in X_{i+1}\times [i+1,i+2])$. We are going to prove the following facts:
\begin{enumerate}
	\item $T(f_1,f_2,\ldots) \simeq T(g_1,g_2,\ldots)$ if $f_i\simeq g_i$ for each $i$.
	
	\item $T(f_1,f_2,\ldots) \simeq T(f_2,f_3,\ldots)$.
	
	\item $T(f_1,f_2,\ldots) \simeq T(f_2f_1,f_4f_3,f_6f_5,\ldots)$.
\end{enumerate}

	To prove \textbf{(2)}, we note that $X_1\times [1,2] \simeq X_1\times {2}$.
	Thus $T(f_1,f_2,\ldots) \simeq T(f_2,f_3,\ldots)$.

	To prove both (1) and (3), we use the \autoref{prop:Hatcher_0.18}.

	First, to prove \textbf{(1)}, let's write the pair $(Y_1,A)$ with $Y_1\coloneqq \sqcup_i (X_i\times [i,i+1])$ and $A\coloneqq \sqcup_i (X_i\times\{i,i+1\})$.
	We define $Y_0\coloneqq \sqcup_{i\in\mathbb{N}} X_i\times\{i\}$ and define the continuous maps
\begin{align*}
f:A\rightarrow Y_0 & \text{ such that }f|_{X_i\times\{i\}}=\Id{X_i\times\{i\}} \text{ and } f|_{X_i\times\{i+1\}}=f_i\times \{i+1\}\\
g:A\rightarrow Y_0 & \text{ such that }g|_{X_i\times\{i\}}=\Id{X_i\times\{i\}} \text{ and } g|_{X_i\times\{i+1\}}=g_i\times \{i+1\}
\end{align*}
	By hypothesis, $f\simeq g$, thus we got the desired result.

	Now, to prove \textbf{(3)}, let's write the pair $(Y_1, A)$
	with
	\begin{align*}
	Y_1\coloneqq & \sqcup_i ((X_{2i-1}\times [2i-1,2i])\sqcup M_{f_{2i}}) \text{ and }\\
	A\coloneqq & \sqcup_i ((X_{2i-1}\times \{2i-1,2i\})\sqcup M_{f_{2i}}).
	\end{align*}
	We define $Y_0\coloneqq \sqcup_i ((X_{2i-1}\times\{2i-1\} )\sqcup M_{f_{2i}})$ and the continuous maps
\begin{align*}
f:A\rightarrow Y_0 &\text{ such that } f\mid_{X_{2i-1}\times\{2i-1\}} = \Id{X_{2i-1}\times\{2i-1\}},\ f\mid_{X_{2i}\times\{2i-1\}}=f_{2i-1}\times\{2i\}\\
& \text{ and }f\mid_{M_{f_{2i}}}=\Id{M_{f_{2i}}}\\
g:A\rightarrow Y_0 & \text{ such that } g\mid_{X_{2i-1}\times\{2i-1\}}=\Id{X_{2i-1}\times\{2i-1\}}, g\mid_{X_{2i}\times\{2i-1\}}=f_{2i}f_{2i-1}\times\{2i+1\}\\
& \text{ and }g\mid_{M_{f_{2i}}}=r_{2i}
\end{align*}
with $r_{2i}$ the retraction of the mapping cylinder to its base.
	Getting the result by the same \autoref{prop:Hatcher_0.18}.

	With these brief results,  we can prove our statement.
	Suppose that the space $Y$ is dominated by the CW complex $X$ via maps $Y \xrightarrow{i} X \xrightarrow{r} Y$ with $ri\simeq \Id{Y}$.
	By (2) and (3) we have $T(ir,ir,ir,\ldots) \simeq T(r,i,r,i,\ldots) \simeq T(i,r,i,r,\ldots) \simeq T(ri,ri,\ldots)$.
	Since $ri=\Id{X}$, $T(ri,ri,\ldots) \simeq Y\times [0,\infty) \simeq Y$.
	On the other hand, the map $ir$ is homotopic to a cellular map $f:X \to X$ by \autoref{theo:Hatcher_4.8}.
	Thus $T(ir,ir,\ldots) \simeq T(f,f,\ldots)$, which is a CW complex.
\end{proof}

	It is direct corollary from \autoref{prop:Hatcher_A.11} that:

\begin{corolario}
\label{coro:RetractPsTopCW}
	Every retract of a pseudotopological CW complex is a pseudotopological CW complex.
\end{corolario}

\section{Simplicial Sets Category and the Geometric Realization}
\label{cap.SSetCat}

	In this section we recall the category of simplicial sets ($\cat{sSet}$) and we define the pseudotopological geometric realization from a simplicial set to a pseudotopological space.
	We define that functor replacing the $\cat{Top}$ structures and colimits for the ones in $\cat{PsTop}$.
	
	The introductory definitions comes from, or are analogous to definitions in, \cite{Heuts_2022_SimplicialAndDendroidal}.
	Ending this section we are going to discover that every geometric realization is a (pseudotopological) CW complex.
	Straightforward, we begin introducing the simplicial category to be able to introduce the $\cat{sSet}$ category and then create a functor to $\cat{PsTop}$\footnote{Spoiler alert: this functor form a Quillen Equivalence, but we observe that much later.}.

\begin{definition}[\cite{Heuts_2022_SimplicialAndDendroidal}, page 49]
\label{def:SimpCategory}
	The simplicial category $\cat{\Delta}$\index{category!simplicial}\nomenclature{$\cat{\Delta}$}{The simplicial category} is the category with objects $[n]\coloneqq \{0,1,\ldots,n\}$ for $n \in \mathbb{N}_0$, and $f\in\Delta([n],[m])$ if and only if $k\leq k'$ implies that $f(k)\leq f(k')$ for every $k,k'\in [n]$.
\end{definition}

	We particularly emphasize the following morphisms: For each $i\in [n]$ the injective monotone function

\begin{align*}
\delta_i:[n-1]\rightarrow [n]
\end{align*}
which skips the value $i$. Also, for each $j\in [n-1]$ there is a surjective function

\begin{align*}
\sigma_j:[n]\rightarrow [n-1]
\end{align*}
given by $\sigma_j(k)=k$ for $k\leq j$ and $\sigma(k)=k-1$ for $k>j$.
	These morphisms are called the elementary faces and elementary degeneracies, respectively.
	Every injective map can be written as a composition of elementary and, in analogous way, every surjective map can be written as a composition of elementary degeneracies.
	It is relative easy to observe the following relations:
\begin{enumerate}
\item $\delta_j\delta_i=\delta_{i}\delta_{j-1}$ for $i<j$.
\item $\sigma_i\sigma_j=\sigma_{j-1}\sigma_i$ for $i<j$.
\item $\sigma_i\delta_j = \left\lbrace\begin{array}{ll}
\delta_{j-1}\sigma_i & \text{ if }i<j-1\\
id & \text{ if }i=j-1 \text{ or }i=j\\
\delta_j\sigma_{j-1}& \text{ if }i>j
\end{array}  \right.$.
\end{enumerate}

	The simplicial set category is a small category which does not have a lot of limits and colimits; however, it satisfies the following interesting properties.

\begin{proposicion}[\cite{Heuts_2022_SimplicialAndDendroidal}, 2.1]
\label{Prop:Heuts2.1}
	Let $0\leq i<j<n$ and consider the pushout
\begin{align*}
\xymatrix{
[n] \ar[r]^{\sigma_{j}} \ar[d]^{\sigma_{i}} & [n-1] \ar[d]^{\sigma_{i}} \\
[n-1] \ar[r]^{\sigma_{j-1}} & [n-2]
}
\end{align*}
\begin{enumerate}[label = (\roman*)]
	\item In the square above, there exist sections (in this context it means right inverses) $\alpha:[n-1]\rightarrow [n]$ of $\sigma_i$ (left) and $\beta:[n-2]\rightarrow [n-1]$ of $\sigma_i$ (right), which are compatible in the sense that $\sigma_j\alpha = \beta\sigma_{j-1}$.
	
	\item Consider a commutative square
\begin{align*}
\xymatrix{
A \ar[r]^{q} \ar[d]^{p} & B \ar[d]^{r} \\
C \ar[r]^{s} & D
}
\end{align*}
in a category $\mathcal{C}$, with $p$ a split epimorphism (i.e., an epimorphism admitting a section) and $q$ an epimorphism.
	If there exist compatible sections $\alpha$ of $p$ and $\beta$ of $r$ (in the sense described in (i)), then the square is a pushout.
	In fact it is an absolute pushout meaning any functor from $\Delta$ to another category sends the square to a pushout square.

	\item Let $[k]\xleftarrow{p} [n] \xrightarrow{q} [l]$ be surjections in $\Delta$.
	The pushout
\begin{align*}
\xymatrix{
[n] \ar[r]^{q} \ar[d]^{p} & [l] \ar[d]^{r} \\
[k] \ar[r]^{s} & [m]
}
\end{align*}
exists in $\Delta$ and is an absolute pushout.
\end{enumerate}
\end{proposicion}

\begin{proof}
\
\begin{enumerate}[label = (\roman*)]
	\item Define $\alpha=\delta_i[n-1]\rightarrow [n]$ and $\beta=\delta_i:[n-2]\rightarrow [n-1]$.
	The equations $\sigma_j\delta_i = \delta_i\sigma_{j-1}$ is discussed above remembering that $i<j$.
	
	\item Consider $X$ an object of $\mathcal{C}$ and the morphisms $f:B\rightarrow X$, $g:C\rightarrow X$ are maps such that $fq=gp$, then ones defines a corresponding map $h:D\rightarrow X$ by $h:f\beta$.
	We should check that $hs=g$ and $hr=f$.
	The first equality is clear from $hs=f\beta s = fq\alpha = gp\alpha = g$.
	For the second equality, given that $q$ is an epimorphism, it suffices to prove that $hrq=fq$. The left-hand side equals $hrq=f\beta s p = fq\alpha p = gp\alpha p=gp=fq$.
	Since $r$ is an epimorphism, if there exists other $h':D\rightarrow X$ such that $hr=h'r$, then $h=h'$.

	Observe that this proof depends only on categorical definitions of morphisms (like epimorphisms and commutative square) therefore it is preserved for any functor.

	\item Since the surjections $p$ and $q$ can both be factored as composition of elementary degeneracies, so that the conclusion follows by repeatedly applying (i) and (ii).\qedhere
\end{enumerate}
\end{proof}

\begin{proposicion}[\cite{Heuts_2022_SimplicialAndDendroidal}, 2.2]\label{Prop:Heuts2.2}~
\begin{enumerate}[label = (\roman*)]
	\item If $f:[m]\rightarrow [n]$ is a monomorphism, then the pullback $g:[k]\rightarrow [n]$ along $f$ exists, provided that the image of $g$ intersects  the image of $f$ non trivially (the intersection is non-empty.)
	
\item If $f:[m]\rightarrow [n]$ is a epimorphism, then the pushout $g:[k]\rightarrow [n]$ along $f$ exists.
\end{enumerate}
\end{proposicion}

\begin{proof}
	Since $f$ is a monomorphism, then it has a right inverse (to say $h$) and the commutative square
\begin{align*}
\xymatrix{
g^{-1}f([m]) \ar[r]^{hg} \ar[d]^{i} & [m] \ar[d]^{f} \\
[k] \ar[r]^{g} & [n]
}
\end{align*}
is a pullback.
	Observe that this does not work with trivial intersection because the empty is not an object in $\Delta$.

	For (ii), it suffices to treat the case where $f$ is an elementary degeneracy $\sigma_i:[m]\rightarrow [m-1]$. Then the pushout of $g$ is the map which collapses the interval $[g(i),g(i+1)]$ to a single point.
\end{proof}

\begin{definition}[\cite{Heuts_2022_SimplicialAndDendroidal}, page 53]
\label{def:SimplicialSets}
	The \emph{category of simplicial sets}, denoted by $\cat{sSet}$,\index{Category!Simplicial sets}\nomenclature{$\cat{sSet}$}{The category of simplicial sets} is defined with functors from $\cat{\Delta}^{op}$ to $\cat{Set}$ as objects and natural transformations between such functors as morphisms.

Let $K\in \ocat{sSet}$, we denote $K([n])$ by $K_n$; we also denote by $\alpha^*:K_n\rightarrow K_m$ the functions $K(\alpha:[n]\rightarrow [m])$.
\end{definition}

	Equivalently, we have the face maps and degeneracy maps of the simplicial object $X$:
\begin{align*}
d_i = (\delta_i)^*:X_n\rightarrow X_{n-1} &,\ i\in [n]\\
s_j = (\sigma_j)^*:X_{n-1}\rightarrow X_n &,\ j\in [n-1]
\end{align*}

	We define $r[n]\coloneqq \{ (x_0,\ldots,x_n) \mid \sum_{i=0}^{n}=1, x_i\geq 0\}$ for every $n\geq 0$.
	Then, any function of sets $f:[m]\rightarrow [n]$ defines an affine map $f_*:r[m]\rightarrow r[n]$ which is uniquely determined by the vertices.

	Trying to simplify some proof in the future, we define the geometric realization as a coend and show that for this specific functor we obtain a quotient from the disjoint union of $K_n$ for $K\in \ocat{sSet}$.
	First, we define what a coend is

\begin{definition}
\label{def:COwedge}
	Let $\mfunt{F}{C^{op}\times C}{X}$ be a functor. A \emph{wedge}\index{wedge} $e:w\rightarrow \funt{F}$ (resp. \emph{cowedge}\index{cowedge} $e:\funt{F}\rightarrow w$) is an object $w$ in $\ocat{X}$ and maps $e_c:w\rightarrow \funt{F}(c,c)$ (resp. $e_c:\funt{F}(c,c)\rightarrow w$) for each $c$, such that given any morphism $f:c\rightarrow c'$, the following diagram commutes:
\begin{align*}
\vcenter{\xymatrix{
	w \ar[r]^{e_{c'}} \ar[d]^{e_c} & \funt{F}(c',c') \ar[d]^{\funt{F}(f,\Id{c'})} \\
	\funt{F}(c,c) \ar[r]^{\funt{F}(\Id{c},f)} & \funt{F}(c,c')
}} & & 
\left(resp.\vcenter{\xymatrix{
	\funt{F}(c,c') \ar[r]^{\funt{F}(f,\Id{c'})} \ar[d]^{\funt{F}(\Id{c},f)} & \funt{F}(c',c') \ar[d]^{e_{c'}} \\
	\funt{F}(c,c) \ar[r]^{e_c} & w
}}\right)
\end{align*}
\end{definition}

\begin{definition}
\label{def:COend}
	Let $\mfunt{F}{C^{op}\times C}{X}$ be a functor. An \emph{end}\index{end} (resp. a \emph{coend}\index{coend}) of $\funt{F}$ is a universal wedge (resp. cowedge), that is a wedge $e:w\rightarrow \funt{F}$ (resp. a cowedge $e:\funt{F}\rightarrow w$) such that any other wedge $e':w'\rightarrow \funt{F}$ (resp. $e':\funt{F}\rightarrow w'$) factor through $e$ via a unique map $w'\rightarrow w$ (resp. $w\rightarrow w'$).
\end{definition}

\begin{definition}[\cite{Robert_Thesis2013}, page 39]
\label{def:GeometricRealization}
	The \emph{geometric realization of a simplicial set} $K\in \ocat{sSet}$ is defined as the coend
\begin{align*}
|K|\coloneqq \int^{[n]\in \Delta} K_n\cdot r([n])
\end{align*}
where $K_n\cdot r([n])$ is the coproduct $\sqcup_{K_n}r([n])$.
	Since $\Delta$ is a small category, and every topological construct is cocomplete, then we can write $|K|$ as the following coequalizer:
\begin{equation}
\label{eq:geomRealCoEq}
\bigsqcup_{[n]\rightarrow [m]\in\Delta} K_m\cdot r[n] \rightrightarrows \bigsqcup_n K_n\cdot r[n] \rightarrow |K|
\end{equation}
where the two parallel arrows given by the maps $K_m \times r[n] \rightarrow K_n\times r[n]$ and $K_m\times r[n]\rightarrow K_m\times r[m]$ induced by the associated functions $[n]\rightarrow [m]$.
\end{definition}

	The coequalizer in \ref{eq:geomRealCoEq} is defined for every topological construct that is a supcategory of $\cat{Top}$.
	That definition is one of the traditional equivalent definitions for the geometric realization into $\cat{Top}$. When it is not obvious in the context, we use $|\cdot|$ to denote the topological realization and $|\cdot|_\cat{PsTop}$ to denote the pseudotopological realization.

	We can find examples where the result of the PsTop realization is not a topology, and then both functors are not the same.
	However, in \autoref{prop:SpheresHomeoGeoReal}, we emphasize particular examples of simplicial sets such that both constructions produce the same object.

\begin{definicion}[\cite{Robert_Thesis2013}, page 38]
\label{def:BoundaryHorns}
	Let $\mfunt{\Delta[\cdot]}{\Delta}{sSet}$ defined by $\funt{\Delta[n]}[k]\coloneqq \{f:[k]\rightarrow [n]\mid x<y\Rightarrow f(x)<f(y)\}$.

	The \emph{boundary} of the simplicial set $\Delta[n]$, denoted by $\partial\Delta[n]$ is obtained by
\begin{align*}
\bigsqcup_{\substack{\text{$[n-2]\rightarrow [n]$ inj}}} \Delta[n-2] \rightrightarrows \bigsqcup_{\substack{\text{$[n-1]\rightarrow [n]$ inj}}} \Delta[n-1]\rightarrow \partial\Delta[n]
\end{align*}

	For $0\leq r \leq n$, we define the \emph{$r$-horn of} $\Delta[n]$, denoted by $\Lambda^r[n]$ and given by
\begin{align*}
\bigsqcup_{\substack{\text{$[n-2]\rightarrow [n]$ inj} \\ \text{$r$ in the image}}} \Delta[n-2] \rightrightarrows \bigsqcup_{\substack{\text{$[n-1]\rightarrow [n]$ inj} \\ \text{$r$ in the image}}} \Delta[n-1]\rightarrow \Lambda^r[n]
\end{align*}
\end{definicion}

\begin{proposicion}
\label{prop:SpheresHomeoGeoReal}
	$|\Delta[n]|_\cat{PsTop} \cong |\Delta[n]|$, $|\partial\Delta[n]|_\cat{PsTop} \cong |\partial\Delta[n]|$ and $|\Lambda^r[n]|_\cat{PsTop} \cong |\Lambda^r[n]|$.
\end{proposicion}

\begin{proof}
	Let $l\in\mathbb{N}_0$.
	Since we can write a coequalizer as a cocone and because $K_m\cdot r[n]$ are topological spaces, we have the following diagram
\begin{align*}
\xymatrix{
\bigsqcup_{[n]\rightarrow [m]} \Delta[l]_m\cdot r[n] \ar@<0.5ex>[rr]\ar@<-0.5ex>[rr] \ar[rd] \ar[rdd] & & \bigsqcup_n \Delta[l]_n\cdot r[n] \ar[ld]\ar[ldd] \\
&|\Delta[l]|_\cat{PsTop} \ar@{.>}[d]_{\Psi}\\
&|\Delta[l]|
}
\end{align*}
	By forgetful functor, we obtain that $\Psi$ is bijective.
	Since $\Delta[l]_n$ is finite for every $n$ and in particular is empty for every natural greater than $l$, we obtain that $\sqcup_n \Delta[l]_n\cdot r[n]$ is compact, then its image $|\Delta[l]|_\cat{PsTop}$ is compact.
	Thus, by proposition 4.25 from \cite{kapulkin_edel2023_SyntheticApproach} paper $\Psi$ is an homemomorphism.

	The other case follows directly of the fact that $|\cdot|_\cat{PsTop}$ is a left adjoint\footnote{We prove this later, after \autoref{def:SingPsTop}.}, then it preserves colimits.
	Thus
\begin{align*}
\left| \bigsqcup \Delta[k] \right|_\cat{PsTop} = \bigsqcup |\Delta[k]|_\cat{PsTop}
\end{align*}
and by similar arguments for the cocone, we obtain that $|\partial\Delta[n]|_\cat{PsTop} \cong |\partial\Delta[n]|$ and $|\Lambda^r[n]|_\cat{PsTop} \cong |\Lambda^r[n]|$.
\end{proof}

	We rephrase the definition of the geometric realization.
	By theorem 1.2.1.8 in \cite{Preuss_2002}, we have that every coequalizer is a quotient which only depends on the relations produced by the maps (as morphisms in $\cat{Set}$.)
	Thus, we can see the coequalizer of the pseudotopological geometric realization as quotients in $\cat{PsTop}$, this is the same process that it is made for the category $\cat{Top}$.
	In other words, let $X\in \ocat{sSet}$, then
\begin{align*}
|X|= \bigsqcup_{n\geq 0} X_n\times r[n]/\sim
\end{align*}
with the identification $(x,\alpha_*t)\sim (\alpha^*x,t)$ for each morphism $\alpha\in \ocat{sSet}$.

	An $n$-simplex $x\in X_n$ is called \emph{degenerated}\index{degenerated} if it lies in the image of one of the degeneracy operators $s_i:X_{n-1}\rightarrow X_n$ for $i\in [n-1]$.
	Since $\mathbb{N}_0$ is lower bounded, then for every $x\in X_n$ which is degenerate, there exists a $z$ which is non-degenerate $k$-simplex and $\beta:[n]\rightarrow [k]$ such that $x=\beta^*z$. Observe that every $0$-simplex is non-degenerate.

\begin{lema}[\cite{Heuts_2022_SimplicialAndDendroidal}, page 57]
\label{lem:UniqueSelection}
	Let $x$ be a degenerated $n$-simplex, the selection of $(\beta,z)$ is unique.
\end{lema}

\begin{proof}
	Suppose that other pair $(\gamma,w)$ such that $\gamma^*w=x$ and $w$ is a non-degenerate $l$-simplex.
	Since $\beta$ is an epimorphism, we can form the pushout
\begin{align*}
\xymatrix{
	[n] \ar[r]^{\beta} \ar[d]^{\gamma} & [k] \ar[d]\\
	[l] \ar[r] & [j]
}
\end{align*}
It is an absolute pushout by \autoref{Prop:Heuts2.1} and therefore the resulting square
\begin{align*}
\xymatrix{
	X_n & X_k \ar[l]^{\beta^*}\\
	X_l \ar[u]^{\gamma^*} & X_j \ar[l] \ar[u]
}
\end{align*}
is a pullback (a pushout in the opposite category.)
	Thus there is an element $v\in X_j$ whose image is $z$ (respectively $w$) in $X_k$ (resp. $X_l$.)
	By the assumption that $w$ and $z$ are non-degenerate, this can only happen if the maps $[k]\rightarrow [j]$ and $[l]\rightarrow [j]$ are identities.
	It follows that $\beta=\gamma$ and $z=w$.
\end{proof}

	Every geometric realization has a CW complex structure.
	We describe that structure through non-degenerate simplices.
	At first, let's denote by $x\otimes t$ for the point of $X$ determined by a pair $(x,t)\in X_k\times \Delta_k$.
	For a $n$-simplex $x$, we write $\hat{x}:\Delta^n\rightarrow |X|$ considering the quotient map restricted to $x\times \Delta^n$ to $|X|$.

\begin{lema}
\label{lem:UniversalPropertyQuotient}
	Let $\cat{C}$ be a topological construct.
	Let $X$ an object in $\cat{C}$ and $\sim$ an equivalence relation on $X$.
	If $f:X\rightarrow Y$ a morphism in $\cat{C}$ such that $a\sim b$ implies $f(a)=f(b)$ for all $a$ and $b$ in $X$, then there exists a unique continuous map $g:X/\sim\rightarrow Y$ such that $g=fq$.
\end{lema}

\begin{proof}
	This follows because every quotient can be written as a coequalizer in a topological construct.
\end{proof}

\begin{proposicion}[Skeletal filtration; \cite{Heuts_2022_SimplicialAndDendroidal}]
\label{prop:SkeletalFiltrarionGeomReal}
	Denote $|X|^{(n)}$ the subspace of $X$ consisting of the points which can be represented as $x\otimes t$ for some $(x,t)\in X_k\times\Delta^k$ with $k\leq n$.
	This describes a filtration of $X$,
\begin{align*}
|X|^{(0)}\subset |X|^{(1)}\subset |X|^{(2)}\subset \cdots \subset \cup_{n\in\mathbb{N}_0} |X|^{(n)} = |X|
\end{align*}
and $|X|$ has the weak structure with respect to these subspaces.
\end{proposicion}

\begin{proof}
	By construction, we have the inclusions.
	For $n\in\mathbb{N}_0$, we can construct the following commutative diagram
\begin{align*}
\xymatrix{
	\bigsqcup\limits_{i\leq n} X_i\cdot r[i] \ar[r]^{\alpha_n} \ar[d]^{\subset_{\sqcup n}} & |X|^{(n)} \ar[d]^{\subset_{|n|}}\\
	\bigsqcup\limits_{i\in\mathbb{N}_0}  X_i\cdot r[i] \ar[r]^{q} & |X|
}
\end{align*}
	Since $|X|^{(k)}$ has the subspace structure of $|X|$, the morphisms $\alpha_k:\sqcup_{i\leq k} X_i r[i] \rightarrow |X|^{(k)}$ are continuous (by construction it is a epimorphism too.)
	Then, if we have the cocone
\begin{align*}
\xymatrix{
	|X|^{(m)} \ar[dd]_{\subset_{|mn|}} \ar[rd]^{\subset_{|m|}} \ar@/^1pc/[rrrd]^{\beta_m} \\
	& |X| & & Y \\
	|X|^{(n)} \ar[ru]^{\subset_{|n|}} \ar@/_1pc/[rrru]^{\beta_n}
}
\end{align*}
we obtain that there exists an unique $h:\sqcup_{i\in\mathbb{N}_0}X_i r[i]\rightarrow$ such that
\begin{align*}
\xymatrix{
	 \bigsqcup\limits_{i\leq m}X_i\cdot r[i]\ar[dd]_{\subset_{\sqcup mn}}  \ar[rd]^{\subset_{\sqcup m}} \ar[r]^{\alpha_m} & |X|^{(m)} \ar@/^1pc/[rrd]^{\beta_m} \\
	& \bigsqcup\limits_{i\in\mathbb{N}_0}X_i\cdot r[i] \ar@{..>}[rr]^{h} & & Y \\
	 \bigsqcup\limits_{i\leq n}X_i\cdot r[i]\ar[ru]^{\subset_{\sqcup n}} \ar[r]^{\alpha_n} & |X|^{(n)} \ar@/_1pc/[rru]^{\beta_n}
}
\end{align*}
is commutative because $\lim_k \sqcup_{i\leq k} X_i\cdot r[i] = \sqcup_{i\in\mathbb{N}_0} X_i\cdot r[i]$.
	Note that $h(x,t)=h(y,s)$ for every $x\otimes t = y\otimes s$.
	Therefore, by the lemma above, there exists a unique $h': |X|\rightarrow Y$ such that $h=h'q$.
	Since $h\sqcup_m = \beta_m\alpha_m$, then we obtain $h'q\sqcup_m=\beta_m\alpha_m$ and this implies that $h'\subset_{|m|}\alpha_m=\beta_m\alpha_m$.
	Thus $h'\subset_{|m|}=\beta_m$ for $\alpha_m$ being an epimorphism.
\end{proof}

	Observe that the previous proof works for every topological constructs such that it has these triangles $r[n]$ as objects; however, we are only interested in $\cat{Top}$ and $\cat{PsTop}$ in this paper.

	In \cite{Heuts_2022_SimplicialAndDendroidal}, it is proved for the category $\cat{Top}$ that the filtration in \autoref{prop:SkeletalFiltrarionGeomReal} serves as the skeletal filtration of the CW structure of $|X|$.
	First we claim that the space the space $|X|^{(0)}$ is discrete and is in fact given by $X_0\times \Delta_0$, so that the elements of $X_0$ will serve as the 0-cells of $|X|$.
	Indeed, it is clear that the map $X_0\times \Delta^0 \rightarrow |X|^{(0)}$ is surjective.
	To see that is a bijection, we should argue that no two distinct 0-simplices of $X$ are identified in $|X|$.
	If $x,y\in X_0$, then $x\otimes 1$ and $y\otimes 1$ represent the same point of $|X|$ only if there exists $z\otimes t$ with $z\in X_n$ and $t\in \Delta^n$, together with morphisms $\alpha, \beta:[0]\rightarrow [n]$ such that $\alpha^*z=x$, $\beta^*=y$ and $\alpha_* 1=\beta_* 1=t$. The last condition implies that $\alpha=\beta$, from which it follows that $x=y$.

	To observe that $|X|^{(n)}$ can be seen as $|X|^{(n-1)}$ attaching some $n-1$ cells, we need the following technical result:
	The commutative square
\begin{align*}
\xymatrix{
\bigsqcup\limits_{x\in nd(X_n)} \partial \Delta^n \ar[r] \ar[d] & |X|^{(n-1)} \ar[d]\\
\bigsqcup\limits_{x\in nd(X_n)} \Delta^n \ar[r] & |X|^{(n)}
}
\end{align*}
	where $nd(X_n)$
	\nomenclature{$nd(X_n)$}{The non-degenerate simplices of $X_n$}
	denoted the subset of $X_n$ consisting of non-degenerate $n$-simplices, is a pushout.
	To prove it, we replicate the proof in 2.4 from \cite{Heuts_2022_SimplicialAndDendroidal}.

\begin{proposicion}[\cite{Heuts_2022_SimplicialAndDendroidal}, 2.4]
\label{prop:Heuts_2.4}
	Let $\xi\in |X|$.
	Choose $x\in X_n$ and $t\in\Delta^n$ with $\xi=x\otimes t$ and with $n$ as small as possible.
	Then $x$ is non-degenerate and if $n\geq 1$ then $t$ is contained in the interior of $\Delta^n$.
	Also, the pair $(x,t)$ representing $\xi$ with $x$ non-degenerated $t$ in the interior of $\Delta^n$ is unique.
\end{proposicion}

\begin{proof}
	First we show $x$ is non-degenerate.
	Suppose that $x$ is degenerate; then there would exist a non-trivial surjection $\alpha:[n]\rightarrow [m]$ and $y\in X_m$ with $\alpha^*y=x$.
	But then $x\otimes t=\alpha^*y\otimes t = y\otimes\alpha_*t$ contradicting the minimality of $n$.
	It is also straightforward to see that $t$ is in the interior of $\Delta^n$ (assuming $n\geq 1$); indeed, if it were on the boundary $\partial \Delta[n]$ then there would exist a nontrivial injection $\beta:[k]\rightarrow [n]$ such that $t$ is in the image of $\beta_*:\Delta^k\rightarrow \Delta^n$, so we may write $t=\beta_*s$.
	In that case $x\otimes t = x\otimes \beta_*s = \beta^*x\otimes s$, which again contradicts the assumption that $n$ is minimal.

	In remain to argue that the representative pair $(x,t)$ of the proposition is unique.
	So suppose $\xi=x\otimes t = y\otimes s$ where both $x$ and $y$ are non-degenerate and $s,t$ are interior points of $\Delta^n$.
	By equivalence relation involved in the definition of $\otimes$, this means that there is a zigzag in $\Delta$ of the form
\begin{align*}
\xymatrix{
 [n] & & [m_1]  & & [m_{N-1}] & & [n] \\
 & [k_1] \ar[lu]^{\alpha_1} \ar[ru]_{\beta_1} & &  \cdots \ar[lu]^{\alpha_3} \ar[ru]_{\beta_{N-1}} & & [k_N] \ar[lu]^{\alpha_N} \ar[ru]_{\beta_N}
}
\end{align*}
and elements $(a_i,u_i)\in X_{k_i}\times\Delta^{k_i}$, $(b_i,v_i)\in X_{m_i}\times \Delta^{m_i}$ for which
\begin{align*}
\alpha_i^*b_{i-1} = a_i, & (\alpha_i)_*u_i=v_{i-1},\\
\beta_i^*b_i = a_i, & (\beta_i)_*u_i=v_i.
\end{align*}
Here we have written $(x,y)=(b_0,v_0)$ and $(y,s)=(b_N,v_N)$.
	Since $t=v_0$ and $s=v_N$ are interior points, the maps $\alpha_1$ and $\beta_N$ must be surjective.
	We proceed by induction on the length of the zig-zag. If $N=1$, then the pushout of the surjections $[n]\leftarrow [k_1] \rightarrow [n]$ exists in $\Delta$ and is absolute (by \autoref{Prop:Heuts2.1}), so that $X$ turns it into a pullback
\begin{align*}
\xymatrix{
	X_l \ar[r]^{\gamma^*} \ar[d] & X_n \ar[d] \\
	X_n \ar[r]^{\beta_1^*} & X_{k_1}
}
\end{align*}
for the appropriate value of $l\leq n$.
	But, because there exists a section of $\gamma$, there is a $z\in X_l$ with $\gamma^*z=x$, meaning that $l$ must equal $n$ (otherwise we reach a contradiction with the minimality of $n$) and $\alpha_1=\beta_1=id$.
	Thus $(x,t)=(y,s)$.
	If $N>1$, factor $\beta_1$ as
\begin{align*}
[k_1] \xrightarrow{\varepsilon} [m'_1] \xrightarrow{\delta} [m_1]
\end{align*}
with $\delta$ a monomorphism and $\varepsilon$ surjective.
	Then one can apply the same argument as above to the pushout of the degeneracies
\begin{align*}
[n] \xleftarrow{\alpha_1} [k_1] \xrightarrow{\varepsilon} [m'_1]
\end{align*}
and the elements $(x,t)\in X_n\times \Delta^n$ and $(\delta^*b_1,\varepsilon_1u_1)\in X_{m'_1}\times \Delta^{m'_1}$, concluding that $\varepsilon$ must be the identity.
	Hence $\beta_1=\delta$ is a monomorphism.
	But then the pullback of $\beta_1$ and $\alpha_2$ exists in $\Delta$ (by \autoref{Prop:Heuts2.2}) and such pullbacks are easily checked to be preserved by the functor $r(\cdot)$:
\begin{align*}
\xymatrix{
	[k_1'] \ar[r]^{\theta} \ar[d]^{\eta} & [k_2] \ar[d]^{\alpha_2} \\
	[k_1] \ar[r]^{\beta_1} & [m_1]
} \hspace{1.5cm}
\xymatrix{
	r[k_1'] \ar[r]^{\theta_*} \ar[d]^{\eta_*} & r[k_2] \ar[d]^{(\alpha_2)_*} \\
	r[k_1] \ar[r]^{(\beta_1)_*} & r[m_1]
}
\end{align*}
So we can shorten the zigzag by replacing the first two spans by the single span
\begin{align*}
[n] \xleftarrow{\alpha_1\eta} [k_1'] \xrightarrow{\beta_2\theta}[m_2]
\end{align*}
and using the element $(c,w)\in X_{k'_1}\times\Delta^{k'_1}$, with $c=\eta^*a_1=\theta^*a_2$ and $w$ the unique point in $\Delta^{k'_1}$ satisfying $\eta_*w=u_1$ and $\theta_*w=u_2$.
	This complete the inductive step.
\end{proof}

	With these results, and supposing that $x\otimes t = y\otimes s$ with $x,y\in X_n$ non-degenerate, we obtain the following:
	If both $s$ and $t$ are in the interior of $\Delta^n$, then the proposition implies $(x,s)=(y,s)$.
	If one of them, say $t$, is on the boundary of $\Delta^n$, then we can write $x\otimes t = z\otimes r$ for $z$ of smaller dimension $k$ and $r$ in the interior of $\Delta^k$ uniquely.
	But then $s$ must also be on the boundary of $\Delta^n$; if it were in the interior this would contradict the uniqueness of representatives expressed bu the proposition.

	Defining $nd(X_n)$ as the subset of $X_n$ consisting of non-degenerate $n$-simplices, we conclude that the following commutative diagram
\begin{align*}
\xymatrix{
	\bigsqcup_{x\in nd(X_n)} \partial \Delta^n \ar[r] \ar[d] & |X|^{(n-1)} \ar[d]\\
	\bigsqcup_{x\in nd(X_n)} \Delta^n \ar[r] & |X|^{(n)}
}
\end{align*}
is a pushout.
	Actually, that commutative diagram is also a pullback.
	
	These results resume in the following result.
	
\begin{teorema}
\label{theo:GeomRealizationCW}
	Let $X$ be a simplicial set. The pseudotopological geometric realization of $X$, $|X|$ is a pseudotopological CW complex.
\end{teorema}

\section{Quillen Equivalence between $\cat{PsTop}$ and $\cat{sSet}$}
\label{sec.PsTopQuillenSSet}

	In the literature, it has been established a (Quillen) model category for $\cat{PsTop}$ \cite{Rieser_arXiv_2022} and found that it is Quillen equivalent to the Quillen model category in $\cat{Top}$ \cite{kapulkin_edel2023_SyntheticApproach}.
	In this section, we are going to show that that model category in $\cat{PsTop}$ is Quillen equivalent to the Quillen model category of the simplicial sets.
	This relation, plus the fibrant and cofibrant object in these categories, is going to give us information about the relations between homotopy and homology in $\cat{PsTop}$.
	In the first part of this section, we are going to cite from \cite{Hovey_1999} the results for Quillen equivalences, and then we give the formulation for the Quillen equivalence between $\cat{PsTop}$ and $\cat{sSet}$.

	In the consecutive, given a category $\cat{C}$, we denote by $\text{Mor}\cat{C}$ the \emph{arrow category}\index{Category!arrow}
	\nomenclature{$\text{Mor}\cat{C}$}{Arrow category of $\cat{C}$} of $\cat{C}$, it means, the category whose objects are all morphisms in $\cat{C}$ and whose morphisms are commutative squares.

\begin{definicion}[Functorial factorization; \cite{Hovey_1999}, 1.1.1]
	A \emph{functorial factorization}\index{factorization!functorial} is an ordered pair $(\alpha,\beta)$ of functors $\text{Mor}\cat{C}\rightarrow \text{Mor}\cat{C}$ such that $f=\beta(f)\circ\alpha(f)$ for all $f\in\text{Mor}\cat{C}$.
\end{definicion}

\begin{definicion}[Model category; 1.1.3 and 1.1.4, \cite{Hovey_1999}]
	A \emph{model category}\index{model category} on a category $\cat{C}$ is composed by three subcategories of $\text{Mor}\cat{C}$ called weak equivalences, cofibrations, and fibrations, and two functorial factorizations $(\alpha,\beta)$ and $(\gamma,\delta)$ satisfying the following properties:
\begin{enumerate}
	\item $\cat{C}$ has all the small limits and small colimits.
	
	\item (\emph{2-out-of-3}) If $f$ and $g$ are morphisms of $\cat{C}$ such that $g\circ f$ is defined and two $f$, $g$ and $g\circ f$ are weak equivalences, then so is the third.
	
	\item (\emph{Retracts})If $g$ and $f$ are morphisms of $\cat{C}$ such that $f$ is a retract of $g$ and $g$ is a weak equivalence, cofibration or fibration, then so is $f$.

	\item (\emph{Lifting}) Then trivial cofibrations have the lifting property with respect to fibrations, and cofibrations have the left lifting property with respect to trivial fibrations.
	
	\item (\emph{Factorization}) For any morphism $f$, $\alpha(f)$ is a cofibration, $\beta(f)$ is a trivial fibration, $\gamma(f)$ is a trivial cofibration, and $\delta(f)$ is fibration.
\end{enumerate}
	Where the map $f$ is said to be a \emph{trivial fibration}\index{trivial fibration} if it is a weak equivalence and a fibration; and  the map $f$ is said to be a \emph{trivial cofibration}\index{trivial cofibration} if it is a weak equivalence and a cofibration.
\end{definicion}

	The categories of our interest are $\cat{PsTop}$, $\cat{sSet}$ and $\cat{Top}$, thus the example we cite are the model categories of them:
\begin{itemize}
	\item For $\cat{Top}$ and $\cat{PsTop}$, we consider the model category with weak homotopy equivalences as weak equivalences, the fibrations are the Serre fibration and the cofibrations are retracts of the cell complexes (see \cite{Rieser_arXiv_2022}).
	
\item In $\cat{sSet}$ we have weak homotopy equivalences of the geometric realization as weak equivalences. The Kan fibrations are its fibrations, i.e., maps $f:X\rightarrow Y$ which have the right lifting property with respect to all horn inclusions
\begin{align*}
\xymatrix{
\Lambda^k[n] \ar[r] \ar[d] & X \ar[d]\\
\Delta[n] \ar[r] & Y
}
\end{align*}
	The cofibrations are the monomorphisms $f: X\rightarrow Y$ which are the level-wise injections, i.e. the morphisms of simplicial sets such that $f_n:X_n\rightarrow Y_n$ is an injection of sets for all $n\in\mathbb{N}$ (see \cite{Robert_Thesis2013}).
\end{itemize}

	Since $\cat{C}$ has small limits and colimits, there exist an initial object, usually denoted by $\varnothing$, and a final object, denoted by $*$.
	Both objects are unique up to isomorphism.
	Then, we say that $X$ is a \emph{fibrant object}\index{fibrant object} if $X\rightarrow *$ is a fibration.
	Analogously, we say that $Y$ is a \emph{cofibrant object}\index{cofibrant object} if $\varnothing\rightarrow Y$ is a cofibration.

\begin{ejemplo}
\label{exam:FibObTopCofObsSet}
	In the model categories of $\cat{Top}$ and $\cat{PsTop}$ mentioned above, all of the objects are fibrant objects and the cofibrant objects are exactly the cell complexes.
	In contrast, in the model category of $\cat{sSet}$ all of the objects are cofibrant objects.
\end{ejemplo}

	The final object is also useful to define a powerful notion in homotopy: the pointed category. We call a model category (or any category with an initial and terminal object) \emph{pointed} if the map from the initial object to the terminal object is an isomorphism.

	Given any model category $\cat{C}$, we can define the category $\cat{C_*}$ whose objects are maps $v:* \to X$, often written $(X,v)$.
	$(X,v)$ is thought as an object $X$ together with a \emph{base point}, in this case $v$.
	The morphism in $\cat{C_*}$ from $(X,v)$ to $(Y,w)$ are morphisms in $\cat{C}$ which sends $v$ to $w$.
	
	Note that $\cat{C_*}$ has arbitrary small limits and colimits because the extra structure just limit the morphisms between the objects. In addition, we can define a functor $\cat{C} \to \cat{C_*}$ such that $X \mapsto X_+\coloneqq X \sqcup *$, with base point $*$ and observe that it is left adjoint to the forgetful functor $\mfunt{U}{C_*}{C}$, which forgets the base point. Moreover, the first functor defines a faithful embedding of $\cat{C}$ into the pointed category $\cat{C_*}$ and both functors define an equivalence of categories between $\cat{C}$ and $\cat{C_*}$ if $\cat{C}$ is already a pointed model category.
	
	In short, we can note that every model category induces a pointed model category:
	
\begin{proposicion}[\cite{Hovey_1999}, 1.1.8]
\label{prop:Riehl_1.1.8}
	Suppose $\cat{C}$ is a model category. Define a map $f$ in $\cat{C_*}$ to be a cofibration (fibration, weak equivalence) if and only if $\funt{U}f$ is a cofibration (fibration, weak equivalence) in $\cat{C}$. Then $\cat{C_*}$ is a model category.
\end{proposicion}

	The factorization axiom in a model category provides us two functors.
	For every $f:\varnothing\rightarrow X$, we obtain that $\alpha(f):\varnothing\rightarrow \funt{Q}(X)$ a cofibration and $\beta(f)$ is a trivial fibration.
	Obtaining a functor $X\mapsto \funt{Q}(X)$ and a natural transformation $\funt{Q}(X)\xrightarrow{q_X} X$.
	The functor $\funt{Q}(\cdot)$ is called the \emph{cofibrant replacement functor}\index{functor!cofibrant replacement} of $\cat{C}$.
	By an analogous argument, there exists a \emph{fibrant replacement functor}\index{functor!fibrant replacement} $\funt{R}(\cdot)$ together with a natural trivial cofibration $X\rightarrow \funt{R}X$.

	The following lemma is crucial in the model categories.
	It is really useful when two model categories have adjoint functors with more structure, for example in Quillen equivalences.

\begin{lema}[Ken Brown's lemma; 1.1.12, \cite{Hovey_1999}]
\label{lem:BrownsLemma_Hovey}
	Suppose $\cat{C}$ is a model category and $\cat{D}$ is a category with a subcategory of weak equivalences which satisfies 2-out-of-3 axiom.
	Suppose $\mfunt{F}{C}{D}$ is a functor which takes trivial cofibrations between cofibrant objects to weak equivalences.
	Then $\funt{F}$ takes all weak equivalences between cofibrant objects to weak equivalences.
	Dually, if $\funt{F}$ takes trivial fibrations between fibrant objects to weak equivalences, then $\funt{F}$ takes all weak equivalences between fibrant objects to weak equivalences.
\end{lema}

	There exists an important relation between categories that preserves fibrations and trivial fibrations by one side and cofibrations and trivial cofibrations by the other side.
	Although it is not an equivalence relation, it is really useful to compare two model categories.
	
\begin{definicion}[Quillen functor; 1.3.1, \cite{Hovey_1999}]
	Suppose $\cat{C}$ and $\cat{D}$ are model categories.
\begin{enumerate}
	\item We call a functor $\mfunt{F}{C}{D}$ a \emph{left Quillen functor}\index{Quillen functor!left} if $\funt{F}$ is a left adjoint and preserves cofibrations and trivial cofibrations.
	\item We call a functor $\mfunt{U}{D}{C}$ a \emph{right Quillen functor}\index{Quillen functor!right} if $\funt{U}$ is a right adjoint and preserves fibrations and trivial fibrations.
	\item Suppose $(\funt{F},\funt{U},\varphi)$ is an adjunction from $\cat{C}$ to $\cat{D}$.
	That is $\dfunt{F}{C}{D}{U}$ are functors, and $\varphi$ is a natural isomorphism $\mcat{\funt{F}A}{B}{D}\to \mcat{A}{\funt{U}B}{C}$.
	We call $(\funt{F},\funt{U},\varphi)$ a \emph{Quillen adjunction}\index{Quillen adjuction} if $F$ is a left Quillen functor.
\end{enumerate}
\end{definicion}

	The more famous example of this kind of relation is this kind of adjunctions is the Quillen adjunction between $\cat{sSet}$ and $\cat{Top}$.
	We will extend it for the category of pseudotopological spaces, but first we need to recall the functor $\funt{Sing}$ (and define the functor $\funt{Sing_\cat{PsTop}}$), that we are going to see that is a right adjoint of the geometric realization $\funt{|\cdot|}$ (resp. $\funt{|\cdot|_\cat{PsTop}}$).

\begin{definicion}[$\funt{Sing}$ modified to $\cat{PsTop}$]
\label{def:SingPsTop}
	Given $X\in \ocat{Top}$ (resp. $X\in \ocat{PsTop}$), we define $\funt{Sing}(X)$
	\nomenclature{$\funt{Sing}(X)$}{The simplicial set having as $n$-simplices the set $\mcat{r[n]}{X}{Top}$}
	($\funt{Sing}_\cat{PsTop}(X)$)
	\nomenclature{$\funt{Sing}_\cat{PsTop}(X)$}{The simplicial set having as $n$-simplices the set $\mcat{r[n]}{X}{PsTop}$}
	to be the simplicial set having as $n$-simplices the set $\mcat{r[n]}{X}{Top}$ (resp $\mcat{r[n]}{X}{PsTop}$).
\end{definicion}

	We also need to recall that \cite{Beattie_Butzmann_2002} define the functor $\mfunt{\iota}{Top}{PsTop}$, seeing $\cat{Top}$ as a subcategory of $\cat{PsTop}$, and the topological modification of pseudotopological spaces $\mfunt{\tau}{PsTop}{Top}$.
	In the same reference we can find that $\tau\dashv\iota$, i.e. $(\tau,\iota)$ is an adjunction.

	Note that $\funt{Sing} = \funt{Sing}(\iota(\cdot))_\cat{PsTop}$ by construction.

	Since we defined the geometric realization $|\cdot|$ as a coend, we can prove that there exists an adjunction functor from $\cat{PsTop}$ to $\cat{sSet}$, the proof is completely analogous to the proof in $\cat{Top}$: The bijection $\varphi_{K,X}:\mcat{|K|_\cat{PsTop}}{X}{PsTop}\rightarrow \mcat{K}{\funt{Sing}_\cat{PsTop}(X)}{sSet}$ is given by the following compositions:

\begin{align*}
\mcat{K}{\funt{Sing}(X)}{sSet} \cong & \int_{[n]\in\Delta} \mcat{K_n}{\funt{Sing}_\cat{PsTop}(X)_n}{Set} = \int_{[n]\in\Delta} \mcat{K_n}{\mcat{r[n]}{X}{PsTop}}{Set}\\
\cong & \int_{[n]\in\Delta} \mcat{K_n\cdot r[n]}{X}{PsTop}\\
\cong & \left[ \int^{[n]\in\Delta} K_n\cdot r[n],X \right]_{\cat{PsTop}} = \mcat{|K|_\cat{PsTop}}{X}{PsTop}
\end{align*}

	In \cite{kapulkin_edel2023_SyntheticApproach} is observed that $\cat{PsTop}$ is a cofibrantly generated model category with $I\coloneqq\{S^{n-1}\rightarrow D^n\mid n\geq 0 \}$ and $J\coloneqq\{ D^n\rightarrow D^n\times I, x\mapsto (x,0) \mid n\geq 0 \}$, i.e.

\begin{definition}[\cite{Robert_Thesis2013}, 2.20], 
	Let $(\cat{M},C,F,W)$ be a model category.
	We call $\cat{M}$ \emph{cofibrantly generated}\index{model category!cofibrantly generated} if there are two sets of maps $I\subset C$ and $J\subset C\cap W$ such that:
\begin{enumerate}
	\item The domains of the maps in $I,J$ are small with respect to $I$-cell and $J$-cell respectively.
	
	\item The class of fibrations is $F=rlp(J)$ and the class of trivial fibration is $F\cap W = rlp(I)$.
\end{enumerate}
\end{definition} 

	Then, the following proposition claims that it is enough to prove that $|i|_\cat{PsTop}$ is a cofibration and $|j|_\cat{PsTop}$ is a fibration for every $i\in I$ and $j\in J$.

\begin{lema}[2.1.20, \cite{Hovey_1999}]
\label{lem:Hovey_2.1.20}
	Suppose $(\funt{F},\funt{U},\varphi):\cat{C}\to \cat{D}$ is an adjunction between model categories.
	Suppose as well that $\cat{C}$ is a cofibrantly generated model category, with generated cofibrations $I$ and generated trivial cofibrations $J$.
	Then $(\funt{F},\funt{U},\varphi)$ is a Quillen adjunction if and only if $\funt{F}f$ is a cofibration for all $i\in I$ and $\funt{F}j$ is a trivial cofibration for all $j\in J$.
\end{lema}

	In \cite{kapulkin_edel2023_SyntheticApproach}, it is noted that $\cat{PsTop}$ is a cofibrantly generated model category with $I\coloneqq\{S^{n-1}\rightarrow D^n\mid n\geq 0 \}$ and $J\coloneqq\{ D^n\rightarrow D^n\times I, x\mapsto (x,0) \mid n\geq 0 \}$.
		Then, the \autoref{lem:Hovey_2.1.20} claims that it is enough to prove that $|i|_\cat{PsTop}$ is a cofibration and $|j|_\cat{PsTop}$ is a trivial cofibration for every $i\in I$ and $j\in J$.
	The category $\cat{Top}$ is also a cofibrantly generated model category with the same $I$ and $J$, and $|\cdot|$ satisfies the conditions of \autoref{lem:Hovey_2.1.20} in this category.
	Thus, for \autoref{prop:SpheresHomeoGeoReal}, $\cat{PsTop}$ satisfies the same conditions. Concluding that:

\begin{teorema}
\label{theo:PsTopSSetQAdj}
	 The functors $|\cdot|_{PsTop}\dashv Sing_{PsTop}$ form a Quillen adjoint.	
\end{teorema}

	The following lemma express that we can define a Quillen adjunction with the left or the right Quillen functor without difference.

\begin{lema}[1.3.4, \cite{Hovey_1999}]
\label{lem:Hovey_1.3.4}
	Suppose $(\funt{F},\funt{U},\varphi):\cat{C}\to \cat{D}$ is an adjunction, and $\cat{C}$ and $\cat{D}$ are model categories.
	Then $(\funt{F},\funt{U},\varphi)$ is a Quillen adjunction if and only if $\funt{U}$ is a right Quillen functor.
\end{lema}

	By definition, it is true that the composition of two left (right) Quillen functors is a new left (right) Quillen functors.
	In the same sense, the composition of two adjunctions is an adjunction.
	Therefore, the composition of Quillen adjunctions is a Quillen adjunction as well.
	Considering these facts, then $(\tau|\cdot|_{PsTop},\funt{Sing}):\cat{sSet}\to \cat{Top}$ is a Quillen adjunction too.

	The conclusion of this section is to observe that $(|\cdot|_\cat{PsTop},\funt{Sing}_\cat{PsTop})$ is a Quillen equivalence, that means that

\begin{definicion}[Quillen Equivalence; 1.3.12 \cite{Hovey_1999}]
\label{def:QuillenEquivalence}
	A Quillen adjunction $(\funt{F},\funt{U},\varphi):\cat{C}\to \cat{D}$ is called a \emph{Quillen equivalence} if and only if, for all cofibrant objects $X\in \ocat{C}$ and fibrant objects $Y\in \ocat{D}$, a map $f:\funt{F}X\to Y$ is a weak equivalence in $\cat{D}$ if and only if $\varphi(f):X\to \funt{U}Y$ is a weak equivalence in $\cat{C}$.
\end{definicion}

	The following equivalence provides us some useful corollaries which delivers the desire outcome.

\begin{proposicion}[1.3.13, \cite{Hovey_1999}]
\label{prop:Hovey_1.3.13}
	The following are equivalent:
\begin{itemize}
	\item $(\funt{F},\funt{U},\varphi)$ is a Quillen Equivalence.

	\item The composite $X\xrightarrow{\eta}\funt{UF}X\xrightarrow{\funt{U}r_{\funt{F}X}} \funt{URF}X$ is a weak equivalence for all cofibrant $X$, and the composite $\funt{FQUX}\xrightarrow{\funt{F}q_{\funt{U}X}} \funt{FU}X \xrightarrow{\varepsilon} X$ is a weak equivalence.
\end{itemize}
\end{proposicion}

	The following two corollaries implies a Quillen equivalence between $\cat{PsTop}$ and $\cat{sSet}$.
	This is the equivalence that we were looking for and it will helps us to understand the interaction between the rational homotopy defined in $\cat{Top}$ and in $\cat{PsTop}$.

\begin{corolario}[1.3.14, \cite{Hovey_1999}]
\label{coro:Hovey_1.3.14}
	Suppose $(\funt{F},\funt{U},\varphi)$ and $(\funt{F},\funt{U'},\varphi')$ are Quillen adjunctions from $\cat{C}$ to $\cat{D}$.
	Then $(\funt{F},\funt{U},\varphi)$ is a Quillen equivalence if and only if $(\funt{F},\funt{U'},\varphi')$ is so.
	Dually, if $(\funt{F'},\funt{U},\varphi'')$ is another Quillen adjunction, then $(\funt{F},\funt{U},\varphi)$ is a Quillen equivalence if and only if $(\funt{F'},\funt{U},\varphi')$ is so.
\end{corolario}

	The above corollary tell us that $(\tau|\cdot|_\cat{PsTop},\funt{Sing}):\cat{sSet}\to \cat{Top}$ is precisely a Quillen equivalence.
	
\begin{corolario}[1.3.15, \cite{Hovey_1999}]
\label{coro:Hovey_1.3.15}
	Suppose $\mfunt{F}{C}{D}$ and $\mfunt{G}{D}{E}$ are left (resp. right) Quillen functors. Then if two out of three of $\funt{F}$, $\funt{G}$ and $\funt{GF}$ are Quillen equivalences, so is the third.
\end{corolario}

	The corollary above of the \autoref{prop:Hovey_1.3.13} have important consequences when we study the relation between weak homotopy equivalence and the singular homology. By the moment, we just write it down.

\begin{corolario}
\label{coro:UnitIsAWeakEq}
Let $(\funt{F},\funt{U},\varphi)$ be a Quillen Equivalence.
\begin{enumerate}
	\item If $\funt{F}X$ is a fibrant object, then $\eta$ is a weak equivalence.
	
	\item If $\funt{U}Y$ is a cofibrant object, then $\varepsilon$ is a weak equivalence.
\end{enumerate}
\end{corolario}

\begin{proof}
	We only prove the first claim, the second affirmation is completely analogous.
	Suppose that $\funt{F}X$ is fibrant. Then, by \autoref{prop:Hovey_1.3.13},
\begin{align*}
X\xrightarrow{\eta}\funt{UF}X\xrightarrow{\funt{U}r_{\funt{F}X}}\funt{URF}X
\end{align*}
is a weak equivalence.
	Since $\funt{F}X$ and $\funt{RF}X$ are fibrant objects, then $r_{\funt{F}X}$ is a weak equivalence too.
	Concluding that $\eta$ is also a weak equivalence by 2-out-of-3 axiom of a model category.
\end{proof}
	
	Since we have the following equality \[(\tau|\cdot|_\cat{PsTop}, \funt{Sing})=(\tau,\iota)\circ (|\cdot|_\cat{PsTop},\funt{Sing}_\cat{PsTop})\], we finally obtain via both corollaries above that:

\begin{teorema}
\label{theo:PsTopsSetQuillenEq}
	$(|\cdot|_\cat{PsTop},\funt{Sing}_\cat{PsTop})$ is a Quillen equivalence.
\end{teorema}

\begin{proof}
	We only prove the first claim, the second affirmation is completely analogous.
	Suppose that $\funt{F}X$ is fibrant. Then, by \autoref{prop:Hovey_1.3.13},
\begin{align*}
X\xrightarrow{\eta}\funt{UF}X\xrightarrow{\funt{U}r_{\funt{F}X}}\funt{URF}X
\end{align*}
is a weak equivalence.
	Since $\funt{F}X$ and $\funt{RF}X$ are fibrant objects, then $r_{\funt{F}X}$ is a weak equivalence too.
	Concluding that $\eta$ is also a weak equivalence by 2-out-of-3 axiom of a model category.
\end{proof}

	A first automatic result from \autoref{theo:PsTopsSetQuillenEq} and \autoref{coro:UnitIsAWeakEq} is that we obtain that $|\funt{Sing}_\cat{PsTop}(X)|_\cat{PsTop}$ is weak homotopy equivalent to $X$ for every $X\in \ocat{PsTop}$ because every object in $\cat{sSet}$ is a cofibrant object, in particular $\funt{Sing}_\cat{PsTop}(X)$.
	We state that claim in the following theorem.

\begin{teorema}
\label{theo:PsTopWECW}
	 Every pseudotopological space is weak homotopy equivalent to some pseudotopological CW complex.
	 Thus, every pseudotopological space is weak homotopy equivalent to some (topological) CW complex.
\end{teorema}

	We can arrive to the same conclusion with Ebel and Kapulkin's paper \cite{kapulkin_edel2023_SyntheticApproach}; however, we will need a bit more for extending rational homotopy theory from topological spaces to pseudotopological spaces.

%%%%%%%%%%%%%%%%%%%%%%%%%%%
%%%% Singular Homology
\section{Singular Homology and Hurewicz Theorem (in PsTop)}
\label{sec.SingHomHureTheo}

	At best, attempting to translate from $\cat{Top}$ to $\cat{PsTop}$ the necessary results in homological algebra to start working on Rational Homotopy would require a whole section.
	However, thanks to the results in \autoref{sec.PsTopQuillenSSet}, we can almost reduce this problem to the topological case.

	We start defining the homology that we use in the rest of the section.
	Naturally, this homology is the singular homology.

\begin{definicion}
\label{def:HomologyPsTop}
	Let $X$ be a topological space (resp. pseudotopological space) and $\Bbbk$ a ring with unit.
	For every $n\geq 0$, we define $C^n(X;\Bbbk)$
	\nomenclature{$C^n(X;\Bbbk)$}{Chains of $X$}
	as the free $\Bbbk$-module generated by $\funt{Sing}(X)_n$ (resp. $(\funt{Sing}_\cat{PsTop}(X))_n$), whose elements are called $n$-chains.
	In the same sense, for $\sigma\in \funt{Sing}(X)$ we define the boundary
\begin{align*}
\partial\sigma = \sum\limits_{i=0}^n \partial_i\sigma,
\end{align*}
and, noting that $\partial\partial\sigma=0$, we define
\begin{align*}
H_n(X;\Bbbk) \coloneqq Ker(\partial_n)/Im(\partial_{n+1})
\end{align*}
calling \emph{$n$-cycles}\index{chain cycles} to the elements of $Ker(\partial_n)$, and \emph{$n$-boundaries}\index{chain boundaries} to the elements of $Im(\partial_{n+1})$.
\end{definicion}

	Note that the long exact sequences of homology, cohomology and homotopy depends only in the structure as modules and domains, thus we have analogous results.

\begin{lema}[Compression lemma;\cite{Hatcher_2002}, 4.6]
\label{lem:CompressionLemma}
	Let $(X,A)$ be a CW pair and let $(Y,B)$ be any pair with $B\neq\varnothing$.
	For each $n$ such that $X-A$ has cells of dimension $n$, assume that $\pi_n(Y,B,y_0)=0$ for all $y_0\in B$.
	Then every map $f:(X,A)\rightarrow (Y,B)$ is homotopic rel $A$ to a map $X\rightarrow B$.

	When $n=0$ the condition that $\pi_n(Y,B,y_0)$ for all $y_0\in B$ is to be regarded as saying that $(Y,B)$ is 0-connected.
\end{lema}

\begin{proof}
	Assume inductively that $f$ has been homotoped to take the skeleton $X^{k+1}$ to $B$.
	If $\Phi$ is the chain map of a cell $e^k$ of $X-A$, the composition $f\circ \Phi$ can be homotoped $B\text{ rel }\partial D^k$ given that $\pi_k(Y,B,y_0)=0$ if $k>0$, or that $(Y,B)$ is $0$-connected if $k=0$.
	This homotopy of $f\circ \Phi$ induces a homotopy of $f$ on the quotient space $X^{k-1}\cup e^k$ \text{rel} $(X^{k-1})$.
	Doing this for all $k$-cells of $X-A$ simultaneously, and taking the constant homotopy on $A$, we obtain a homotopy of $f\mid_{X^k\cup A}$ to a map into $B$.
	By the homotopy extension property (\autoref{prop:Hatcher0.16}), this homotopy extends to a homotopy defined on all of $X$, and the induction step is completed.

	Considering any skeleton, we make the same trick with defining the homotopy in $[1-1/2^k,1-1/2^{k-1}]$ for the $k$ complex in $X$.
\end{proof}

\begin{proposicion}[\cite{Hatcher_2002}, 4.21]
\label{prop:Hatcher_4.21}
	A weak homotopy equivalence $f:X\rightarrow Y$ induces isomorphisms $H_*(f;G):H_n(X;G)\rightarrow H_n(Y,G)$ and $H^*(f;G):H^*(Y;G)\rightarrow H^*(X;G)$ for all $n$ and all coefficient groups $G$.
\end{proposicion}

\begin{proof}
	We know that every strong homotopy equivalence induces isomorphisms in singular homology groups and singular cohomology groups (the proof is the same as \cite{Hatcher_2002}, 2.10 without any change), then we can replace $Y$ by the mapping cylinder $M_f$.
	Therefore, observing the long exact sequences of homotopy, homology and cohomology, it is enough to prove that:

	If $(Z,X)$ is an $n$-connected pair of path-connected spaces, then  $H_i(Z,X;G)=0$ and $H^i(Z,X;G)=0$ for all $i\leq n$ and all $G$.

	We prove that affirmation. Let $\alpha\coloneqq\sum_j n_j\sigma_j \in H_k(Z,X;G)$, for singular $k$-simplices $\sigma_j:\Delta^k\rightarrow Z$.
	Build a $\Delta$-complex $K$ from a disjoint union of $k$-simplices, one for each $\delta_j$, by identifying all $(k-1)$ faces of these $k$-simplices for which the corresponding restriction of the $\sigma_j$s are equal.
	Thus the $\sigma_j$s induce a map $\sigma:K\rightarrow Z$.
	Since $\alpha$ is a relative cycle, $\partial\alpha$ is a chain in $X$.
	Let $L\subset K$ be the subcomplex consisting of $(k-1)$-simplices corresponding to the singular $(k-1)$ simplices in $\partial\alpha$, so $\sigma(L)\subset X$.
	The chain $\alpha$ is the image under the chain map $\sigma_\#$ of a chain $\tilde{\alpha}$ in $K$, with $\partial\tilde{\alpha}$ a chain in $L$.
	In relative homology we then have $\sigma_*[\tilde{\alpha}]=[\alpha]$.
	If we assume $\pi_i(Z,X)=0$ for $i\leq k$, then $\sigma$ is homotopic relative to $L$ to a map with image in $X$, by \autoref{lem:CompressionLemma}.
	Hence $\sigma_*[\tilde{\alpha}]$ is in the image of the map $H_k(X,X;G)\rightarrow H_k(Z,X;G)$ and since $H(X,X;G)=0$ we conclude that $[\alpha]=\sigma_*[\tilde{\alpha}]=0$.
	This proves the result for homology, and the result for cohomology then follows by the universal coefficient theorem.

	To conclude the proof, observe that $(M_f,X)$ is $n$-connected for every $n$ because they homotopy groups are isomorphic.
\end{proof}

	\autoref{prop:Hatcher_4.21} and \autoref{coro:UnitIsAWeakEq} allow us to conclude the following theorem which is key in the proofs of rational homotopy theory for $\cat{PsTop}$.

\begin{teorema}
\label{theo:PsTopTopWHE}
	Let $f:X\rightarrow Y$ be a morphism in $\cat{PsTop}$, then there exists a topological morphism $f^\circ:X^\circ\rightarrow Y^\circ$ such that $f$ is a weak homotopy equivalence if and only if $f^\circ$ is a weak homotopy equivalence.
	The same affirmation is true for replacing weak homotopy equivalence condition for quasi-isomorphisms.
\end{teorema}

\begin{proof}
	Let $f:X\rightarrow Y$ be a morphism in $\cat{PsTop}$.
	Since every object is cofibrant in $\cat{sSet}$ and every object is fibrant in $\cat{PsTop}$, by \autoref{coro:UnitIsAWeakEq} we obtain that the image of the natural transformation $\varepsilon_{\cat{PsTop},\cat{sSet}}:|\funt{Sing}(\cdot)_\cat{PsTop}|_\cat{PsTop}\Rightarrow \Id{\cat{PsTop}}$ is a weak homotopy equivalence.
	Analogously, since every object is fibrant in $\cat{Top}$ and $|\funt{Sing}(\cdot)_\cat{PsTop}|_\cat{PsTop}$ is a CW complex, and then cofibrant, by \autoref{coro:UnitIsAWeakEq} we obtain that the image of the natural transformation $\eta_{\cat{PsTop},\cat{Top}}:\Id{\cat{PsTop}}\Rightarrow \iota\tau(\cdot)$ is a weak homotopy equivalence in $\varepsilon(X)$, i.e., we have the following commutative diagram
\begin{align*}
\xymatrixcolsep{10mm}\xymatrix{
X \ar[d]_{f} & \varepsilon X \ar[r]^{\eta_{\varepsilon X}}_{\text{w.e.}} \ar[l]_{\varepsilon_X}^{\text{w.e.}} \ar[d]_{\varepsilon f} & \eta\varepsilon X \ar[d]_{\eta\varepsilon f} \\
Y  & \varepsilon Y \ar[r]_{\eta_{\varepsilon Y}}^{\text{w.e.}} \ar[l]^{\varepsilon_Y}_{\text{w.e.}} & \eta\varepsilon Y
}
\end{align*}
	Define $f^\circ \coloneqq \eta\varepsilon f$, $X^\circ \coloneqq \eta\varepsilon X$ and $Y^\circ \coloneqq \eta\varepsilon Y$.
	Suppose that $f$ is a weak homotopy equivalence (resp. a quasi-isomorphism), then
\begin{align*}
\xymatrixcolsep{15mm}\xymatrix{
\pi_n(X) \ar[d]_{\pi_n(f)} & \pi_n(\varepsilon X) \ar[r]^{\pi_n(\eta_{\varepsilon X})}_{\cong} \ar[l]_{\pi_n(\varepsilon_X)}^{\cong} \ar[d]_{\pi_n(\varepsilon f)} & \pi_n(X^\circ) \ar[d]_{\pi_n(f^\circ)} \\
\pi_n(Y)  & \pi_n(\varepsilon Y) \ar[r]_{\pi_n(\eta_{\varepsilon Y})}^{\cong} \ar[l]^{\pi_n(\varepsilon_Y)}_{\cong} & \pi_n(Y^\circ)
}
\end{align*}
\begin{align*}
\left(resp.~~~~ \vcenter{\xymatrixcolsep{15mm}\xymatrix{
H_n(X;\Bbbk) \ar[d]_{H_n(f;\Bbbk)} & H_n(\varepsilon X;\Bbbk) \ar[r]^{H_n(\eta_{\varepsilon X};\Bbbk)}_{\cong} \ar[l]_{H_n(\varepsilon_X;\Bbbk)}^{\cong} \ar[d]_{H_n(\varepsilon f;\Bbbk)} & H_n(X^\circ;\Bbbk) \ar[d]_{H_n(f^\circ;\Bbbk)} \\
H_n(Y;\Bbbk)  & H_n(\varepsilon Y;\Bbbk) \ar[r]_{H_n(\eta_{\varepsilon Y};\Bbbk)}^{\cong} \ar[l]^{H_n(\varepsilon_Y;\Bbbk)}_{\cong} & H_n(Y^\circ;\Bbbk)
}}\right)
\end{align*}
for every $n\geq 0$.
	Thus $\pi_n(f^\circ)$ (resp. $H_n(f^\circ;\Bbbk)$) is an isomorphism because is a bijection and a homomorphism, concluding that $f^\circ$ is a weak homotopy equivalence (resp. a quasi-isomorphism).

	On the other hand, suppose that $f^\circ$ is a weak homotopy equivalence (resp. a quasi-isomorphism).
	It is completely analogous to prove that that $f$ is as well.
\end{proof}

With the previous result, we obtain as a corollary the analogous for pseudotopological spaces of the following important result in topological spaces. 

\begin{corolario}[Whitehead;\cite{FelixHalperinThomas_2001_RatHom}, 4.19,  modified for $\cat{PsTop}$]
\label{theo:FHT_4.19}
	Let $(X,x_0)$ be a pointed pseudotopological space for which $\pi_i(X,x_0)=0$, $i\leq r$.
\begin{enumerate}
	\item If $r=0$, then the homomorphism
\begin{align*}
	hur_X:\pi_{1}(X,x_0)\rightarrow H_{1}(X;\mathbb{Z})
\end{align*}
is surjective and its kernel is the subgroup of $\pi_1$ generated by the commutators $\alpha\beta\alpha^{-1}\beta^{-1}$.

	\item If $r\geq 1$, then $H_i(X;\mathbb{Z})=0$, $1<i\leq r$ and
\begin{align*}
hur_X:\pi_{r+1}(X,x_0)\rightarrow H_{r+1}(X;\mathbb{Z})
\end{align*}
is an isomorphism.
\end{enumerate}
\end{corolario}

\begin{proof}
	Let $(X,x_0)$ be a pointed pseudotopological space.
	Given the construction in the proof of \autoref{theo:PsTopTopWHE}, we have that $\eta_X$ and $\varepsilon_X$ provide isomorphisms between the homotopy groups of $X$ and $X^\circ$ and the homological groups of $X^\circ$, with $X^\circ$ a topological space.
	Since $X^\circ$ is already a topological space, we know that this space satisfies properties (1) and (2).
	Thus, the same is satisfied for $H_i^{-1}(\eta_X\varepsilon_X)hur_{(X^\circ)}\pi_i^{-1}(\eta_X\varepsilon_X)$.
\end{proof}

	Also using \autoref{theo:PsTopTopWHE}, we obtain the following two results which makes the definition of rational pseudotopological spaces.

\begin{corolario}[Whitehead; \cite{FelixHalperinThomas_2001_RatHom}, 8.6, modified for $\cat{PsTop}$]
\label{theo:FHT_8.6}
	Suppose that $\Bbbk$ is a subring of $\mathbb{Q}$ and $g:X\rightarrow Y$ is a morphism in $\cat{PsTop}$ of simply connected spaces.
	Then the following assertions are equivalent:
\begin{enumerate}
\item $\pi_*(g)\otimes\Bbbk$ is an isomorphism.
\item $H_*(g;\Bbbk)$ is an isomorphism.
%\item $H_*(\Omega g;\Bbbk)$ is an isomorphism.
\end{enumerate}
\end{corolario}

	Now we start to discuss our subject of interest in this section: the rational homotopy.
	First we begin by introducing $\mathcal{P}$-local spaces, and then we go directly to the definition of the $A_{PL}$ functor and the minimal models.
	
	Fix a set $\mathcal{P}$ of prime numbers and define $\mathcal{I}$ as the set of integers which are relatively prime to the elements of $\mathcal{P}$.
	We said that an abelian group $G$ is \emph{$\mathcal{P}$-local}\index{$\mathcal{P}$-local} if the multiplication by $k$ in $G$ is an isomorphism for every $k\in\mathcal{I}$.
	For $k\in\mathcal{I}$, the rational numbers $m/k$ are a subring $\Bbbk\subset \mathbb{Q}$, thus $G$ is $\mathcal{P}$-local if and only if it is a $\Bbbk$-module.
	By convention, when $\mathcal{P}=\varnothing$ we take $\Bbbk= \mathbb{Q}$.
	
	Observe that for any abelian group $G$, $G\otimes_\mathbb{Z}\Bbbk$ is $\mathcal{P}$-local.
	Moreover, if $G$ is $\mathcal{P}$-local then $G\rightarrow G\otimes_\mathbb{Z}\Bbbk$ is an isomorphism.
	For general any abelian group $G$, the morphism \[G\rightarrow G\otimes_\mathbb{Z} \Bbbk,\ \gamma\mapsto \gamma\otimes 1\] is called the $\mathcal{P}$-localization of $G$. 
	
\begin{definicion}[\cite{FelixHalperinThomas_2001_RatHom}, page 102]
\label{def:RationalSpace}
	A simply connected pseudotopological space $X$ is a \emph{$\mathcal{P}$-local space} if $\pi_*(X)$ is $\mathcal{P}$-local.
	When $\Bbbk= \mathbb{Q}$ is called a rational space.
\end{definicion}

	Direct from this definition, we can note that a pseudotopological spaces is a $\mathcal{P}$-local spaces if and only if the topological modification of the $|\funt{Sing}_\cat{PsTop}(\cdot)|_\cat{PsTop}$ does so. 
	
	The following result is a direct corollary from \autoref{theo:PsTopTopWHE} and provides us of equivalent definitions of $\mathcal{P}$-local spaces.

\begin{corolario}[\cite{FelixHalperinThomas_2001_RatHom}, 9.3 modified for $\cat{PsTop}$]
\label{theo:FHT_9.3}
	Let $X$ be a simply connected pseudotopological space.
	Then the following conditions are equivalent:
\begin{enumerate}
	\item $\pi_*(X)\otimes\Bbbk$ is $\mathcal{P}$-local.
	
	\item $H_*(X,pt;\Bbbk)$ is $\mathcal{P}$-local.

	%\item $H_*(\Omega X, pt;\Bbbk)$ is $\mathcal{P}$-local.
\end{enumerate}
\end{corolario}

	Now that we know what a rational space is, we define how we can go from a pseudotopological space to a rational pseudotopological space.
	In the subsequent section we observe that, for every simply connected space, we can always do that step.
	By a modern algebra result, if $\varphi: X \to Y$ is a continuous map between simply connected spaces and $Y$ is already $\mathcal{P}$-local, then $\pi_*(\varphi)$ extends uniquely to a morphism $\pi_*(X)\otimes_\Z\Bbbk \to \pi_*(Y)$.
	This result allows the next definition.

\begin{definicion}[\cite{FelixHalperinThomas_2001_RatHom}, page 108]
\label{def:Razionalization}
	A \emph{$\mathcal{P}$-localization}\index{$\mathcal{P}$-localization} of a simply connected space $X$ is a map $\varphi: X \to X_\mathcal{P}$ to a simply connected $\mathcal{P}$-local space $X_\mathcal{P}$ such that $\varphi$ induces an isomorphism \[\pi_*(X)\otimes_\Z \Bbbk \xrightarrow{\cong} \pi_*(X_\mathcal{P}).\]
	When $\Bbbk=\Q$ this is called a \emph{rationalization} and is denoted by $X\to X_\Q$
\end{definicion}

	We have the following equivalence of a $\mathcal{P}$-localization in terms of the functor $H_*(\cdot)$ instead of the homotopy groups.

\begin{teorema}[\cite{FelixHalperinThomas_2001_RatHom}, 9.6]
\label{theo:FHT_9.6}
	A continuous map $\varphi: X \to Y$ between simply connected spaces is a $\mathcal{P}$-localization if and only if $Y$ is $\mathcal{P}$-local and $H_*(\varphi,\Bbbk)$ is an isomorphism.
\end{teorema}

	Considering $\Bbbk = \Q$, \autoref{theo:FHT_8.6} implies that, for every continuous map $f:X\to Y$ between simply connected spaces, the following conditions are equivalent
	\begin{itemize}
	\item $\pi_*(f)\otimes\Q$ is an isomorphism.
	\item $H_*(f;\Q)$ is an isomorphism.
	\end{itemize}
	In this case $f$ is called a \emph{rational homotopy equivalence}\index{rational homotopy equivalence}.

\begin{definicion}[\cite{FelixHalperinThomas_2001_RatHom}, page 111]
\label{def:RationalHomotopyType}
	The weak homotopy type of $X_\mathbb{Q}$ is the \emph{rational homotopy type} of $X$.
\end{definicion}

	The following proposition gives an equivalent description of a rational homotopy type which only depends on rational homotopy equivalences, and then involves the homotopy and the homology.
	This is going to be useful when we observe that two spaces with the same rational homotopy type have the same $A_{PL}(\cdot)$ in \autoref{theo:FHT_10.9}.

\begin{proposicion}[\cite{FelixHalperinThomas_2001_RatHom}, 9.8]
\label{theo:FHT_9.8}
	Simply connected spaces $X$ and $Y$ have the same rational homotopy type if and only if there is a chain of rational homotopy equivalences
	\begin{align*}
	X \leftarrow Z(0) \rightarrow \cdots \leftarrow Z(k)\rightarrow Y.
	\end{align*}
\end{proposicion}

	For the previous proposition and the construction of $X^\circ$ in \autoref{theo:PsTopTopWHE}, we obtain that $X$, $\varepsilon X$ and $X^\circ$ have the same rational homotopy type because $\varepsilon_X$ and $\eta_{\varepsilon X}$ have weak homotopy equivalence and then they are a zigzag of rational homotopy equivalences.
	Furthermore, we make the following remark.

\begin{observacion}
	Let $X$ and $Y$ be simply connected pseudotopological spaces such that $f^\circ: X^\circ \to Y^\circ$ is a rational homotopy equivalence, then $f:X\to Y$ is a rational homotopy equivalence too.
	In addition, if $X^\circ$ and $Y^\circ$ have the same rational homotopy type, then $X$ and $Y$ do.
\end{observacion}

\section{Rational Homotopy in $\cat{PsTop}$}
\label{sec.RatHom}

	In this section, we aim to prove that we can compute the rational homotopy groups of a pseudotopological space in the same way that they are computed for topological spaces.
	Most of these constructions are not made in the category of topological spaces, but in the category of simplicial sets $\cat{sSet}$ or in the category of differential graded commutative algebras $\cat{DGCA}$.
	
	In order to simplify notation, and following the notation of \cite{FelixHalperinThomas_2001_RatHom}, in this section we denote by $S_*(\cdot)$
	\nomenclature{$S_*(\cdot)$}{Short for $\funt{Sing}_\cat{PsTop}(\cdot)$}
	the functor $\funt{Sing}_\cat{PsTop}(\cdot)$ of \autoref{def:HomologyPsTop}, which is the simplicial set having as $n$-simplices the set $\mcat{r[n]}{\cdot}{PsTop}$.
	
	We briefly recall that a \emph{differential graded algebra}\index{differential graded algebra} is a graded algebra $R$ together with a differential in $R$ that is a derivation $d$ of grade $k$, i.e., a linear map $d:R\to R$ of degree $-1$ such that $d(xy)=(dx)y+(-1)^{k deg(x)}xdy$.
	We can define the \emph{homology of the algebra}\index{homology!of an algebra} $H(R,d)$ as the graded algebra $H(R,d)= ker(d)/Im(d)$.
	An \emph{augmentation}\index{augmentation} is a morphism $\epsilon: (R,d)\to \Bbbk$.
	
	A morphism of differential graded algebras $f:(R,d)\to (S,d)$ is a morphism of graded algebras satisfying $fd=df$, inducing a morphism $H(f):H(R)\to H(S)$ of graded algebras.
	This forms a category, which is denoted by $\cat{DGA}$.
	\nomenclature{$\cat{DGA}$}{The category of differential graded algebras}
	In addition, if the objects satisfies that $xy=(-1)^{deg(x)deg(y)}yx$ for every $x,y\in A$, we say that the objects are \emph{commutative} and we denote the category as $\cat{DGCA}$.
	\nomenclature{$\cat{DGCA}$}{The category of differential graded commutative algebras}
	
	In $\cat{DGA}$ we use the notation $\xrightarrow{\simeq}$ to denote a \emph{quasi-isomorphism}\index{quasi-isomorphism}, i.e., a map $f:(R,d)\to (S,d)$ of differential graded algebras which induces an isomorphism $H(f):H(R)\to H(S)$ of graded algebras. 
	
	In the same reference, \cite{FelixHalperinThomas_2001_RatHom}, we have the next definitions:
	let $V$ be any free graded module, the \emph{tensor algebra}\index{tensor algebra} $TV$
	\nomenclature{$TV$}{The tensor algebra of a free graded algebra $V$}
	is defined by
	\begin{equation}
	TV\coloneqq \oplus_{q=0}^\infty T^q V,\ T^q V\coloneqq \underbrace{V\otimes \cdots \otimes V}_q
	\end{equation}
	and the \emph{free commutative grade algebra on $V$} is denoted by $\Lambda V$
	\nomenclature{$\Lambda V$}{Free commutative graded algebra on $V$}
	and defined by
	\[\Lambda V \coloneqq TV/I\]
	with $I$ the ideal generated by
	\begin{equation}
	v\otimes w - (-1)^{deg(v)deg(w)}w\otimes v,\ v,w\in V.
	\end{equation}

	We may use the following notations and results described in \cite{FelixHalperinThomas_2001_RatHom} without further reference:
\begin{itemize}
	\item $\Lambda V = \text{symmetric algebra}(V^\text{even})\otimes \text{exterior algebra}(V^\text{odd})$.
	
	\item If $\{v_\alpha\}$ or $\{v_1, v_2\ldots \}$ is a basis for $V$ we write $\Lambda(\{v_\alpha\})$ or $\Lambda (v_1,v_2,\ldots)$ for $\Lambda V$.
	
	\item $\Lambda^q V$ is the linear span of elements of the form $v_1 \wedge\cdots\wedge v_q$, $v_i\in V$. Elements in $\Lambda^q V$ have word-length $q$.
	
	\item $\Lambda V = \oplus_q \Lambda^q (V)$ and we write $\Lambda^{\geq q} (V) = \oplus_{i\geq q}\Lambda^i (V)$ and $\Lambda^+ (V) = \Lambda^{\geq 1} (V)$.
	
	\item If $V=\oplus_\lambda V_\lambda$ then $\Lambda V = \otimes \Lambda V_\lambda$.
	
	\item Any linear map of degree zero from $V$ to a commutative graded algebra $A$ extends to a unique graded algebra morphism $\Lambda V\to A$.
	
	\item Any linear map of degree $k\in\Z$ from $V$ to $\Lambda V$ extends to a unique derivation of degree $k$ in $\Lambda V$.
\end{itemize}
	
	We define the differential graded commutative algebra $A_{PL} = \{ (A_{PL})_n \}_{n\geq 0}$ by: The cochain algebra $(A_{PL})_n$
	\nomenclature{$A_{PL}$}{Piecewise linear algebra}
	is given by
	\begin{equation}
	(A_{PL})_n = \frac{\Lambda(t_0,\ldots,t_n,y_0,\ldots,y_n)}{(\sum t_i -1, \sum y_i)}
	\end{equation}
	such that $dt_i = y_i$ and $dy_i=0$.
	
	The face degeneracy morphisms are the unique cochain algebra morphisms \[ \partial_i:(A_{PL})_{n+1} \to (A_{PL})_{n} \text{ and } s_j: (A_{PL})_{n} \to (A_{PL})_{n+1} \] satisfying \[ \partial_i(t_k) = \begin{cases} t_k & k<i\\ 0 & k=i\\ t_{k-1} & k>i \end{cases} \text{ and } s_j(t_k) = \begin{cases} t_k & k<j\\ t_k + t_{k+1} & k=j\\ t_{k+1} & k>j \end{cases} \]
	
	In general, for every simplicial set $K$ and for every simplicial cochain complex or a simplicial cochain complex $A = \{ A_n \}_{n\geq 0}$ we can define the ordinary cochain complex (or cochain complex $A(K) = \{ A^p(K) \}_{p\geq 0}$ where $A^p(K)$ is the set of simplicial set morphisms from $K$ to $A_p$ (for details of this definition, see \cite{FelixHalperinThomas_2001_RatHom} page 118).
	In particular, for pseudotopological spaces $X$ and continuous maps $f$ we apply this construction to the simplicial set $S_*(X)$ and to $S_*(f)$.
	This define a contravariant functor from $\cat{PsTop}$ to commutative cochain algebras, which will be denoted \[A_{PL}(X) \text{ and } A_{PL}(f).\] 

	On the other hand, in \autoref{def:HomologyPsTop} we denote $C^*(X;\Bbbk)$ the singular cochain algebra.
	As the construction only depends on $S_*(X)$, it generalizes to any simplicial set $K$ to give the cochain algebra $C^*(K,\Bbbk)$, where, in particular we will denote $C^*(S_*(X))$ by $C^*(X)$.
	
	In addition, we can define the simplicial cochain algebra $C_{PL}$ using the simplicial sets $\Delta[n] \subset S_*(\Delta^n)$ as $C_{PL} = \{ (C_{PL})_n \}_{n\geq 0}$ by $(C_{PL})_n$ is the cochain algebra $C^*(\Delta[n])$, and the face and degeneracy morphisms are the $C^*([d_i])$ and $C^*([s_j])$, with $d_i$ and $s_j$ defined as \autoref{def:SimplicialSets}.
	
	Finally, let
	\[C_{PL}\otimes A_{PL} \coloneqq \{ (C_{PL})_n\otimes (A_{PL})_n; \partial_i\otimes \partial_i; s_j\otimes s_j \}\]
	be the tensor product simplicial cochain algebra.
	Morphisms \[ C_{PL} \rightarrow C_{PL}\otimes A_{PL} \leftarrow A_{PL} \] are defined by $\gamma\mapsto \gamma \otimes 1$ and $\Phi \mapsto 1\otimes \Phi$.
	
	With those constructions, we recall the following sentence, that remains true for the category of pseudotopological spaces because only depends on simplicial sets.
	
\begin{teorema}[\cite{FelixHalperinThomas_2001_RatHom}, 10.9]
\label{theo:FHT_10.9}
	Let $K$ be a simplicial set. Then
	\begin{enumerate}
	\item There is a natural isomorphism $C_{PL}(K) \xrightarrow{\cong} C^*(K)$ of cochain algebras.
	
	\item The natural morphisms of cochain algebras, \[C_{PL}(K) \rightarrow (C_{PL}\otimes A_{PL})(K) \leftarrow A_{PL}(K)\] are quasi-isomorphisms.
	\end{enumerate}
\end{teorema}

	It is a corollary to say that, for any pseudotopological space $X$ there exist natural quasi-isomorphisms between the cochain algebras
	\[C_{PL}(X) \xrightarrow{\simeq} (C_{PL}\otimes A_{PL})(X) \xleftarrow{\simeq} A_{PL}(X)\]
	and this give an isomorphism $H^*(X)\cong H(A_{PL}(X))$.
	This allows to identify the singular homology groups with the the homology in the differential graded commutative algebra $A_{PL}(X)$, and thus $A_{PL}(\cdot)$ defines a map from rational homotopy types to weak equivalence classes of commutative cochain algebra by \autoref{theo:FHT_9.8}.

	Now we introduce the Sullivan algebra and the notion of model for a pseudotopological space.
	The minimal Sullivan model for a pseudotopological spaces $X$ is the central piece in the rational homotopy theory to simplify the computations of the rational homotopy groups.
	
\begin{definicion}[\cite{FelixHalperinThomas_2001_RatHom}, page 138]
\label{def:SullivanAlgebras}
	A \emph{Sullivan algebra}\index{Sullivan algebra} is a commutative cochain algebra of the form $(\Lambda V, d)$, which satisfies the following conditions

\begin{itemize}
	\item $V= \{ V^p \}_{p\geq 1}$ and $\Lambda V$ denotes the free graded commutative algebra on $V$;.
	
	\item (\emph{Nilpotence condition on $d$}) $V = \cup_{k=0}^\infty V(k)$, where $V(0)\subset V(1) \subset \ldots$ is an increasing sequence of graded subspaces such that $d(0)=0$ in $V(0)$ and $d:V(k)\to \Lambda V(k-1)$ for every $k\geq 1$.
\end{itemize}
\end{definicion}
	
\begin{definicion}[\cite{FelixHalperinThomas_2001_RatHom}, definition 1 in page 138]
\label{def:SullivanModels}
	A \emph{Sullivan model}\index{Sullivan model} for a commutative cochain algebra $(A,d)$ is a quasi-isomorphism \[m:(\Lambda V,d) \to (A,d)\] from a Sullivan algebra $(\Lambda V, d)$.
	
	If $X$ is a path connected topological space then a Sullivan model for $A_{PL}(X)$, \[m: (\Lambda V,d) \xrightarrow{\simeq} A_{PL}(X)\] is called a \emph{Sullivan model for $X$}.
	
	A Sullivan algebra (or model), $(\Lambda V,d)$ is called \emph{minimal} if \[Im(d)\subset \Lambda^+ V \cdot \Lambda^+ V.\]
\end{definicion}

\begin{observacion}
\label{obs:QIsomModel}
	Given $f:X\to Y$ a quasi-isomorphism of two pseudotopological spaces and $m:(\Lambda V,d) \to A_{PL}(X)$ a Sullivan model for $X$.
	Via the identification of $H(A_{PL}(f))$ and $H^*(f)$, we can observe that $A_{PL}(f)\circ m$ is also a quasi-isomorphism and then it is a Sullivan model for $Y$.
\end{observacion}

	The first big question is about the existence of a minimal Sullivan model for a pseudotopological space.
	The combination of two results tells us that there always exists a minimal Sullivan model for a specific, but pretty general, kind of pseudotopological space. 

\begin{proposicion}[\cite{FelixHalperinThomas_2001_RatHom}, 12.1 and 12.2]
\label{prop:FHT_12.1}
	Any $(A,d)$ commutative cochain algebra $(A,d)$ satisfying that $H^0(A)=\Bbbk$ has a Sullivan model \[m:(\Lambda V,d) \xrightarrow{\simeq} (A,d).\]
	In addition, if $H^1(A)=0$, then it has a minimal Sullivan model.
	Further, if $X$ satisfies that each $H_i(X;\Q)$ is finite dimensional, then $X$ has a minimal Sullivan model \[m:(\Lambda V,d) \xrightarrow{\simeq} A_{PL}(X)\] such that $V=\{V^i\}_{i\geq 2}$ and each $V^i$ is finite dimensional.
\end{proposicion}

	The following result says when we can lift a morphism from a Sullivan algebra $(\Lambda V,d)$.
	We use this result later to compare a minimal Sullivan model of a pseudotopological space $X$ and a minimal Sullivan model of $X^\circ$, where $X^\circ$ denote the topological space defined in \autoref{theo:PsTopTopWHE}.

\begin{lema}[\cite{FelixHalperinThomas_2001_RatHom}, 12.4]
\label{lem:LiftingLemma}
	Let $\psi: (\Lambda V,d) \to (C,d)$ and $\eta: (A,d) \to (C,d)$ be morphisms of commutative cochain algebras such that $(\Lambda V,d)$ is a Sullivan algebra and $\eta$ is a surjective quasi-isomorphisms.
	There is a morphism $\varphi: (\Lambda V,d) \to (A,d)$ such that $\eta\varphi = \psi$.
\end{lema}

	In \cite{FelixHalperinThomas_2001_RatHom}, on page 139, it is defined a homotopy for morphisms whose domain is the same Sullivan algebra:
	Two morphisms $\varphi_0,\varphi_1:(\Lambda V,d) \to (A,d)$ from a Sullivan algebra are \emph{homotopic} if there exists a morphism \[\Phi: (\Lambda V,d)\rightarrow (A,d)\otimes \Lambda(t,dt)\] such that $(id\cdot \epsilon_i)\Phi = \varphi_i$, $i\in \{0,1\}$, with the augmentations $\epsilon_0,\epsilon_1:\Lambda(t,dt)\to \Bbbk$ such that $\epsilon_0(t)=0$ and $\epsilon_1(t)=1$.
	The homotopy relation is an equivalence relation in the set of morphisms from a Sullivan algebra.
	
	Given continuous maps $f_0,f_1: X\to Y$ and a morphism $\psi: (\Lambda V,d) \to A_{PL}(Y)$ from a Sullivan algebra $(\Lambda V,d)$, we have the following result:
	
\begin{proposicion}[\cite{FelixHalperinThomas_2001_RatHom}, 12.6]
\label{prop:FHT_12.6}
	If $f_0\simeq f_1: X\to Y$ then $A_{PL}(f_0)\psi \simeq A_{PL}(f_1)\psi: (\Lambda V,d) \to A_{PL}(X)$.
\end{proposicion}

\begin{definicion}[\cite{FelixHalperinThomas_2001_RatHom}, page 152]
\label{def:LinearPart}
	The \emph{linear part}\index{linear part} of a morphism $\varphi: (\Lambda V,d)\to (\Lambda W,d)$ between Sullivan algebras is the linear map \[Q\varphi: V\to W\] defined by: $\varphi v - Q\varphi v \in \Lambda^{\geq 2} W$, $v\in V$.
\end{definicion}

	If we have two morphisms $\varphi: (\Lambda V,d) \to (\Lambda W,d)$ and $\psi: (\Lambda W, d) \to (\Lambda U,d)$, it follows that $Q(\psi \circ \varphi) = Q(\psi)\circ Q(\varphi)$.
	This follows for the uniqueness of $Q(\psi\circ\varphi)$ and the definition of morphism of algebras; remind that every Sullivan algebra comes from a graded vector space $V$ which least degree $1$.
	
	Weaker than \autoref{lem:LiftingLemma}, we have the following result that observes the behavior of the homotopy classes of morphisms.

\begin{proposicion}[\cite{FelixHalperinThomas_2001_RatHom}, 12.9]
\label{prop:FHT_12.9}
	Let $\psi: (\Lambda V,d) \to (C,d)$ and $\eta: (A,d) \to (C,d)$ be morphisms of commutative cochain algebras such that $(\Lambda V,d)$ is a Sullivan algebra and $\eta$ is a quasi-isomorphisms.
	There is a unique homotopy class of morphisms $\varphi: (\Lambda V,d) \ to (A,d)$ such that $\eta \varphi \simeq \psi$.
	Thus \[\eta_\#: [(\Lambda V,d),(A,d)] \xrightarrow{\cong} [(\Lambda V,d),(C,d)]\] is a bijection.
\end{proposicion}

	Suppose $\alpha: (A,d)\to (A',d)$ is an arbitrary morphism of commutative algebras that satisfy $H^0=\Bbbk$.
	Let $m: (\Lambda V,d) \to (A,d)$ and $m': (\Lambda V',d)\to (A',d)$ be Sullivan models, there is a unique homotopy class of morphisms $\varphi:(\Lambda V,d) \to (\Lambda V',d)$ such that $m'\varphi \sim \alpha m$.
	A morphism $\varphi$ which satisfies that conditions is called a \emph{Sullivan representative for $\alpha$}.
	If $f:X\to Y$ is a continuous map then a Sullivan representative of $A_{PL}(f)$ is called a \emph{Sullivan representative of $f$}. 

\begin{proposicion}[\cite{FelixHalperinThomas_2001_RatHom}, 12.10]\
\label{prop:FHT_12.10}
	\begin{enumerate}
	\item A quasi-isomorphism between minimal Sullivan algebras is an isomorphism if both cohomology algebras vanish in degree one.
	\item If $(A,d)$ is a commutative cochain algebra and $H^0(A)=\Bbbk$, $H^1(A)=0$, then the minimal models of $(A,d)$ are all isomorphic.
	\end{enumerate}
\end{proposicion}

	In some sense, this proposition is extended for spaces that does not have a trivial first cohomology group.
	The extension of this section have been long enough and we are just working with the trivial first cohomology group, but we invite the reader to check the Chapter 14 of the references of rational homotopy \cite{FelixHalperinThomas_2001_RatHom}.

	Let $f:X\to Y$ be a quasi-isomorphism of simply connected pseudotopological spaces, and let $m_X:(\Lambda V,d)\to A_{PL}(X)$ and $m_Y: (\Lambda W,d) \to A_{PL}(Y)$ be minimal Sullivan models of $X$ and $Y$, respectively.
	\autoref{prop:FHT_12.10} joint with \autoref{obs:QIsomModel} implies that $(\Lambda V,d) \cong (\Lambda W,d)$.
	This is because $(\Lambda V,d)$ is a minimal model for both $X$ and $Y$.

	When $\Bbbk = \Q$, the paragraph above can be stated in the following claim. 
	
\begin{teorema}
\label{theo:RatHomTypMinSullAlg}
	Let $X$ and $Y$ be a pair of simply connected pseudotopological spaces with the same rational homotopy type.
	Then they have the same minimal Sullivan model up to isomorphism.
\end{teorema}

	Considering a minimal Sullivan model $m_X: (\Lambda V, d)\to A_{PL}(X)$, now we are capable to define a natural pairing $\langle\cdot ;\cdot \rangle$ between $V$ and $\pi_*(X)$ which only depends on the choice of $m_X$.
	
	For this construction, we use the fact that the following are minimal models: $m_k:(\Lambda(e),0)\to A_{PL}(S^k)$, when $k$ is odd, and $m_k:(\Lambda (e,e'),de'=e^2)\to A_{PL}(S^k)$, when $k$ is even.
	
	Suppose $\alpha\in \pi_k(X)$ is represented by $a:(S^k,*)\to (X,*)$.
	Then $Q(a):V^k \to \Bbbk \cdot e$ depends only on $\alpha$ and the choice of $m_X$.
	Define the pairing \[\langle\cdot;\cdot\rangle: V\times \pi_*(X)\to \Bbbk\] by the equation \[\langle v;\alpha \rangle = \begin{cases} Q(a)v &, v\in V^k\\ 0 &, deg(v)\neq deg(\alpha)\end{cases}.\]
	It follows that $\langle\cdot;\cdot\rangle$ bilinear and, for $f:(Y,*)\to (X,*)$ a continuous map between simply connected spaces, \[\langle Q(f)v;\beta\rangle = \langle v; \pi_*(f)\beta\rangle,\ v\in V,\ \beta\in\pi_*(Y).\]
	
	This is important because it shows that $\langle \cdot; \cdot\rangle$ induces a natural linear map \[\nu: V\to Hom_\Z(\pi_*(X),\Bbbk),\ v\mapsto \langle v; \cdot \rangle\]

\begin{teorema}[\cite{FelixHalperinThomas_2001_RatHom}, 15.11; adapted for $\cat{PsTop}$]
\label{theo:FHT_15.11}
	Suppose $X$ is simply connected and $H_*(X,\Bbbk)$ has finite type.
	Then the bilinear map $V_X\times \pi_*(X)\to \Bbbk$ is non-degenerate.
	Equivalently, \[\nu_X: V_X \xrightarrow{\cong} Hom_\Z (\pi_*(X),\Bbbk)\] is an isomorphism.
\end{teorema}

\begin{proof}
	Let $X$ be a simply connected pseudotopological space such that $H^*(X,\Bbbk)$ has finite type.
	Then $X^\circ$, defined in \autoref{theo:PsTopTopWHE}, is a topological space such that it is simply connected and $H^*(X,\Bbbk)$ has finite type; therefore $X^\circ$ satisfies theorem 15.11 from \cite{FelixHalperinThomas_2001_RatHom}, i.e., \[\nu_{X^\circ}: V_{X^\circ} \xrightarrow{\cong} Hom_\Z (\pi_*(X^\circ),\Bbbk)\] is an isomorphism.
	
	It follows from \autoref{theo:RatHomTypMinSullAlg} that $X$ and $X^\circ$ have the same minimal Sullivan model up to isomorphism.
	Let's say that the minimal model is $m_{X}:(\Lambda V, d) \to A_{PL}(X)$.
	Additionally, $|Sing(X)_\cat{PsTop}|_\cat{PsTop}$ also have the same minimal Sullivan model.
	By the construction in \autoref{theo:PsTopTopWHE}, we have the weak homotopy equivalences \[\varepsilon_X: \varepsilon X \to X\text{ and } \eta_{\varepsilon X}: \varepsilon X \to X^\circ,\] with $\varepsilon X = |Sing(X)_\cat{PsTop}|_\cat{PsTop}$ and $\eta X = \eta_{\varepsilon X}$.
	
	Observe that $\Id{\Lambda V} m_X = m_X$, and then $Q(\varphi_{\Id{\Lambda V}})= \Id{V}$. 
	Now, let $\beta \in \pi_n(\varepsilon X)$ be represented by $b:S^n\to \varepsilon X$, then $[\varepsilon_X \circ b] \in \pi_n(X)$ and $[\eta_X\circ b]\in \pi_n(X^\circ)$.
	Thus \[ \langle v; \beta \rangle = \langle v ; \pi_*(\varepsilon_X)\beta \rangle = \langle v ; \pi_*(\eta_X)\beta \rangle,\ v\in V.\]
	Given that $\nu_{X^\circ}: V_{X^\circ} \xrightarrow{\cong} Hom_\Z (\pi_*(X^\circ), \Bbbk)$, then \[\nu_{X}: V_{X} \xrightarrow{\cong} Hom_\Z (\pi_*(X), \Bbbk)\] is an isomorphism.
\end{proof}

	It follows from the theorem above that, if $H_*(X,\Q)$ has finite type, $\pi_*(X)\otimes\Q$ has finite type, since in this case the case the Sullivan model also has a finite type.

\begin{definicion}[\cite{FelixHalperinThomas_2001_RatHom}, Page 247]
\label{def:SullivanRealization}
	The \emph{Sullivan realization}\index{Sullivan realization}
	\nomenclature{$\langle A,d\rangle$}{Sullivan realization of a differential graded algebra}
	is the contravariant functor $(A,d)\mapsto \langle A,d\rangle$ from commutative cochain algebra to simplicial sets, given by:
\begin{itemize}
	\item The $n$-simplices of $\langle A, d\rangle$ are the dga morphisms $\sigma: (A,d)\to (A_{PL})_n$.
	
	\item The face and degeneracy operators are given by $\partial_i\sigma = \partial_i\circ\sigma$ and $s_j\sigma = s_j\circ \sigma$.
	
	\item If $\varphi:(A,d)\to (B,d)$ is a morphism of commutative cochain algebras then $\langle\varphi\rangle:\langle A,d\rangle\leftarrow \langle B,d\rangle$ is the simplicial morphism given by $\langle\varphi\rangle(\sigma)=\sigma\circ\varphi$ with $\sigma\in \langle B,d\rangle_n$.
\end{itemize}
\end{definicion}
	We can define a natural bijection \[\mcat{(A,d)}{A_{PL}(K)}{DGA} \xrightarrow{\cong} \mcat{K}{\langle A,d \rangle}{sSet},\ \varphi\mapsto f\] is defined by \[f(\sigma)(a) = \varphi(a)(\sigma),\ a\in A,\sigma \in K_n, n\geq 0 .\]
	
	Establishing that $A_{PL}(\cdot)$ and $\langle \cdot \rangle$ as adjoint functors between commutative cochain algebras and simplicial sets.
	More specifically, the pair forms a Quillen adjoint \cite{hess2007rationalhomotopytheorybrief}.
	In particular, adjoint to the identity of $\langle A,d\rangle$ is the canonical dga morphism \[\eta_A:(A,d)\to A_{PL}\langle A,d\rangle.\]
	
\begin{definicion}[\cite{FelixHalperinThomas_2001_RatHom}, Page 248]
\label{def:SpatialRealization}
	The \emph{spatial realization of a commutative cochain algebra $(A,d)$}\index{Spatial realization}
	\nomenclature{$|A,d|_\cat{PsTop}$}{$|\langle A,d\rangle|_\cat{PsTop}$}
	is (the CW complex) $|A,d|_\cat{PsTop} = |\langle A,d\rangle|_\cat{PsTop}$. The \emph{spatial realization of a morphism $\varphi: (A,d)\to (B,d)$} is the continuous map $|\varphi| = |\langle\varphi\rangle|_\cat{PsTop}$.
\end{definicion}

	Defining $\zeta_n: \pi_n(|\Lambda V,d|_\cat{PsTop}) \to Hom_\Bbbk(V,\Bbbk)$ as $\zeta_n(\alpha)(v) = (-1)^n\langle v;\alpha \rangle$, with $\alpha\in \pi_n(X)$, $\zeta_n$ are morphisms of abelian groups if $n\geq 2$ and we have the following result

\begin{teorema}[	\cite{FelixHalperinThomas_2001_RatHom}, 17.10; $\cat{PsTop}$ adapted]
\label{theo:FHT_17.10}
	Let $(\Lambda V,d)$ be a Sullivan algebra such that $H^1(\Lambda V,d) = 0$ and each $H^p(\Lambda V,d)$ is finite dimensional. Then
\begin{enumerate}
	\item $|\Lambda V,d|_\cat{PsTop}$ is simply connected and $\zeta_n: \pi_n(|\Lambda V,d|_\cat{PsTop})\xrightarrow{\simeq} Hom_\Bbbk(V,\Bbbk)$ is an isomorphism, $n\geq 2$.
	
	\item If $\Bbbk = \Q$ then $m_{(\Lambda V,d),\cat{PsTop}}$ is a quasi-isomorphism \[m_{(\Lambda V,d),\cat{PsTop}}: (\Lambda V,d) \xrightarrow{\simeq} A_{PL}(|\Lambda V,d|_\cat{PsTop}).\]
\end{enumerate}
\end{teorema}

	When $\Bbbk=\Q$ then the theorem above exhibits any simply connected Sullivan algebra of finite type as the Sullivan model of a pseudotopological CW complex, and then a CW complex too.

\begin{proof}
	Let $(\Lambda V,d)$ be a Sullivan algebra such that $H^1(\Lambda V,d) = 0$ and each $H^p(\Lambda V,d)$ is finite dimensional.
\begin{enumerate}
	\item Recall $\eta_{\cat{PsTop},\cat{Top}}: \Id{\cat{PsTop}} \to \iota\tau(\cdot)$ from \autoref{theo:PsTopTopWHE} and denote this functor by $\eta_\cat{Ps}$ in this proof.
	Since $\tau$ is a left adjoin, thus it preserve colimits, and $|\Lambda V,d|_\cat{PsTop}$ is a cofibrant object, we have that $\iota\tau(|\Lambda V,d|_\cat{PsTop}) = |\Lambda V,d|$ is weak homotopy equivalent to $|\Lambda V,d|_\cat{PsTop}$.
	
	Since the object $|\Lambda V,d|$ is a topological space, applying theorem 17.10 from \cite{FelixHalperinThomas_2001_RatHom}, we have that $\zeta_n: \pi_n(|\Lambda V,d|)\xrightarrow{\simeq} Hom_\Bbbk(V,\Bbbk)$ and $|\Lambda V,d|$ is simply connected.
	In addition we know that both spaces have the same minimal Sullivan model up to isomorphisms.
	Thus, using an analogous argument to \autoref{theo:FHT_15.11}, we obtain the desired result because $\langle \cdot,\cdot\rangle_\cat{PsTop}$ is essentially the same map as $\langle \cdot,\cdot\rangle$.
	
	\item This follows directly from $|\Lambda V,d|$ being weak homotopy equivalent to $|\Lambda V,d|_\cat{PsTop}$ and apply the lifting \autoref{lem:LiftingLemma}.
	To prove it, we have the following commutative diagram
	\begin{align*}
	\xymatrix{
	A_{PL}(|\Lambda V,d|_\cat{PsTop}) \ar@{<->}[rr] \ar[rrdd]^<<<<<<{A_{PL}(\xi_\cat{PsTop})} & {} & A_{PL}(|\Lambda V,d|) \ar[dd]^{A_{PL}(\xi)} \\
	{} & {} & {}\\
	(\Lambda V,d) \ar@{.>}[uu]^{m_{(\Lambda V,d),\cat{PsTop}}} \ar[rr]_{\eta_{(\Lambda V,d)}} \ar@{.>}[rruu]_<<<<<<<{m_{(\Lambda V,d)}} & {} & A_{PL}\langle \Lambda V,d \rangle
	}
	\end{align*}
	In the diagram above, $A_{PL}(\xi_\cat{PsTop})$, $A_{PL}(\xi)$, $A_{PL}(\xi_\cat{PsTop})$ and $m_{(\Lambda V,d)}$ are quasi-isomorphisms.
	Then the morphism $\eta_{(\Lambda V,d)}$ is a quasi-isomorphism because it is composition of two.
	The morphism $m_{(\Lambda V,d),\cat{PsTop}}$ is a quasi-isomorphism because composed with a quasi-isomorphism is a quasi-isomorphism.
	Thus $m_{(\Lambda V,d),\cat{PsTop}}$ is a quasi-isomorphism as desired.
	\end{enumerate}
\end{proof}

	The quasi-isomorphism $m_{(\Lambda V,d),\cat{PsTop}}: (\Lambda V,d) \xrightarrow{\simeq} A_{PL}(|\Lambda V,d|_\cat{PsTop})$ is called a \emph{canonical Sullivan model}.\index{Sullivan model!canonical}

	Now let $X$ be any simply connected pseudotopological CW complex with rational homology of finite type, and let \[m_X: (\Lambda W,d) \to A_{PL}(X)\] be a minimal Sullivan model.
	Since $A_{PL}(X)$ is defined as $A_{PL}(S_*(X))$, the adjointness formula produces a natural simplicial map
	\[\gamma_X: S_*(X) \to \langle \Lambda W,d \rangle,\]
	is adjoint to $m_X$.
	On the other hand, $s_X: |S_*(X)| \to X$ is a natural weak homotopy equivalence, and then it is a natural homotopy equivalence because $X$ is a pseudotopological CW complex.
	Let $t_X$ be the inverse homotopy equivalence (defined uniquely up to homotopy) and set
	\[h_X \coloneqq |\gamma_X|\circ t_X : X \to |\Lambda W,d|.\]

\begin{teorema}[	\cite{FelixHalperinThomas_2001_RatHom}, 17.12; $\cat{PsTop}$ adapted]
\label{theo:FHT_17.12}
	With the notation and hypothesis above,
	\begin{enumerate}
	\item The diagram
	\begin{align*}
	\xymatrix{
	A_{PL}\langle \Lambda W, d \rangle & {} & A_{PL}| \Lambda W, d | \ar[ll]_{A_{PL}(\xi_{\Lambda W, d})} \ar[rr]^{A_{PL}(h_X)} & {} & A_{PL}(X)\\
	{} & {} & {} & {} & {}\\
	{} & {} & (\Lambda W, d) \ar[lluu]^{\eta_{(\Lambda W,d)}} \ar[uu]_{m_{(\Lambda W,d)}} \ar[uurr]_{m_X} & {} & {}
	}
	\end{align*}
	is homotopy commutative.
	
	\item If $\Bbbk = \Q$ then all the morphisms in the diagram are a quasi-isomorphism.
	In particular, $h_X$ is a rationalization of $X$.
	\end{enumerate}
\end{teorema}

\begin{proof}
	Since we have already proved \autoref{theo:FHT_9.6}, \autoref{prop:FHT_12.6} and \autoref{theo:FHT_17.10} in the pseudotopological space context, the proof is exactly the same as \cite{FelixHalperinThomas_2001_RatHom}.
\end{proof}

	With \autoref{theo:FHT_17.10} and \autoref{theo:FHT_17.12}, joint with \autoref{theo:RatHomTypMinSullAlg}, we have completed the following bijection \[ \left\lbrace \begin{array}{c}
\text{rational homotopy}\\ \text{types}	
\end{array}\right\rbrace \xrightarrow{\cong} \left\lbrace \begin{array}{c}
\text{isomorphism classes of}\\ \text{minimal Sullivan algebras}\\ \text{over }\Q	
\end{array} \right\rbrace \]
	Observe that if $\varphi: (\Lambda V,d) \to (\Lambda W, d)$ is an isomorphism, then the maps $(\Lambda V,d) \to A_{PL}(|\Lambda V,d|_\cat{PsTop})$ and $(\Lambda W,d) \to A_{PL}(|\Lambda V,d|_\cat{PsTop})$ are minimal Sullivan models, and then $A_{PL}(|\Lambda V,d|_\cat{PsTop})$ has the same rational homotopy of $A_{PL}(|\Lambda W,d|_\cat{PsTop})$.

In \autoref{prop:FHT_12.6}, we have proved that a pair of homotopic continuous maps $f_0 \simeq f_1: X \to Y$ have homotopic Sullivan representative. We have the following reverse direction:

\begin{teorema}[	\cite{FelixHalperinThomas_2001_RatHom}, 17.13; $\cat{PsTop}$ adapted]
\label{theo:FHT_17.13}
	If $\varphi_0\simeq \varphi: (\Lambda V,d)\to (\Lambda W,d)$ are homotopic morphisms of Sullivan algebras then \[|\varphi_0|_\cat{PsTop}\simeq |\varphi_1|_\cat{PsTop}:|\Lambda W,d|_\cat{PsTop} \to |\Lambda V,d|_\cat{PsTop}.\]
\end{teorema}

\begin{proof}
	Since $|\cdot|_\cat{PsTop}$ preserves the cartesian product, the proof is exactly the same as 17.13 in \cite{FelixHalperinThomas_2001_RatHom}.
\end{proof}

\begin{teorema}[	\cite{FelixHalperinThomas_2001_RatHom}, 17.15; $\cat{PsTop}$ adapted]
\label{theo:FHT_17.15}\
	\begin{enumerate}
	\item Let $\varphi: (\Lambda V,d) \to (\Lambda W,d)$ be a morphism of Sullivan algebras. Then
	\begin{align*}
	\xymatrix{
	(\Lambda V,d) \ar[rr]^\varphi \ar[d]_{m_{(\Lambda V,d)}} & & Y \ar[d]^{m_{\Lambda W,d)}}\\
	A_{PL}(|\Lambda W, d|_\cat{PsTop}) \ar[rr]_{A_{PL}(|\varphi|_\cat{PsTop})} & & A_{PL}(|\Lambda V, d|_\cat{PsTop})
	}
	\end{align*}	
	is homotopy commutative.
	In particular, if $\varphi$ is a morphism of simply connected minimal Sullivan algebras of finite type, then $\varphi$ is a Sullivan representative of $|\varphi|_\cat{PsTop}$.
	
	\item Let $f:X \to Y$ be a continuous map between simply connected pseudotopological CW complexes with finite type rational homology.
	Let $m_X: (\Lambda W,d) \to A_{PL}(X)$ and $m_Y: (\Lambda V,d) \to A_{PL}(Y)$ be minimal rational Sullivan models and let $\varphi: (\Lambda V,d) \to (\Lambda W, d)$ be a Sullivan representative for $f$.
	Then the diagram
	\begin{align*}
	\xymatrix{
	X \ar[rr]^f \ar[d]_{h_X} & & Y \ar[d]^{h_Y}\\
	|\Lambda W, d|_\cat{PsTop} \ar[rr]_{|\varphi|_\cat{PsTop}} & & |\Lambda V, d|_\cat{PsTop}
	}
	\end{align*}
	is homotopy commutative.
	\end{enumerate} 
\end{teorema}

\begin{proof}
	The proof is the same as \cite{FelixHalperinThomas_2001_RatHom} without any change up the notation of $\cat{PsTop}$.
\end{proof}	
	
	With \autoref{theo:FHT_17.13} and \autoref{theo:FHT_17.15}, we have completed the following bijection \[ \left\lbrace \begin{array}{c}
\text{homotopy classes of}\\ \text{continuous maps of}\\ \text{rational spaces}	
\end{array}\right\rbrace \xrightarrow{\cong} \left\lbrace \begin{array}{c}
\text{homotopy classes of}\\ \text{morphisms of Sullivan algebras}\\ \text{over }\Q	
\end{array} \right\rbrace \]
	
	This is given because:
	\begin{itemize}
	\item If we map $[f]$ to a Sullivan representative $\varphi$, this is well defined because all the Sullivan representatives are homotopic (\autoref{prop:FHT_12.9}) and if $f'\in [f]$ it has the same Sullivan representative given \autoref{prop:FHT_12.6} and $m_X \alpha \simeq A_{PL}(f)m_Y \simeq A_{PL}(f')m_Y$.
	\item If we map a homotopy class of morphism of Sullivan algebra $|\varphi|$ to $|\varphi|_\cat{PsTop}$, the map is well defined by \autoref{theo:FHT_17.13}.
	
	\item Lastly, we can observe that this an injection because $\varphi$ is a representative for both $f$ and $|\varphi|_\cat{PsTop}$.
	Thus $h_Y|\varphi|_\cat{PsTop} \simeq |\varphi|_\cat{PsTop} h_X \simeq h_Y f$.
	By definition of $h_Y$, and given it is a homotopy between CW complexes, $f\simeq |\varphi|_\cat{PsTop}$.
	\end{itemize}		

	With this, we conclude the relation between homotopies of spaces and homotopies of minimal Sullivan models.
	
	In a previous work \cite{Trevino_2024_GraphsVRPsTop}, we have seen the relation between the graphs, with its canonical \v{C}ech structure, and the geometric realization of their Vietoris-Rips complexes.
	We can observe then that, given \autoref{theo:RatHomTypMinSullAlg}, this two spaces have the same minimal Sullivan model, up to isomorphism, whenever the geometric realization is simply connected.
	In particular, every digital sphere $S^n$ has the same Sullivan model for $n\geq 2$ as the ordinary topological sphere $S^n$.

\section*{Discussion}

	This paper works as a brief first approach of the relation between pseudotopological spaces, simplicial sets and then rational homotopy.
	Atypically, we start the discussion section talking about the results in the appendixes.
	In the \autoref{Appx:LimPushoutEmb}, we introduce for this paper the pseudotopological spaces and the main results is
	
\begin{customcoro}{\ref{prop:EmbPushoutPsTop}}
	Let $X,Y,Z$ be pseudotopological spaces, $f:X\rightarrow Y$ and $g:X\rightarrow Z$ continuous maps and the following commutative diagram the pushout of $f$ and $g$ in the corresponding category:
\begin{align*}
\xymatrix{
Y \ar[r]^{g'} & P\\
X \ar[u]^{f} \ar[r]^{g} & Z \ar[u]_{f'}
}
\end{align*}
	If $g$ is an embedding, that is a injective map which is homeomorphism to its image, then $g': Y\rightarrow P$ is an embedding.
\end{customcoro}

This works as an opening for the subsequent \autoref{Appx:ProductCWComplex}. As the name indicates, the main result in this section is

\begin{customthm}{\ref{theo:CartesianProductCW}}
	Let $X$ and $Y$ be a pair of pseudotopological CW complexes, then the categorical product in $\cat{PsTop}$ is a pseudotopological CW complex.
\end{customthm}
	
	And the result follows from a chasing through different commutative diagrams with most of the arrows being embeddings.
	
	Having a bit context of the pseudotopological spaces and their CW complexes. In \autoref{cap.PsTopCW}, we rewrite the results using $\cat{PsTop}$ to prove that every retraction of a CW complex is homotopy equivalent to a CW complex:
	
\begin{customprop}{\ref{prop:Hatcher_A.11}}
	A space dominated by a CW complex is homotopy equivalent to a CW complex.
\end{customprop}

	The intention with that proposition was comparing basic properties of the two kinds of CW complexes and having a first outlook of the cofibrant objects in $\cat{PsTop}$.
	
	We know the geometric realization in from simplicial sets into topological spaces.
	In \autoref{cap.SSetCat}, we recall the simplicial category $\cat{\Delta}$ and the category of simplicial sets $\cat{sSet}$ to redesign the geometric realization but now into pseudotopological spaces.
	
	The changes are subtle, but we consider that include once again this construction can be useful if someone is reading it for first time and to warranty we do not need extreme changes.
	The main result of this section is

\begin{customthm}{\ref{theo:GeomRealizationCW}}
	Let $X$ be a simplicial set. The pseudotopological geometric realization of $X$, $|X|$ is a pseudotopological CW complex.
\end{customthm}	

%Section 3
	\autoref{sec.PsTopQuillenSSet} establish the adjoint functors between $\cat{PsTop}$ and $\cat{sSet}$ which results to be a Quillen equivalence.
	There we recall the definition of model categories the recover the Quillen model categories of $\cat{Top}$ and $\cat{sSet}$, also citing the Quillen model category of $\cat{PsTop}$ studied in \cite{Rieser_arXiv_2022}.
	
	We cite some results of cofibrantly generated model categories to observe that
	
\begin{customthm}{\ref{theo:PsTopSSetQAdj}}
	 The functors $|\cdot|_{PsTop}\dashv Sing_{PsTop}$ form a Quillen adjoint.
\end{customthm}	

\begin{customthm}{\ref{theo:PsTopsSetQuillenEq}}
	$(|\cdot|_{PsTop},Sing_{PsTop})$ is a Quillen equivalence.
\end{customthm}	

	Additionally, we use \autoref{coro:UnitIsAWeakEq} to observe that
	
\begin{customthm}{\ref{theo:PsTopWECW}}
	 Every pseudotopological space is weak homotopy equivalent to some pseudotopological CW complex.
\end{customthm}	

%Section 4	

	Now we can start studying a bit deeper the singular homology and its relations with homotopy in $\cat{PsTop}$.
	We proved that weak homotopy equivalences, and then strong homotopy equivalences, induces isomorphisms in the homology and cohomology groups in our category of interest.
	
\begin{customprop}{\ref{prop:Hatcher_4.21}}
	A weak homotopy equivalence $f:X\rightarrow Y$ induces isomorphisms $f:H_n(X;G)\rightarrow H_n(Y,G)$ and $f:H^*(Y;G)\rightarrow H^*(X;G)$ for all $n$ and all coefficient groups $G$.
\end{customprop}

	That important observation, together with a corollary in Quillen equivalences, allows us to show an important theorem in the translation of the rational homotopy theory from $\cat{Top}$ to $\cat{PsTop}$.

\begin{customthm}{\ref{theo:PsTopTopWHE}}
	Let $f:X\rightarrow Y$ in $\cat{PsTop}$, then there exists a topological morphisms  $f^\circ:X^\circ\rightarrow Y^\circ$ such that $f$ is a weak homotopy equivalence if and only if $f^\circ$ is a weak homotopy equivalence.
	That implies the same result replacing weak homotopy equivalence by quasi-isomorphisms.
\end{customthm}

	Firstly, we can obtain a version for pseudotopological space of the following Whitehead theorem.
	
\begin{customcoro}{\ref{theo:FHT_8.6}}
	Suppose that $\Bbbk$ is a subring of $\mathbb{Q}$ and $g:X\rightarrow Y$ is a morphism in $\cat{PsTop}$ between simply connected spaces.
	Then the following assertions are equivalent:
\begin{enumerate}
\item $\pi_*(g)\otimes\Bbbk$ is an isomorphism.
\item $H_*(g;\Bbbk)$ is an isomorphism.
%\item $H_*(\Omega g;\Bbbk)$ is an isomorphism.
\end{enumerate}
\end{customcoro}

	We define a morphism $f:X \to Y$ as a rational homotopy equivalence in the way as the category of topological spaces, i.e., $\pi_*(f)\otimes \Q$ is an isomorphism or, equivalently, $H_*(f;\Q )$ is an isomorphism. Showing later that
	
\begin{customprop}{\ref{theo:FHT_9.8}}
	Simply connected spaces $X$ and $Y$ have the same rational homotopy type if and only if there is a chain of rational homotopy equivalences
	\begin{align*}
	X \leftarrow Z(0) \rightarrow \cdots \leftarrow Z(k)\rightarrow Y.
	\end{align*}
\end{customprop}
	
%Section 5

	The exploration presented earlier allows us to obtain the definitions and results necessary for rational homotopy theory and to corroborate the same initial assertions that hold true for $\cat{Top}$.

	First, we obtain a formulation to send the a class of simply connected pseudotopological spaces with the same rational homotopy type

\begin{customthm}{\ref{theo:RatHomTypMinSullAlg}}
	Let $X$ and $Y$ be a pair of simply connected pseudotopological spaces with the same rational homotopy type. Then they have the same minimal Sullivan model up to isomorphism.
\end{customthm}

	Next we verify that there exist an isomorphism from a vector space to $Hom_\Z(\pi_*(X),\Bbbk)$.
	
\begin{customthm}{\ref{theo:FHT_15.11}}
	Suppose $X$ is simply connected and $H_*(X,\Bbbk)$ has finite type.
	Then the bilinear map $V_X\times \pi_*(X)\to \Bbbk$ is non-degenerate.
	Equivalently, \[\nu_X: V_X \xrightarrow{\cong} Hom_\Z (\pi_*(X),\Bbbk)\] is an isomorphism.
\end{customthm}
	
	This vector space is important because comes from the minimal Sullivan model and, when $H_*(X,\Q)$ is finite type, then every dimension of $V_X$ has finite dimension.
	Concluding that $\pi_*(X)\otimes \Q$ has finite type.
	
	With the same bilinear operator, we can define $\zeta_n$, and extend the following theorem from $\cat{Top}$ to $\cat{PsTop}$.
	
\begin{customthm}{\ref{theo:FHT_17.10}}
	Let $(\Lambda V,d)$ be a Sullivan algebra such that $H^1(\Lambda V,d) = 0$ and each $H^p(\Lambda V,d)$ is finite dimensional. Then
	\begin{enumerate}
	\item $|\Lambda V,d|_\cat{PsTop}$ is simply connected and $\zeta_n: \pi_n(|\Lambda V,d|_\cat{PsTop})\xrightarrow{\simeq} Hom_\Bbbk(V,\Bbbk)$ is an isomorphism, $n\geq 2$.
	\item If $\Bbbk = \Q$ then $m_{(\Lambda V,d),\cat{PsTop}}$ is a quasi-isomorphism \[m_{(\Lambda V,d),\cat{PsTop}}: (\Lambda V,d) \xrightarrow{\simeq} A_{PL}(|\Lambda V,d|_\cat{PsTop}).\]
	\end{enumerate}
\end{customthm}

	And concluding the following bijection in the category of pseudotopological spaces
	
\[ \left\lbrace \begin{array}{c}
\text{rational homotopy}\\ \text{types}	
\end{array}\right\rbrace \xrightarrow{\cong} \left\lbrace \begin{array}{c}
\text{isomorphism classes of}\\ \text{minimal Sullivan algebras}\\ \text{over }\Q	
\end{array} \right\rbrace \]

	Reducing our objects to pseudotopological simply connected CW complexes. We extend the following last result
	
\begin{customthm}{\ref{theo:FHT_17.15}}\
	\begin{enumerate}
	\item Let $\varphi: (\Lambda V,d) \to (\Lambda W,d)$ be a morphism of Sullivan algebras. Then
	\begin{align*}
	\xymatrix{
	(\Lambda V,d) \ar[rr]^\varphi \ar[d]_{m_{(\Lambda V,d)}} & & Y \ar[d]^{m_{\Lambda W,d)}}\\
	A_{PL}(|\Lambda W, d|_\cat{PsTop}) \ar[rr]_{A_{PL}(|\varphi|_\cat{PsTop})} & & A_{PL}(|\Lambda V, d|_\cat{PsTop})
	}
	\end{align*}	
	is homotopy commutative.
	In particular, if $\varphi$ is a morphism of simply connected minimal Sullivan algebras of finite type, then $\varphi$ is a Sullivan representative of $|\varphi|_\cat{PsTop}$.
	
	\item Let $f:X \to Y$ be a continuous map between simply connected pseudotopological CW complexes with rational homology of finite type.
	Let $m_X: (\Lambda W,d) \to A_{PL}(X)$ and $m_Y: (\Lambda V,d) \to A_{PL}(Y)$ be minimal rational Sullivan models and let $\varphi: (\Lambda V,d) \to (\Lambda W, d)$ be a Sullivan representative for $f$.
	Then the diagram
	\begin{align*}
	\xymatrix{
	X \ar[rr]^f \ar[d]_{h_X} & & Y \ar[d]^{h_Y}\\
	|\Lambda W, d|_\cat{PsTop} \ar[rr]_{|\varphi|_\cat{PsTop}} & & |\Lambda V, d|_\cat{PsTop}
	}
	\end{align*}
	is homotopy commutative.
	\end{enumerate}
\end{customthm}
	
	With that theorem we finally find a bijection in $\cat{PsTop}$ with the homotopy classes

\[ \left\lbrace \begin{array}{c}
\text{homotopy classes of}\\ \text{continuous maps of}\\ \text{rational spaces}	
\end{array}\right\rbrace \xrightarrow{\cong} \left\lbrace \begin{array}{c}
\text{homotopy classes of}\\ \text{morphisms of Sullivan algebras}\\ \text{over }\Q	
\end{array} \right\rbrace \]

\appendix
\section{Limit Spaces and The Pushout of an Embedding}
\label{Appx:LimPushoutEmb}

By completeness, we begin this section given an introduction of limit spaces and pseudotopological spaces. First, we introduce the notion of topological construct:

\begin{definicion}[\cite{Preuss_2002}; 1.1, 1.2]
\label{def:TopConstructs}
By a construct we mean a concrete category $\cat{C}$ whose objects are structures sets, i.e., pairs $(X,\xi)$ where $X$ is a set and $\xi$ a ``$\cat{C}$-structure on $X$'', whose morphisms $f:(X,\xi)\rightarrow (Y,\eta)$ are suitable maps between $X$ and $Y$ and whose composition laws are the usual composition of ($\cat{Sets}$) maps.

A construct $\cat{C}$ is called topological if and only if satisfies the following conditions:
\begin{enumerate}
\item \label{def:InitialStructure} Existence of initial structures:

{\noindent For any set $X$, any family $((X_i,\xi_i))_{i\in I}$ of $\cat{C}$-objects indexed by a class $I$ and any family $(f_i: X\rightarrow X_i)_{i\in I}$ of maps indexed by $I$ there exists a unique $\cat{C}$-structure $\xi$ on $X$ which is \emph{initial} with respect to $(X, f_i, (X_i,\xi_i), I)$, i.e., such that for any $\cat{C}$-object $(Y,\eta)$ a map $g:(Y,\eta)\rightarrow (X,\xi)$ is a $\cat{C}$-morphism if and only if for every $i\in I$ the composite map $f_i\circ g:(Y,\eta)\rightarrow (X_i,\xi_i)$ is a $\cat{C}$-morphism.}

\item For any set $X$, the class $\{ (Y,\eta)\in \ocat{C} \mid X=Y \}$ of all $\cat{C}$-objects with underlying set $X$ is a set.

\item For any set $X$ with cardinality at most one, there exists exactly one $\cat{C}$-object with underlying set $X$.
\end{enumerate}
\end{definicion}

The category of topological spaces are the obvious example of a topological construct; however, there are some categories that contain topological spaces and have richer categorical properties and contain examples which are interesting in discrete homotopy, for example, the graphs. These categories are the limit spaces and the pseudotopological spaces.

\begin{definicion}[\cite{Preuss_2002}, 2.3.1.1]
For each set $X$, let $F(X)$ be the set of all filters on $X$. A limit space is a pair $(X,q)$ where $X$ is a set and $q\subset X\times F(X)$ such that the following are satisfied.
\begin{enumerate}
\item $([x],x)\in q$ for each $x\in X$, where $[x]=\{A\subset X: x\in A\}$,

\item $(\mathcal{G},x)\in q$ whenever $(\mathcal{F},x)\in q$ and $\mathcal{G}\supset \mathcal{F}$,

\item $(\mathcal{F}\cap \mathcal{G})\in q$ whenever $(\mathcal{F},x)\in q$ and $(\mathcal{G},x)\in q$.
\end{enumerate}
A map $f:(X,q)\rightarrow (X',q')$ between limit spaces is called continuous provided that $(f(\mathcal{F},f(x))\in q'$ for each $(\mathcal{F},x)\in q$.
\end{definicion}

\begin{definicion}[\cite{Preuss_2002}, 2.3.1.1]
\label{def:PsTopSpaces}
A pseudotopological spaces is a limit space which satisfies that $(\mathcal{F},x)\in q$ whenever $(\mathcal{U},x)\in q$ for every ultrafilter $\mathcal{U}\supset \mathcal{F}$.
\end{definicion}

The categories of limit spaces and pseudotopological spaces with the continuous maps are denoted as $\cat{Lim}$ and $\cat{PsTop}$, respectively.

	By Theorem 1.2.1.10 in \cite{Preuss_2002}, every topological construct has all of the small limits and small colimits.
	This properties make them great candidates to look for a model category for them and ensure the existence of familiar objects like: pullbacks, pushouts, direct limits and inverse limits.
	In this section, we extensively use the concept of pushout to build CW complexes and, later in this section, we observe a property for the embeddings in $\cat{Lim}$ and $\cat{PsTop}$. Before that, we remind this important definition from category theory which is used a lot during the section.

\begin{definicion}[\cite{Riehl_2016}, 4.1.1]
\label{def:adjoint}
An \emph{adjunction} consists of a pair of functors $\mfunt{F}{C}{D}$ and $\mfunt{G}{D}{C}$ together with an isomorphism
\begin{align*}
\mcat{\funt{F}c}{d}{D} \cong \mcat{c}{\funt{G}d}{C}
\end{align*}
for each $c\in \cat{C}$ and $d\in \cat{D}$ that is natural in both variables. Here $F$ is \emph{left adjoint} to $\cat{G}$ and $\cat{G}$ is \emph{right adjoint} to $\cat{F}$.
\end{definicion}

	It is a well-known result that in topology the pushout map of an open or closed embedding is, respectively an open or closed embedding (see, for example, Proposition 1.2.4 in \cite{tom_Dieck_2008}).
	This result is used in topology to note that the pushout of an embedding is actually an embedding.
	Furthermore, in \cite{bubenik2024_RelCellCompClos} Bubenik has proved that this result is preserved in $\cat{PrTop}$.
	
	The result in $\cat{Top}$ for closed maps is specially well-received in the construction of CW complexes to observe that $X_{i}\rightarrow X_{i+1}$ is an embedding.
	In this brief appendix section, we are going to observe that limit spaces also preserves embedding through pushout, without using the properties of open or closed map, and give the arguments to obtain that this result follows in pseudotopological spaces.
	It is important to mention that, if we only attach a finite amount of cells in every step $n$, this results for CW complexes follows directly for theorem 1.4.10 from \cite{Beattie_Butzmann_2002}.

\begin{proposicion}
\label{prop:EmbPushoutLim}
Let $X,Y,Z$ be limit spaces, $f:X\rightarrow Y$ and $g:X\rightarrow Z$ continuous maps and the following commutative diagram the pushout of $f$ and $g$:
\begin{align*}
\xymatrix{
Y \ar@{.>}[r]^{g'} & P \pushoutur\\
X \ar[u]^{f} \ar[r]^{g} & Z \ar@{.>}[u]_{f'}
}
\end{align*}
If $g$ is an embedding, that is a injective map which is homeomorphism to its image, then $g': Y\rightarrow P$ is an embedding.
\end{proposicion}

\begin{proof}
Since the limit spaces are a topological construct, then $P$ is isomorphic as a limit space to $Y\sqcup_X Z$, that is the limit space $Y\sqcup Z$ with the equivalence relation $g(x)\sim f(x)$ for every $x\in X$.

We want to show that $Y\simeq g'(Y)$. We use without proof that the pushout in $\cat{Sets}$ of an injective map is an injective map, thus $g'$ is an injective map as a $\cat{Sets}$ morphisms and bijective over its image.

	For the pushout, we already know that $g'$ is continuous, in particular, it is continuous onto its image;
	therefore, it is enough to prove that for every $p\in g'(Y)$ and every filter $\mathcal{F}\rightarrow p$ in $g'(Y)$ there exits a filter $\mathcal{J}\rightarrow (g')^{-1}(p)$ such that $g'(\mathcal{J})\subset \mathcal{F}$.

	Let $\mathcal{F}\rightarrow p$ be a convergent filter in $g'(Y)$.
	By the subspace structure, we have that $[\mathcal{F}]_P \rightarrow p$ in $Y\sqcup_X Z$.
	Thus, there exists $y\in Y$, $z_{k+1},z_{k+2},\ldots,z_n\in Z$, filters $\mathcal{G}_1,\ldots,\mathcal{G}_k$ in $Y$ and filters $\mathcal{H}_{k+1},\ldots,\mathcal{H}_n$ in $Z$ such that $g'(y)=p$, $f'(z_{k+1})=\ldots =f'(z_n)=p$,
\begin{align*}
	\mathcal{G}_i\rightarrow y & \text{ for every }i\in \{1,\ldots,k\}\\
	\mathcal{H}_i\rightarrow z_i & \text{ for every }i\in \{k+1,\ldots,n\}
\end{align*}
	and $g'(\mathcal{G}_1)\cap\ldots g'(\mathcal{G}_k)\cap f'(\mathcal{H}_{k+1})\cap\ldots\cap f'(\mathcal{H}_n)\subset [\mathcal{F}]_P$.
	Since $g'(Y)\in [\mathcal{F}]_P$, and by the construction of $P$, we can consider that $g(X)\in \mathcal{H}_i$.

	For the subspace structure of $g(X)$ in $Z$, and considering that $f'(\mathcal{H}_i)$ must intersect $g'(Y)$, we have that $\mathcal{H}_i\cap g(X)\rightarrow z_i$ in $g(X)$ and $[\mathcal{H}_i\cap g(X)]_Z= \mathcal{H}_i$.
	Given that $g$ is an embedding, then $g^{-1}(\mathcal{H}_i\cap g(X))\rightarrow g^{-1}(z_i)$ in $X$, and thus $fg^{-1}(\mathcal{H}_i\cap g(X))\rightarrow fg^{-1}(z_i)$ in $Y$.

	Since $g'fg^{-1}(z_i) = f'gg^{-1}(z_i)=p$, then $fg^{-1}(z_i)=y$ because $g'$ is injective.
	It is only necessary to prove that
\begin{align*}
	g'(\mathcal{G}_1)\cap\ldots g'(\mathcal{G}_k)\cap g'f(g^{-1}(\mathcal{H}_{k+1}\cap g(X)))\cap\ldots\cap g'f(g^{-1}(\mathcal{H}_n\cap g(X)))\subset [\mathcal{F}]_P.
\end{align*} 
	However, this follows because $g'f(X)=f'g(X)$ and $g'(Y)\in [\mathcal{F}]_P$, i.e., $[\mathcal{F}]_P$ is omitting all of the part of $Z$ which are not in $g'(Y)$.
%\todo{Fill the proof.}
\end{proof}

	In the page 58 of their book \cite{Beattie_Butzmann_2002}, Beattie and Butzmann claim that the Choquet (or in our case, pseudotopological) modification $\chi$ has a right adjoint, which is the ultrafilter modification, and thus the pseudotopological modification is a left adjoint.
	This result implies that the pseudotopological modification preserves colimits in the following sense:
\begin{align*}
	\chi\left( \lim\limits_{\rightarrow}  (X_i) \right) \simeq \lim\limits_{\rightarrow} (\chi(X_i)) 
\end{align*}

	In the particular case of the pushout, which is a colimit, we have the following commutative diagram for every $X,Y,Z$ pseudotopological spaces is the pushout:
\begin{align*}
	\xymatrix{
	Y \ar[r]^{\chi(g')} & \chi(P) \pushoutur\\
	X \ar[u]^{f} \ar[r]^{g} & Z \ar[u]_{\chi(f')}
	}
\end{align*}
	where $P$ is the pushout in the category of limit spaces.

	Since both $\cat{Lim}$ and $\cat{PsTop}$ are concrete categories and since every functor preserves isomorphisms, we obtain that every one-to-one map in $\cat{Lim}$ whose domain is homeomorphic to its image is mapped through $\chi$ to an one-to-one map in $\cat{PsTop}$ whose domain is homeomorphic to its image, too.
	Thus, we obtain the following corollary.

\begin{corolario}
\label{prop:EmbPushoutPsTop}
	Let $X,Y,Z$ be pseudotopological spaces, $f:X\rightarrow Y$ and $g:X\rightarrow Z$ continuous maps and the following commutative diagram the pushout of $f$ and $g$ in the corresponding category:
\begin{align*}
\xymatrix{
Y \ar[r]^{g'} & P \pushoutur\\
X \ar[u]^{f} \ar[r]^{g} & Z \ar[u]_{f'}
}
\end{align*}
	If $g$ is an embedding, that is a injective map which is homeomorphism to its image, then $g': Y\rightarrow P$ is an embedding.
\end{corolario}

\section{The Cartesian Product of two PsTop CW Complexes}
\label{Appx:ProductCWComplex}

	In his book \cite{Hatcher_2002} (page 524, after the proof of Theorem A.6), Hatcher illustrates that the $\cat{Top}$ product of two CW complex is not necessarily a CW complex.
	He provides sufficient conditions on a CW complex to guarantee that the $\cat{Top}$ product is actually a CW complex, illustrates an example which is not a CW complex, and provides a topological  structure in the Haursdoff compactly generated category (page 524, \cite{Hatcher_2002}) to see the $\cat{CGHaus}$ product of any pair of CW complexes as a CW complex in that category.

The category of pseudotopological spaces is very special for many reasons, one of which has been repeated in this document many times: the category is cartesian closed. A consequence of this is that the functors $A\times-$ and $-\times A$ have right adjoints for every $A\in \ocat{PsTop}$, and thus they preserves colimits.

In this appendix, we will see that this condition implies that the cartesian product of a pair of pseudotopological CW complexes is a pseudotopological CW complex. Specifically, we show the following.

\begin{teorema}
\label{theo:CartesianProductCW}
	Let $X$ and $Y$ be a pair of pseudotopological CW complexes, then the categorical product in $\cat{PsTop}$ is a pseudotopological CW complex.
\end{teorema}

Let $X$ and $Y$ be a pair of CW complexes with the following skeletons:
\begin{align*}
\xymatrix{
X_0 \ar[r]^{\rho_{0,1}} \ar[drr]_{\rho_0} & X_1 \ar[r]^{\rho_{1,2}} \ar[r] \ar[dr]^{\rho_1} & \cdots \ar[r]^{\rho_{n-1,n}} & X_n \ar[r]^{\rho_{n,n+1}} \ar[dl]^{\rho_n} & \cdots\\
& & X=\lim\limits_{\xrightarrow[i]{}} X_i \\
Y_0 \ar[r]^{\varrho_{0,1}} \ar[drr]_{\varrho_0} & Y_1 \ar[r]^{\varrho_{1,2}} \ar[r] \ar[dr]^{\varrho_1} & \cdots \ar[r]^{\varrho_{n-1,n}} & Y_n \ar[r]^{\varrho_{n,n+1}} \ar[dl]^{\varrho_n} & \cdots\\
& & Y=\lim\limits_{\xrightarrow[i]{}} Y_i
}
\end{align*}
	for simplicity in the notation, assume that if $X$ or $Y$ are finite-dimensional, then $X_k=X_{dim(X)}$ or $Y_\ell=Y_{dim(Y)}$ for $k\geq dim(X)$ and $\ell\geq dim(Y)$, and define $\rho_{i,i}\coloneqq \Id{X_i}$, $\varrho_{i,i}\coloneqq \Id{X_i}$, \[\rho_{i,j}\coloneqq\rho_{i,i+1}\circ \ldots \circ \rho_{j-1,j} \text{ and }
	\varrho_{i,j}\coloneqq\varrho_{i,i+1}\circ \cdots \circ \varrho_{j-1,j}\] for every $i\leq j$

	Since $\cat{PsTop}$ is cartesian closed, we can write $X\times Y$ in the following two ways.
\begin{align*}
	X\times Y = & \lim\limits_{\xrightarrow[i]{}} X_i\times Y = \lim\limits_{\xrightarrow[i]{}} \lim\limits_{\xrightarrow[j]{}} X_i\times Y_j\\
	 &  \lim\limits_{\xrightarrow[j]{}} X\times Y_j = \lim\limits_{\xrightarrow[j]{}} \lim\limits_{\xrightarrow[i]{}} X_i\times Y_j
\end{align*}
	In the following, we construct a limit which is isomorphic to $X\times Y$.
	First we define $A_i\coloneqq X_i\times Y_i$ and $\alpha_{i,i+1}\coloneqq\rho_{i,i+1}\times\varrho_{i,i+1}:A_i\rightarrow A_{i+1}$.
	Denote by $A\coloneqq \lim\limits_{\xrightarrow[i]{}} A_i$ and $\alpha_i:A_i\rightarrow A$, and observe that the diagrams
\begin{align*}
\xymatrix{
X_j\times Y_i \ar[rr]^{\Id{X_i}\times \varrho_{i,i+1}} \ar[dr]_{\rho_{j,i+j}\times\varrho_{i,i+j}} & & X_{j}\times Y_{i+1} \ar[d]^{\rho_{j,i+j+1}\times\varrho_{i+1,i+j+1}}\\ 
& A_{i+j} \ar[rd]_{\alpha_{i+j}} & A_{i+j+1} \ar[d]^{\alpha_{i+j+1}}\\
& & A
}
\end{align*}
are commutative for every $i,j\in \N$. Thus, since $Y$ is a colimit, there exists a unique $\Lambda^X_j:X_j\times Y\rightarrow A$ such that the following diagram commutes for every
\begin{align*}
\xymatrix{
X_j\times Y_i \ar[rrr] \ar[dd] \ar[dr] & & & A_{i+j} \ar[dd] \ar[dl]\\
& X_j\times Y \ar@{.>}[r]^{\Lambda^X_j} & A\\
X_j \times Y_{i+1} \ar[rrr] \ar[ru] & & & A_{i+j+1} \ar[lu]
}
\end{align*}
	Now we claim that the diagram
\begin{align*}
\xymatrix{
X_j\times Y \ar[rr]^{\rho_{j,j+1}\times \Id{Y}} \ar[rrd]_{\Lambda_{j}^X} & & X_{j+1}\times Y \ar[d]^{\Lambda_{j+1}^X}\\
& & A
}
\end{align*}
	commutes.
	To prove that claim, we have the following diagram that is commutative in the {\color{NavyBlue}inner diagrams}, and then the outside diagram is also commutative:
\begin{align*}
\xymatrix{
X_j\times Y_i \ar[ddd]_{\rho_{j,j+i}\times \varrho_{i,j+i}} \ar[rrr]^{\Id{X_j}\times\varrho_i} \ar@[NavyBlue][rd]^{\color{NavyBlue}\rho_{j,j+1}\times \Id{Y_i}} & & & X_j\times Y \ar[d]^{\rho_{j,j+1}\times \Id{Y}}\\
& {\color{NavyBlue}X_{j+1}\times Y_i} \ar@[NavyBlue][rr]^{\color{NavyBlue}\Id{X_{j+1}}\times\varrho_i} \ar@[NavyBlue][d]^{\color{NavyBlue}\rho_{j+1,j+i+1}\times\varrho_{i,j+i+1}} & & X_{j+1}\times Y \ar[dd]^{\Lambda^X_{j+1}}\\
& {\color{NavyBlue}A_{j+i+1}} \ar@[NavyBlue][drr]^{\color{NavyBlue}\alpha_{j+i+1}}\\
A_{j+i} \ar@[NavyBlue][ru]_{\color{NavyBlue}\alpha_{j+i,j+i+1}} \ar[rrr]^{\alpha_{j+1}} & & & A
}
\end{align*}
	for every $i$ and $j$.
	For the uniqueness of $\Lambda^X_j$, we obtain that $\Lambda^X_j = \Lambda_{j+1}^X (\rho_{j,j+1}\times \Id{Y})$. Thus, there exists a unique $\Lambda:X\times Y\rightarrow A$ such that $\Lambda_j^X = \Lambda (\rho_{j}\times \Id{Y})$.

Because $A$ is a colimit, there exists also  a unique $\Lambda': A\rightarrow X\times Y$ such that this given the commutative diagram
\begin{align*}
\xymatrix{
A_i \ar[r] \ar[dd] \ar[rd] & X_{i}\times Y \ar[rd]\\
& A \ar@{.>}^{\Lambda'}[r] & X\times Y\\
A_{i+1} \ar[ru]\ar[r] & X_{i+1}\times Y \ar[ru]
}
\end{align*}

	Since $X\times Y$ is the colimit of $X_k\times Y$, then there exists a unique $\Psi$ such that
\begin{align*}
\xymatrix{
X_j\times Y \ar[r]\ar[dd]\ar[rd] & A \ar[rd]^{\Lambda'} \\
& X\times Y \ar@{.>}[r]^{\Psi} & X\times Y\\
X_{j+1}\times Y \ar[ru]\ar[r] & A \ar[ru]_{\Lambda'}
}
\end{align*}

	Noting that the diagram below is commutative for every $i$ for the construction of $\Lambda_i^X$ and $\Psi$.
\begin{align*}
\xymatrix{
A_i \ar[rr]\ar[d] & & X_i\times Y \ar[d] \ar[dl] \ar[dll]_{\Lambda_i^X}\\
A \ar[r]_{\Lambda'} \ar@/_1.5pc/[rr]_{\Lambda'} & X\times Y \ar[r]^{\Psi} & X\times Y
}
\end{align*}
	We conclude that $\Psi\Lambda'=\Lambda'$.
	We have that the following diagrams
\begin{align*}
\xymatrix{
X\times Y \ar[r]  & A \ar[r]  & X\times Y \\
 X_i\times Y \ar[u] \ar[ru]\ar@/_1pc/[rru]
}\\
\xymatrix{
A_i \ar[rr]\ar[d] & & A \ar[d]^{\Lambda'} \\
X_i\times Y_j \ar[d]\ar[r] & X_i\times Y \ar[r] & X\times Y \ar[d]^{\Lambda}\\
A_{i+j} \ar[rr] & & A
}
\end{align*}
are commutative, then $\Lambda\circ \Lambda' = \Id{A}$ and $\Lambda'\circ \Lambda = \Id{X\times Y}$; concluding that $A\cong X\times Y$ in the category of pseudotopological spaces.

Now we define the pseudotopological spaces $B_i$ as follows:
\begin{align*}
B_0 \coloneqq & X_0\times Y_0\\
B_1 \coloneqq & (X_1\times Y_0) \bigsqcup\limits_{X_0\times Y_0} (X_0\times Y_1)\\
B_2 \coloneqq & (X_2\times Y_0) \bigsqcup\limits_{X_1\times Y_0} (X_1\times Y_1) \bigsqcup\limits_{X_0\times Y_1} (X_0\times Y_2)\\
B_i \coloneqq & (X_i\times Y_0) \bigsqcup\limits_{X_{i-1}\times Y_0} (X_{i-1}\times Y_1) \bigsqcup\limits_{X_{i-2}\times Y_1} \ldots \bigsqcup\limits_{X_0\times Y_{i-1}} (X_0\times Y_i)\
\end{align*}
It follows by definition that we have the inclusions $A_i\rightarrow B_i \rightarrow A_{2i}$, thus we have that $B\coloneqq \lim_i B_i \cong A$

The last step in to observe that $B_{i+1}$ is equal to $B_i$ attaching some $(i+1)$-cells, obtaining that $B$ is a CW complex.  Since $X$ and $Y$ are CW complexes, for every $i\in\mathbb{N}_0$ we have that
\begin{align*}
\xymatrix{
\bigsqcup\limits_{\alpha\in I^X_i} S^i_\alpha \ar[r]\ar[d] & \bigsqcup\limits_{\alpha\in I^X_i} D^{i+1}_\alpha \ar[d] \\
X_i \ar[r] & X_{i+1} \pushoutdr
}\hspace{5mm}
\xymatrix{
\bigsqcup\limits_{\beta\in I^Y_i} S^i_\beta \ar[r]\ar[d] & \bigsqcup\limits_{\beta\in I^Y_i} D^{i+1}_\beta \ar[d] \\
Y_i \ar[r] & Y_{i+1} \pushoutdr
}
\end{align*}

Since $\cat{PsTop}$ is cartesian closed we have the following commutative diagram (\autoref{Fig:HugeCommDiag})
\begin{figure}[h!]
{\small \begin{align*}
\xymatrixrowsep{6mm}
\xymatrixcolsep{5mm}
\xymatrix{
& *+[o][F]{1} \ar[d] & *+[o][F]{2} \ar[d] & *+[o][F]{3} \ar[d] & *+[o][F]{4} \ar[d]\\
 *+[o][F]{1'} & \sqcup S^i_\alpha \times Y_{j} \ar[r] \ar[l] \ar[d] & \sqcup S^i_\alpha \times Y_{j+1} \ar[r]\ar[d] \pushoutdr & X_i\times Y_{j+1} \ar[d] \pushoutdl & X_i\times Y_j \ar[l]\ar[d] & *+[o][F]{1'} \ar[l]  \\
*+[o][F]{2'} & \sqcup D_{\alpha}^{i+1}\times Y_j \ar[l]\ar[r] & \sqcup S_\alpha^i\times Y_j \ar[r] \pushoutur & X_{i+1}\times Y_{j+1}\pushoutul \pushoutdr  & X_{i+1}\times Y_j \ar[l] \pushoutdl & *+[o][F]{2'} \ar[l]\\
*+[o][F]{3'} & \sqcup D_\alpha^{i+1}\times \sqcup S_\beta^j \ar[l]\ar[r]\ar[u] & \sqcup D_\alpha^{i+1}\times \sqcup D_\beta^{j+1} \ar[r]\ar[u] & X_{i+1}\times \sqcup D_\beta^{j+1} \ar[u] \pushoutur & X_{i+1}\times \sqcup S_\alpha^j \ar[l]\ar[u] \pushoutul & *+[o][F]{3'} \ar[l] \\
*+[o][F]{4'} & \sqcup S_\alpha^i \times \sqcup S_\beta^j \ar[r]\ar[u]\ar[l]\ar[d] & \sqcup S_\alpha^i \times \sqcup D_\beta^{j+1} \ar[r]\ar[u]\ar[d] & X_i \times \sqcup D_\beta^{j+1} \ar[u] \ar[d] & \sqcup S_\beta^{j+1} \ar[u] \ar[d] \ar[l] & *+[o][F]{4'} \ar[l] \\
& *+[o][F]{1} & *+[o][F]{2} & *+[o][F]{3} & *+[o][F]{4} &
}
\end{align*}}
\caption{This relations are giving by construction and the pushouts.}\label{Fig:HugeCommDiag}
\end{figure}

We have the following pushout, where $S^{i+j+1}$ is the $i+j+1$-sphere,
\begin{align*}
\xymatrix{
D^{i+1}\times S^j \ar[r] & S^{i+j+1} \pushoutur\\
S^i\times S^j \ar[r]\ar[u] & S^i\times D^{j+1} \ar[u]
}
\end{align*}
Then, the we obtain
\begin{align*}
\xymatrix{
\sqcup D^{i+1}_\alpha \times \sqcup S^j_\beta \ar[r] & \pushoutur \sqcup S^{i+j+1}_{(\alpha,\beta)} \\
\sqcup S^{i}_\alpha \times \sqcup S^j_\beta \ar[r] \ar[u] & \sqcup S^{i}_\alpha \times \sqcup D^{j+1}_\beta \ar[u]
}
\end{align*}
	Let's denote the space $(X_{i}\times Y_{j+1})\sqcup_{X_{i}\times Y_{j}}(X_{i+1}\times Y_{j})$ by $X_iY_j^\times$ in the rest of the section looking for simplify a bit the writing.
	 Because $\sqcup S^{i+j+1}_{(\alpha,\beta)}$ and $X_iY_j^\times$ are pushout, there exists unique morphisms $\Psi:\sqcup S^{i+j+1}_{(\alpha,\beta)}\rightarrow \sqcup D^{i+1}_\alpha \times \sqcup D^{j+1}_\beta$ and $\Psi':X_iY_j^\times\rightarrow X_{i+1}\times Y_{j+1}$ such that
\begin{align*}
\xymatrix{
{} & {} & \sqcup D^{i+1}_\alpha\times \sqcup D^{j+1}_\beta\\
\sqcup D_\alpha^{i+1} \times \sqcup S_\beta^{j} \ar[r] \ar@/^1pc/[rru] & \pushoutur \sqcup S^{j+i+1}_{(\alpha,\beta)} \ar@{.>}[ru]^{\Psi} \\
\sqcup S_\alpha^{i} \times \sqcup S_\beta^{j} \ar[u]\ar[r] & \sqcup S_\alpha^{i} \times \sqcup D_\beta^{j+1} \ar[u] \ar@/_1pc/[ruu]
}
\end{align*}
\begin{align*}
\xymatrix{
{} & {} & X_{i+1}\times Y_{j+1}\\
X_{i}\times Y_{j+1} \ar[r] \ar@/^2pc/[rru] & \pushoutur X_iY_j^\times \ar@{.>}[ru]^{\Psi'} \\
X_{i}\times Y_{j} \ar[u]\ar[r] & X_{i+1}\times Y_{j} \ar[u] \ar@/_2pc/[ruu]
}
\end{align*}
Note that $\Psi$ and $\Psi'$ must be monomorphisms because the legs are monomorphisms and we have that the disjoint union of the pushout morphisms are onto the pushout. Then, given the following compositions of functions which are equal for \autoref{Fig:HugeCommDiag}
\begin{align*}
\sqcup S_\alpha^{i} \times \sqcup S_\beta^{j} \rightarrow & \sqcup D_\alpha^{i+1} \times \sqcup S_\beta^{j} \rightarrow \sqcup D_\alpha^{i+1} \times Y_j \rightarrow X_{i+1} \times Y_j \rightarrow X_iY_j^\times \\
\sqcup S_\alpha^{i} \times \sqcup S_\beta^{j} \rightarrow & \sqcup S_\alpha^{i} \times Y_j \rightarrow \sqcup D_\alpha^{i+1} \times Y_j \rightarrow X_{i+1}\times Y_j \rightarrow X_iY_j^\times \\
\sqcup S_\alpha^{i} \times \sqcup S_\beta^{j} \rightarrow & \sqcup S_\alpha^{i} \times Y_j \rightarrow X_i \times Y_j \rightarrow X_{i+1} \times Y_j \rightarrow X_iY_j^\times \\
\sqcup S_\alpha^{i} \times \sqcup S_\beta^{j} \rightarrow & \sqcup S_\alpha^{i} \times Y_j \rightarrow X_i \times Y_j \rightarrow X_{i} \times Y_{j+1} \rightarrow X_iY_j^\times \\
\sqcup S_\alpha^{i} \times \sqcup S_\beta^{j} \rightarrow & \sqcup S_\alpha^{i} \times Y_j \rightarrow \sqcup S_\alpha^{i} \times Y_{j+1} \rightarrow X_{i} \times Y_{j+1} \rightarrow X_iY_j^\times\\
\sqcup S_\alpha^{i} \times \sqcup S_\beta^{j} \rightarrow & \sqcup S_\alpha^{i} \times \sqcup D_\beta^{j+1} \rightarrow X_i \times \sqcup D_\beta^{j+1} \rightarrow X_{i} \times Y_{j+1} \rightarrow X_iY_j^\times
\end{align*}
there exists a unique $\Psi'':\sqcup S^{j+i+1}_{(\alpha,\beta)}\rightarrow X_iY_j^\times$ which makes the following diagram commutes

\begin{figure}[h!]
\begin{align*}
\xymatrix{
& \sqcup D_\alpha^{i+1} \times \sqcup S_\beta^{j} \ar[r] \ar[d] & \sqcup D_\alpha^{i+1} \times Y_j \ar[r] & X_{i+1} \times Y_j \ar[d] \\
\sqcup S_\alpha^{i} \times \sqcup S_\beta^{j} \ar[ru] \ar[r] \ar[rd] & \sqcup S^{j+i+1}_{(\alpha,\beta)} \ar@{.>}[rr]^{\Psi''}  & & X_iY_j^\times\\
& \sqcup S_\alpha^{i} \times \sqcup D_\beta^{j+1} \ar[r] \ar[u] & X_i \times \sqcup D_\beta^{j+1} \ar[r] & X_j \times Y_{i+1} \ar[u]
}
\end{align*}
\caption{}\label{Fig:PsiBiPrime}
\end{figure}

	The next step is to observe that the push out of $\Psi$ and $\Psi''$ is actually $X_{i+1}\times Y_{j+1}$.
	We will support this in a huge bunch of pushouts in \autoref{Fig:HugeCommDiag}.
	First we observe that the diagram
\begin{equation}
\label{Eq:PushOutFinal}
\xymatrix{
X_iY_j^\times \ar[r] & X_{i+1} \times Y_{j+1} \\
\sqcup S_{(\alpha,\beta)^{i+j+1}}\ar[r] \ar[u] & \sqcup D_\alpha^{i+1} \times \sqcup D_\beta^{j+1} \ar[u]
}
\end{equation}
is commutative. This follows by the uniqueness of the morphisms from a pushout and for the following equalities
\begin{align*}
\sqcup D_\alpha^{i+1} \times \sqcup S_\beta^{j} \rightarrow & \sqcup S_{(\alpha,\beta)}^{i+j+1} \xrightarrow{\Psi} \sqcup D_\alpha^{i+1} \times \sqcup D_\beta^{j+1} \rightarrow X_{i+1} \times \sqcup D_\beta^{j+1} \rightarrow X_{i+1} \times Y_{j+1}\\
\sqcup D_\alpha^{i+1} \times \sqcup S_\beta^{j} \rightarrow & \sqcup D_\alpha^{i+1} \times \sqcup D_\beta^{j+1} \rightarrow X_{i+1} \times \sqcup D_\beta^{j+1} \rightarrow X_{i+1} \times Y_{j+1}\\
\sqcup D_\alpha^{i+1} \times \sqcup S_\beta^{j} \rightarrow & \sqcup D_\alpha^{i+1} \times Y_j \rightarrow \sqcup D_\alpha^{i+1} \times Y_{j+1} \rightarrow X_{i+1} \times Y_{j+1}\\
\sqcup D_\alpha^{i+1} \times \sqcup S_\beta^{j} \rightarrow & \sqcup D_\alpha^{i+1} \times Y_j \rightarrow X_{i+1} \times Y_j \rightarrow X_{i+1} \times Y_{j+1}\\
\sqcup D_\alpha^{i+1} \times \sqcup S_\beta^{j} \rightarrow & \sqcup D_\alpha^{i+1} \times Y_j \rightarrow X_{i+1} \times Y_j \rightarrow X_iY_j^\times \rightarrow X_{i+1} \times Y_{j+1}\\
\sqcup D_\alpha^{i+1} \times \sqcup S_\beta^{j} \rightarrow & \sqcup S_{(\alpha,\beta)}^{i+j_1} \xrightarrow{\Psi''} X_iY_j^\times \rightarrow X_{i+1} \times Y_{j+1}
\end{align*}
and the equalities
\begin{align*}
\sqcup S_\alpha^{i} \times \sqcup D_\beta^{j+1} \rightarrow & \sqcup S_{(\alpha,\beta)}^{i+j+1} \xrightarrow{\Psi} \sqcup D_\alpha^{i+1} \times \sqcup D_\beta^{j+1} \rightarrow \sqcup D_\alpha^{i+1} \times Y_{j+1} \rightarrow X_{i+1} \times Y_{j+1}\\
\sqcup S_\alpha^{i} \times \sqcup D_\beta^{j+1} \rightarrow & \sqcup D_\alpha^{i+1} \times \sqcup D_\beta^{j+1} \rightarrow \sqcup D_\alpha^{j+1} \times Y_{j+1} \rightarrow X_{i+1} \times Y_{j+1}\\
\sqcup S_\alpha^{i} \times \sqcup D_\beta^{j+1} \rightarrow & X_i \times \sqcup D_\beta^{j+1} \rightarrow X_{i+1} \sqcup D_\beta^{i+1} \rightarrow X_{i+1} \times Y_{j+1}\\
\sqcup S_\alpha^{i} \times \sqcup D_\beta^{j+1} \rightarrow & X_i \times \sqcup D_\beta^{i+1} \rightarrow X_{i} \times Y_{j+1} \rightarrow X_{i+1} \times Y_{j+1}\\
\sqcup S_\alpha^{i} \times \sqcup D_\beta^{j+1} \rightarrow & X_i \times \sqcup D_\beta^{i+1} \rightarrow X_{i} \times Y_{j+1} \rightarrow X_iY_j^\times \rightarrow X_{i+1} \times Y_{j+1}\\
\sqcup S_\alpha^{i} \times \sqcup D_\beta^{j+1} \rightarrow & \sqcup S_{(\alpha,\beta)}^{i+j_1} \xrightarrow{\Psi''} X_iY_j^\times \rightarrow X_{i+1} \times Y_{j+1}
\end{align*}

The following step is to show that the commutative square \ref{Eq:PushOutFinal} is a pushout. Consider $W$ a pseudotopological space and two continuous maps $F:X_iY_j^\times \rightarrow W$ and $G: \sqcup D_\alpha^{i+1}\times \sqcup D_\beta^{j+1}\rightarrow W$ such that $F\Psi''=G\Psi$. In the subsequent diagrams, we chase a map from $X_{i+1}\times Y_{j+1}$ to $W$ through some of the pushouts from \autoref{Fig:HugeCommDiag}.

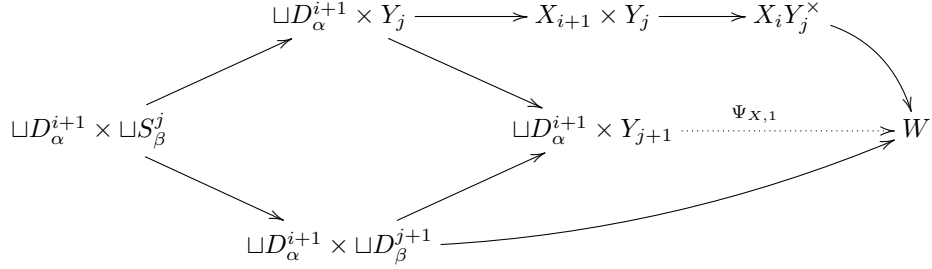
\begin{figure}[h!]
\begin{align*}
\xymatrix{
{} & \sqcup D_\alpha^{i+1} \times Y_j \ar[r] \ar[rd] & X_{i+1}\times Y_j \ar[r] & X_iY_j^\times \ar@/^1pc/[rd] & {}\\
\sqcup D_\alpha^{i+1}\times \sqcup S_\beta^{j} \ar[ru]\ar[rd] & {} & \sqcup D_\alpha^{i+1}\times Y_{j+1} \ar@{.>}[rr]^{\Psi_{X,1}} & {} & W\\
{} & \sqcup D_\alpha^{i+1} \times \sqcup D_\beta^{j+1} \ar[ru] \ar@/_1pc/[rrru] & {} & {} & {}
}
\end{align*}
\caption{There exists a unique $\Psi_{X,1}$ that makes the diagram commute}\label{Fig:DY1}
\end{figure}

\begin{figure}[h!]
\begin{align*}
\xymatrix{
{} & \sqcup S_\alpha^i \times Y_j \ar[rd] \ar[r] & \sqcup D_\alpha^{i+1} \times Y_j \ar[r] & X_{i+1}\times Y_j \ar[r] & X_iY_j^\times \ar[d]\\
\sqcup S_\alpha^i \times \sqcup S_\beta^j \ar[ru]\ar[rd] & {} & \sqcup S_\alpha^i \times Y_{j+1} \ar@{.>}[rr]^{\Psi_{X,2}} & {} & W\\
{} & \sqcup S_\alpha^i \times \sqcup D_\beta^{j+1} \ar[ru]\ar[r] & \sqcup D_\alpha^{i+1} \times \sqcup D_\beta^{j+1} \ar@/_1pc/[rru] & {} & {}
}
\end{align*}
\caption{There exists a unique $\Psi_{X,2}$ that makes the diagram commute}\label{Fig:SY1}
\end{figure}
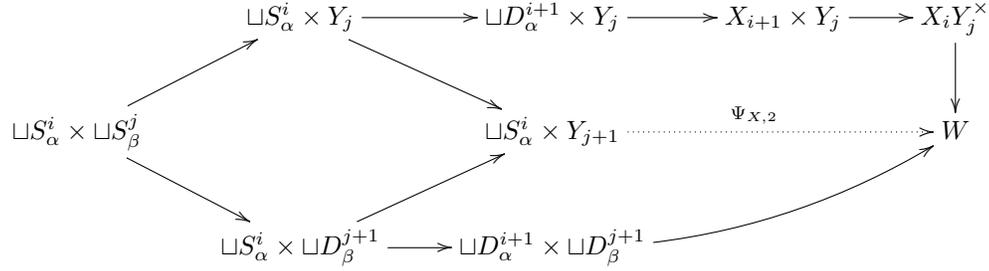

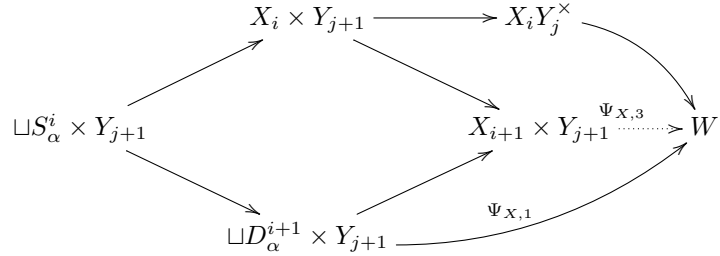
\begin{figure}[h!]
\begin{align*}
\xymatrix{
{} & X_i \times Y_{j+1} \ar[rd] \ar[r] & X_iY_j^\times \ar@/^1pc/[rd] & {}\\
\sqcup S_\alpha^i \times Y_{j+1} \ar[ru]\ar[rd] & {} & X_{i+1}\times Y_{j+1} \ar@{.>}[r]^{\Psi_{X,3}} & W\\
{} & \sqcup D_\alpha^{i+1} \times Y_{j+1} \ar@/_1.5pc/[rru]^{\Psi_{X,1}} \ar[ru] & {} & {}
}
\end{align*}
\caption{There exists a unique $\Psi_{X,3}$ that makes the diagram commute}\label{Fig:X1Y1fromX0Y1}
\end{figure}

	Giving the \autoref{Fig:DY1} and the uniqueness of $\Psi_{X,2}$ in \autoref{Fig:SY1}, we can note that the following exterior square commutes and we can use the pushout property to obtain \autoref{Fig:X1Y1fromX0Y1}.

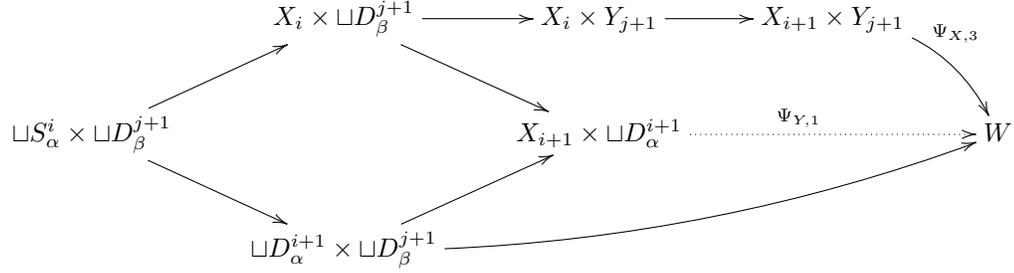
\begin{figure}[h!]
\begin{align*}
\xymatrix{
{} & X_i\times \sqcup D_\beta^{j+1} \ar[r] \ar[rd] & X_{i}\times Y_{j+1} \ar[r] & X_{i+1}\times Y_{j+1} \ar@/^1pc/[rd]^{\Psi_{X,3}} & {}\\
\sqcup S_\alpha^{i}\times \sqcup D_\beta^{j+1} \ar[ru]\ar[rd] & {} & X_{i+1} \times \sqcup D_\alpha^{i+1} \ar@{.>}[rr]^{\Psi_{Y,1}} & {} & W\\
{} & \sqcup D_\alpha^{i+1} \times \sqcup D_\beta^{j+1} \ar[ru] \ar@/_1pc/[rrru] & {} & {} & {}
}
\end{align*}
\caption{There exists a unique $\Psi_{Y,1}$ that makes this diagram commutes}\label{Fig:X1D}
\end{figure}

On the other hand, we can use the morphism $\Psi_{X,3}$ to build up the following commutative diagrams. Observing that $(\sqcup S_\alpha^{i} \times \sqcup D_\beta^{j+1} \rightarrow X_i \times \sqcup D_\beta^{j+1}\rightarrow X_i \times Y_{j+1})=(\sqcup S_\alpha^{i} \times \sqcup D_\beta^{j+1} \rightarrow \sqcup S_\alpha^{i} \times Y_{j+1}\rightarrow X_i \times Y_{j+1})$ and $(X_i \times Y_{j+1} \rightarrow X_{i+1} \times Y_{j+1} \xrightarrow{\Psi_{X,3}} W)=(X_i \times Y_{j+1} \rightarrow X_iY_j^\times \rightarrow W)$ we have the \autoref{Fig:X1D}.

\begin{figure}[h!]
\begin{align}
\xymatrix{
{} & X_i \times \sqcup S_\beta^j \ar[rd] \ar[r] & X_i \times \sqcup D_\beta^{j+1} \ar[r] & X_{i}\times Y_{j+1} \ar[r] & X_{i}\times Y_{j+1} \ar[d]_{\Psi_{X,3}}\\
\sqcup S_\alpha^i \times \sqcup S_\beta^j \ar[ru]\ar[rd] & {} & X_{i+1} \times \sqcup S_\alpha^i \ar@{.>}[rr]^{\Psi_{Y,2}} & {} & W\\
{} & \sqcup D_\alpha^{i+1} \times \sqcup S_\beta^{j} \ar[ru]\ar[r] & \sqcup D_\alpha^{i+1} \times \sqcup D_\beta^{j+1} \ar@/_1pc/[rru] & {} & {}
}
\end{align}
\caption{There exists a unique $\Psi_{Y,2}$ that makes this diagram commutes}\label{Fig:X1S}
\end{figure}
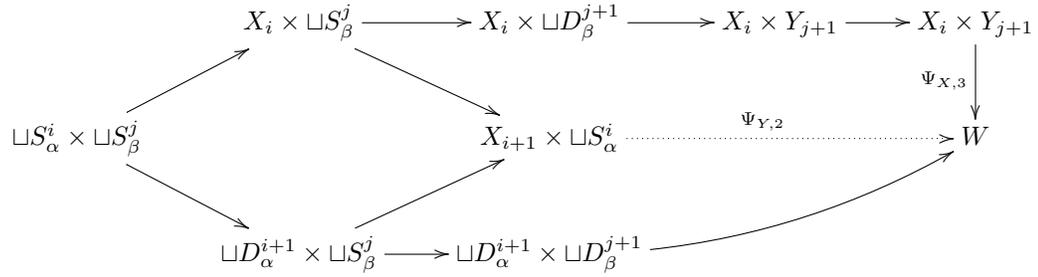

Now, recalling that $(\sqcup S_\alpha^{i} \times \sqcup S_\beta^{j} \rightarrow X_i \times \sqcup S_\beta^{j}\rightarrow X_i \times \sqcup D_\beta^{j+1})=(\sqcup S_\alpha^{i} \times \sqcup S_\beta^{j} \rightarrow \sqcup S_\alpha^{i} \times \sqcup D_\beta^{j+1}\rightarrow X_i \times \sqcup D_\beta^{j+1})$, we have the \autoref{Fig:X1S}.

\begin{figure}[h!]
\begin{align*}
\xymatrix{
{} & X_{i+1} \times Y_{j} \ar[rd] \ar[r] & X_iY_j^\times \ar@/^1pc/[rd] & {}\\
X_{i+1} \times \sqcup S_\beta^j \ar[ru]\ar[rd] & {} & X_{i+1}\times Y_{j+1} \ar@{.>}[r]^{\Psi_{Y,3}} & W\\
{} & X_{i+1} \times \sqcup D_\beta^{j+1} \ar@/_1.5pc/[rru]^{\Psi_{X,1}} \ar[ru] & {} & {}
}
\end{align*}
\caption{There exists a unique $\Psi_{Y,3}$ that makes this diagram commutes}\label{Fig:X1Y1fromX1Y0}
\end{figure}
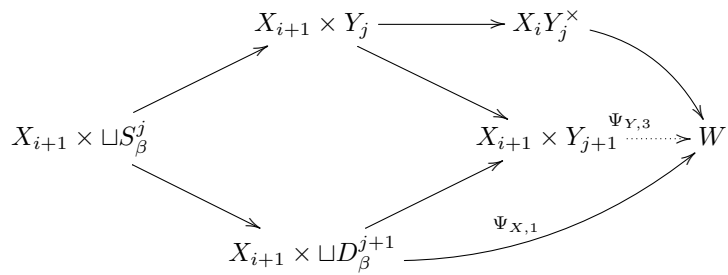

In an analogous way to \autoref{Fig:DY1}, \ref{Fig:SY1} and \ref{Fig:X1Y1fromX0Y1}, we obtain the exterior square of \autoref{Fig:X1Y1fromX1Y0} commutes thanks to \autoref{Fig:X1D} and \autoref{Fig:X1S}.

With these all of diagrams we can observe that $(X_{i+1} \times \sqcup D_\beta^{j+1} \rightarrow X^{i+1} \times Y^{j+1} \xrightarrow{\Psi_{X,3}} W)=(X_{i+1} \times \sqcup D_\beta^{j+1} \rightarrow X^{i+1} \times Y^{j+1} \xrightarrow{\Psi_{Y,3}} W)$ and $(\sqcup D_\alpha^{i+1} \times Y_{j+1} \rightarrow X^{i+1} \times Y^{j+1} \xrightarrow{\Psi_{X,3}} W)=(\sqcup D_\alpha^{i+1} \times Y^{j+1} \rightarrow X^{i+1} \times Y^{j+1} \xrightarrow{\Psi_{Y,3}} W)$. Thus $\Psi_{X,3}=\Psi_{Y,3}$ because we have embeddings onto $X^{i+1} \times Y^{j+1}$. Thus, we have a pushout.

To conclude this appendix, we are ready to observe that the sequence of $B_i$ actually forms a CW complex and, thus, $X\times Y$ has a structure as CW complex.

First, if we define $X_{-1}=Y_{-1}=\varnothing$, we observe that $B_i$ is homemorphic to
\begin{align*}
B'_j \coloneqq B'_{j,0} \bigsqcup\limits_{X_{j}\times Y_{0}} B'_{j,1} \bigsqcup\limits_{X_{j-1}\times Y_{1}} \ldots \bigsqcup\limits_{X_{1}\times Y_{j-1}}  B'_{j,j} \bigsqcup\limits_{X_{0}\times Y_{j}} B'_{j,j+1}
\end{align*}
for every $j\in\mathbb{N}_0$, where we define in addition
\begin{align*}
B'_{j,i} \coloneqq (X_{j-1+1} \times Y_{i-1}) \bigsqcup\limits_{X_{j-i}\times Y_{i-1}} (X_{j-i}\times Y_{i})
\end{align*}
for every $i\in \{0,1,\ldots,j+1\}$. Those constructions of $B'_j$, in addition to the following result, are the ingredients to obtain our desired result.

\begin{proposicion}
\label{prop:PushoutOfPushout}
Let $Z_1,Z_2,\ldots,Z_n$ be a collection of pseudotopological spaces and, for some $k\in \{1,\ldots,n\}$, and let the following commutative square
\begin{align}\label{diag:PushoutOfPushout1}
\xymatrix{
A \ar[r]^{f_k} \ar[d]_g & Z_k \ar[d]^{g'_k}\\
B \ar[r]^{f'_k} & W_k 
}
\end{align}
be a pushout with $g$ an embedding. Then the following commutative square
\begin{align}\label{diag:PushoutOfPushout2}
\xymatrix{
A \ar[r]^{f_k}\ar[d]_{g} & Z_1 \sqcup \ldots \sqcup Z_n \ar[d]^{g'_1\sqcup\ldots \sqcup g'_n}\\
B \ar[r]_{f'_k} & W_1 \sqcup \ldots \sqcup W_n
}
\end{align}
is a pushout, where $W_i=Zi$ and $g'_i=\Id{Z_i}$ if $k\neq i$. In addition, if we have a pair of equivalence relations $\sim_Z$ and $\sim_W$ in $Z_1 \sqcup \ldots \sqcup Z_n$ and $W_1 \sqcup \ldots \sqcup W_n$, respectively, such that the relation of $W_1 \sqcup \ldots \sqcup W_n$ are exactly the same as the relations of $Z_1 \sqcup \ldots \sqcup Z_n$ seen as a subspace. Then, the exterior square of the following commutative diagram
\begin{align}\label{diag:PushoutOfPushout3}
\xymatrix{
A \ar[r]\ar[d] & Z_1 \sqcup \ldots \sqcup Z_n \ar[d] \ar[r]^{q_z} & (Z_1 \sqcup \ldots \sqcup Z_n)/\sim_Z \ar[d]^{(g'_1\sqcup\ldots \sqcup g'_n)'}\\
B \ar[r] & W_1 \sqcup \ldots \sqcup W_n \ar[r]^{q_w} & (W_1 \sqcup \ldots \sqcup W_n)/\sim_W}
\end{align}
is a pushout.
\end{proposicion}

\begin{proof}
	Consider the pushout of the diagram \ref{diag:PushoutOfPushout1} described in the sentence of the theorem.
	We are going to verify that the diagram \ref{diag:PushoutOfPushout2} is a pushout.
	Let $C$ be a pseudotopological space and $h:B\rightarrow C$ and $h_i:Z_i\rightarrow C$ such that \[hg=(h_1 \sqcup \ldots \sqcup h_n)f_k,\] then we have that there exists a unique $\alpha_k:W_k\rightarrow C$ such that $\alpha_k f'_k = h$ and $\alpha_k g'_k = h_k$.
	If we define $\alpha_i\coloneqq h_i$ for every $i\neq k$, we obtain that \[(\alpha_1 \sqcup \ldots \sqcup \alpha_n)f'_k = h,\ (\alpha_1 \sqcup \ldots \sqcup \alpha_n)(g'_1 \sqcup \ldots \sqcup g'_n)=(h_1 \sqcup \ldots \sqcup h_n). \]
	It only rests to prove that this continuous maps is the only one which makes both triangles commutative.
	Observe that none of the $\alpha_i$ with $i\neq k$ can be changed by construction; if we ask for $\alpha_i\neq h_i$ then $\alpha_i g'_i = \alpha_i \neq h_i$.
	By the uniqueness condition on the pushout, $h_k$ neither can be changed.
	Thus the diagram \ref{diag:PushoutOfPushout2} is a pushout.
	
	For the second part of the theorem, we use that $g$ is an embedding, and thus $g'_1\sqcup \ldots \sqcup g'_n$ is also an embedding by \autoref{prop:EmbPushoutPsTop}, and we ask for two relations $\sim_W$ and $\sim_Z$ such that $w_1 \sqcup \ldots \sqcup w_n \sim_W w'_1 \sqcup \ldots \sqcup w'_n$ if and only if $w_k,w'_k \in g'_k(Z_k)$ and $z_1 \sqcup \ldots \sqcup z_n \sim_W z'_1 \sqcup \ldots \sqcup z'_n$ with $z_i=w_i$ and $z'_i=w'_i$ if $i\neq k$ and $z_k=(g'_k)^{-1}(w_k)$ and $z'_k=(g'_k)^{-1}(w'_k)$.
	
	Using the notation of the diagram \ref{diag:PushoutOfPushout3}, considering a pseudotopological space $C$ and maps $h:B\rightarrow C$ and $h':(Z_1 \sqcup \ldots \sqcup \ldots Z_n)/\sim_Z \rightarrow C$ such that $h' q_z f_k = h g$, since the diagram \ref{diag:PushoutOfPushout2} is a pushout, we have that there exists a unique $\beta_1:W_1 \sqcup \ldots \sqcup W_n \rightarrow C$ such that $\beta_1 f'_k = h$ and $\beta_1 (g'_1 \sqcup \ldots \sqcup g'_n) = h' q_z$.
	By the way we define the relations and given $\beta_1 (g'_1 \sqcup \ldots \sqcup g'_n) = h' q_z$, by the universal property of the quotient we have that there exists a unique $\beta_2: (W_1 \sqcup \ldots \sqcup W_n)/\sim_W \rightarrow C$ such that $\beta_2 q_w = \beta_1$.
	
	Lastly, since $h' q_z = \beta_2 (g'_1 \sqcup \ldots \sqcup g'_n)' q_z$ and $q_z$ is surjective, then $h' = \beta_2 (g'_1 \sqcup \ldots \sqcup g'_n)'$. We conclude that the diagram \ref{diag:PushoutOfPushout3} is a pushout.
\end{proof}

With these two constructions, we can observe inductively over $(j,i)$ that the sequence of $B'_j$ is actually a CW complex give \autoref{def:PsTopCWComplex}.

\printnomenclature
\bibliographystyle{amsalpha}
\bibliography{all}

\begin{bibdiv}
\begin{biblist}

\bib{Babson_etal_2006}{article}{
      author={Babson, Eric},
      author={Barcelo, H{\'e}l{\`e}ne},
      author={de~Longueville, Mark},
      author={Laubenbacher, Reinhard},
       title={Homotopy theory of graphs},
        date={2006},
        ISSN={0925-9899},
     journal={J. Algebraic Combin.},
      volume={24},
      number={1},
       pages={31\ndash 44},
         url={http://dx.doi.org/10.1007/s10801-006-9100-0},
      review={\MR{2245779}},
}

\bib{Beattie_Butzmann_2002}{book}{
      author={Beattie, R.},
      author={Butzmann, H.-P.},
       title={Convergence structures and applications to functional analysis},
   publisher={Kluwer Academic Publishers, Dordrecht},
        date={2002},
        ISBN={1-4020-0566-0},
         url={http://dx.doi.org/10.1007/978-94-015-9942-9},
      review={\MR{2327514}},
}

\bib{bubenik2024_RelCellCompClos}{misc}{
      author={Bubenik, Peter},
       title={Relative cell complexes in closure spaces},
        date={2024},
         url={https://arxiv.org/abs/2312.13227},
}

\bib{Bubenik_Milicevic_2021}{misc}{
      author={Bubenik, Peter},
      author={Milićević, Nikola},
       title={Eilenberg-steenrod homology and cohomology theories for Čech's
  closure spaces},
   publisher={arXiv},
        date={2021},
         url={https://arxiv.org/abs/2112.13421},
}

\bib{Carlsson_2009}{article}{
      author={Carlsson, Gunnar},
       title={Topology and data},
        date={2009},
        ISSN={0273-0979},
     journal={Bull. Amer. Math. Soc. (N.S.)},
      volume={46},
      number={2},
       pages={255\ndash 308},
         url={http://dx.doi.org/10.1090/S0273-0979-09-01249-X},
      review={\MR{2476414}},
}

\bib{Carlsson_et_al_2012}{article}{
      author={Carlsson, Gunnar},
      author={de~Silva, Vin},
      author={Ishkhanov, Tigran},
      author={Zomorodian, Afra},
       title={{On the local behavior of spaces of natural images}},
        date={2008},
     journal={International Journal of Computer Vision},
      volume={76},
       pages={1\ndash 12},
}

\bib{Chih_Scull_2021}{article}{
      author={Chih, Tien},
      author={Scull, Laura},
       title={A homotopy category for graphs},
        date={2021},
        ISSN={0925-9899},
     journal={J. Algebraic Combin.},
      volume={53},
      number={4},
       pages={1231\ndash 1251},
         url={https://doi.org/10.1007/s10801-020-00960-5},
      review={\MR{4263649}},
}

\bib{Demaria_1987}{inproceedings}{
      author={Demaria, Davide~Carlo},
       title={On some applications of homotopy to directed graphs},
        date={1987},
   booktitle={Proceedings of the {G}eometry {C}onference ({M}ilan and
  {G}argnano, 1987)},
      volume={57},
       pages={183\ndash 202 (1989)},
         url={http://dx.doi.org/10.1007/BF02925050},
      review={\MR{1017928}},
}

\bib{Dochtermann_2009}{article}{
      author={Dochtermann, Anton},
       title={Hom complexes and homotopy theory in the category of graphs},
        date={2009},
        ISSN={0195-6698},
     journal={European J. Combin.},
      volume={30},
      number={2},
       pages={490\ndash 509},
         url={https://doi.org/10.1016/j.ejc.2008.04.009},
      review={\MR{2489282}},
}

\bib{kapulkin_edel2023_SyntheticApproach}{misc}{
      author={Ebel, Sterling},
      author={Kapulkin, Chris},
       title={Synthetic approach to the quillen model structure on topological
  spaces},
        date={2023},
}

\bib{FelixHalperinThomas_2001_RatHom}{book}{
      author={F\'{e}lix, Yves},
      author={Halperin, Stephen},
      author={Thomas, Jean-Claude},
       title={Rational homotopy theory},
      series={Graduate Texts in Mathematics},
   publisher={Springer-Verlag, New York},
        date={2001},
      volume={205},
        ISBN={0-387-95068-0},
  url={https://doi-org.biblio.cimat.remotexs.co/10.1007/978-1-4613-0105-9},
      review={\MR{1802847}},
}

\bib{Hatcher_2002}{book}{
      author={Hatcher, Allen},
       title={Algebraic topology},
   publisher={Cambridge University Press, Cambridge},
        date={2002},
        ISBN={0-521-79540-0},
      review={\MR{1867354}},
}

\bib{hess2007rationalhomotopytheorybrief}{incollection}{
      author={Hess, Kathryn},
       title={Rational homotopy theory: a brief introduction},
        date={2007},
   booktitle={Interactions between homotopy theory and algebra},
      series={Contemp. Math.},
      volume={436},
   publisher={Amer. Math. Soc., Providence, RI},
       pages={175\ndash 202},
         url={https://doi.org/10.1090/conm/436/08409},
      review={\MR{2355774}},
}

\bib{Heuts_2022_SimplicialAndDendroidal}{book}{
      author={Heuts, Gijs},
      author={Moerdijk, Ieke},
       title={Simplicial and dendroidal homotopy theory},
      series={Ergebnisse der Mathematik und ihrer Grenzgebiete. 3. Folge. A
  Series of Modern Surveys in Mathematics [Results in Mathematics and Related
  Areas. 3rd Series. A Series of Modern Surveys in Mathematics]},
   publisher={Springer, Cham},
        date={[2022] \copyright 2022},
      volume={75},
        ISBN={978-3-031-10446-6; 978-3-031-10447-3},
  url={https://doi-org.biblio.cimat.remotexs.co/10.1007/978-3-031-10447-3},
      review={\MR{4485749}},
}

\bib{Hovey_1999}{book}{
      author={Hovey, Mark},
       title={Model categories},
      series={Mathematical Surveys and Monographs},
   publisher={American Mathematical Society, Providence, RI},
        date={1999},
      volume={63},
        ISBN={0-8218-1359-5},
      review={\MR{1650134}},
}

\bib{Lupton_etal_2022a}{article}{
      author={Lupton, Gregory},
      author={Oprea, John},
      author={Scoville, Nicholas~A.},
       title={Homotopy theory in digital topology},
        date={2022},
        ISSN={0179-5376,1432-0444},
     journal={Discrete Comput. Geom.},
      volume={67},
      number={1},
       pages={112\ndash 165},
  url={https://doi-org.biblio.cimat.remotexs.co/10.1007/s00454-021-00336-x},
      review={\MR{4356246}},
}

\bib{Lupton_Scoville_2022}{article}{
      author={Lupton, Gregory},
      author={Scoville, Nicholas~A.},
       title={Digital fundamental groups and edge groups of clique complexes},
        date={2022},
        ISSN={2367-1726,2367-1734},
     journal={J. Appl. Comput. Topol.},
      volume={6},
      number={4},
       pages={529\ndash 558},
  url={https://doi-org.biblio.cimat.remotexs.co/10.1007/s41468-022-00095-5},
      review={\MR{4496690}},
}

\bib{Preuss_2002}{book}{
      author={Preuss, Gerhard},
       title={Foundations of topology},
   publisher={Kluwer Academic Publishers, Dordrecht},
        date={2002},
        ISBN={1-4020-0891-0},
         url={https://doi.org/10.1007/978-94-010-0489-3},
        note={An approach to convenient topology},
      review={\MR{2033142}},
}

\bib{Riehl_2016}{book}{
      author={Riehl, Emily},
       title={Category theory in context},
   publisher={Dover Publications, Inc},
        date={2016},
        ISBN={978-0-486-80903-8},
}

\bib{Rieser_2021}{article}{
      author={Rieser, Antonio},
       title={\v{C}ech closure spaces: {A} unified framework for discrete and
  continuous homotopy},
        date={2021},
        ISSN={0166-8641},
     journal={Topology and its Applications},
      volume={296},
       pages={107613},
  url={https://www.sciencedirect.com/science/article/pii/S0166864121000249},
}

\bib{Rieser_arXiv_2022}{misc}{
      author={Rieser, Antonio},
       title={Cofibration and model category structures for discrete and
  continuous homotopy},
   publisher={arXiv},
        date={2022},
         url={https://arxiv.org/abs/2209.13510},
}

\bib{rieser2023semicoarse}{misc}{
      author={Rieser, Antonio},
      author={Treviño-Marroquín, Jonathan},
       title={Semi-coarse spaces, homotopy and homology},
        date={2025},
      volume={167},
  url={https://www.sciencedirect.com/science/article/pii/S0196885825000326},
}

\bib{Robert_Thesis2013}{thesis}{
      author={Robert-Nicoud, Daniel},
       title={Quillen equivalence of topological spaces and simplicial sets},
        type={Ph.D. Thesis},
        date={2013},
}

\bib{tom_Dieck_2008}{book}{
      author={tom Dieck, Tammo},
       title={Algebraic topology},
      series={EMS Textbooks in Mathematics},
   publisher={European Mathematical Society (EMS), Z\"urich},
        date={2008},
        ISBN={978-3-03719-048-7},
         url={http://dx.doi.org/10.4171/048},
      review={\MR{2456045}},
}

\bib{Trevino_2024_GraphsVRPsTop}{misc}{
      author={Treviño-Marroquín, Jonathan},
       title={Graphs and their vietoris-rips complexes have the same
  pseudotopological weak homotopy type},
        date={2024},
         url={https://arxiv.org/abs/2406.00149},
}

\bib{Trevino_2025_FundGroupoidSC}{misc}{
      author={Treviño-Marroquín, Jonathan},
       title={Semi-coarse spaces: Fundamental groupoid and the van kampen
  theorem},
        date={2025},
         url={https://arxiv.org/abs/2404.08874},
}

\end{biblist}
\end{bibdiv}

%\addcontentsline{toc}{chapter}{Bibliografía}

\end{document}